\documentclass[12pt,a4paper]{amsart}%
\usepackage{amsfonts}
\usepackage{amsmath}
\usepackage{amssymb}
\usepackage{graphicx}%

\usepackage{tikz}

\usepackage[margin=20mm]{geometry}

\usepackage{multirow}
\usepackage{array}

\usepackage[colorlinks, linkcolor=blue,  citecolor=blue, urlcolor=blue, bookmarks=false]{hyperref}%
\hypersetup{pdfstartview=FitH}
\usepackage{url}
\usepackage{xspace}

\usepackage{dsfont}

\usepackage{verbatim}
\usepackage{alltt}
\usepackage{graphicx}
\graphicspath{{../figures/}{../figuresOld/}}

\newtheorem{theorem}{Theorem}[section]

\newtheorem{lemma}[theorem]{Lemma}

\newtheorem{proposition}[theorem]{Proposition}
\newtheorem{remark}[theorem]{Remark}

\numberwithin{equation}{section}

\def\XXint#1#2#3{{\setbox0=\hbox{$#1{#2#3}{\int}$}
     \vcenter{\hbox{$#2#3$}}\kern-.5\wd0}}

\let\originalleft\left
\let\originalright\right
\renewcommand{\left}{\mathopen{}\mathclose\bgroup\originalleft}
\renewcommand{\right}{\aftergroup\egroup\originalright}

\renewcommand{\ref}[1]{\autoref{#1}}

\AtBeginDocument{

}

\usepackage{enumitem}

\setlist[description]{align=right, labelwidth=2cm, leftmargin=!, font=\normalfont\bfseries}

\begin{document}
\title[Simultaneous coordinate-wise CLT]{Simultaneous coordinate-wise Wasserstein central limit theorem in high dimensions}
\author{Yuanlong Ruan}
\thanks{Copyright \copyright~2026 The Author. All rights reserved.}
\address{Beihang University, Beijing, China}

\subjclass[2020]{60F05, 49Q22, 65C99}
\keywords{Multivariate central limit theorem, stochastic localization, fourth moment, Wasserstein metric, Hermite expansion, Bessel's inequality}

\date{\today}

\begin{abstract}

For independent centered random vectors $(X_i)_{i=1}^n$ of finite fourth moments, we obtain dimension free Wasserstein bound between its normalized sum $W_n =  n^{-1/2}\sum X_i$ and a Gaussian whose covariance equals that of $W_n$. The coupling between $W_n$ and the Gaussian is strong in the sense that a single coupling controls controls every projection onto a subset of those eigen coordinates. The proof combines stochastic localization, Hermite expansion on a Gaussian prior with exchangeable pairing. This proves the Conjecture surmised in \cite{fang2026sharp} and improves the coefficient on the fourth-moment term from $2$ to $197/108$. 

The same Wasserstein error bound extends to real separable Hilbert space. For Hilbert-valued i.i.d. observations, the $L^2$ prediction error bound has an order of $O(n^{-1/2})$, with no dimension factor or restriction on covariance eigenvalue decay, and the same coupling retains bounds for every spectral projection. In an $\ell^2$ example, we demonstrate that the coordinate-wise error gives sharp error order of $n^{-1/2}$.

\end{abstract}

\maketitle
\setcounter{tocdepth}{1}
\tableofcontents

\section{Introduction and main results}

For centered i.i.d $X_{i},i=1,...,n,$ write $W_{n}=n^{-1/2}\sum_{i=1}^{n}
X_{i}$ and $\Sigma=\operatorname{Cov}\left(  W_{n}\right)  .$In a seminal
work, Berry-Esseen \cite{esseen1956moment} proves that in one dimension
$W_{n}$ converges to standard Gaussian distribution weakly with the optimal
rate $n^{-1/2},$ the error is quanitfied. It leaves open whether such rate is
attainable in high dimension \cite{bentkus2003dependence,bentkus2005lyapunov}.
In the effort of extending Berry-Esseen bounds to high dimensions, quantiative
central limit theorems are developed rapidly. In one dimension, a quantiative
rate $n^{-1/2}$ in quadratic Wasserstein distance is given for normalized
i.i.d. with finite fourth moments \cite[Theorem~4.1, $r=2$]{rio2009upper}. In
high dimension, \ under the assumption that all $X_{i}$ are bounded by a
common constant, \cite{zhai2018high,eldan2020clt} proved multivariate central
limit theorems that quantitatively bound the error $W_{2}^{2}(\mathcal{L}
(W_{n}),N(0,\Sigma)).$ Both bounds contain a factor involving $1+\log n.$
\cite{eldan2020clt} uses martingale embedding to make coupling between the
random vector and a Brownian motion with inspiration from
\cite{eldan2016skorokhod}. For centered i.i.d.\ vectors of covariance $I_{d}$,
\cite[Theorem~1]{bonis2020stein} obtained a bound for $W_{2}(\mathcal{L}
(W_{n}),\gamma_{d})$ by $O\left(  n^{-1/2}\right)  $ under finite fourth
moments assumption. The logarithm factor in earlier work is removed. Bonis
smooths the law along the Ornstein-Uhlenbeck flow, and uses exchangeability to
expand a Gaussian score in Hermite polynomials. In high dimensions, several
other works are available, e.g. \cite{fang2021high,lopes2022central}.

The recent \cite[Theorem~2.1]{fang2026sharp} proves the aggregate estimate
\begin{equation}
W_{2}^{2}(\mathcal{L}(W_{n}),N(0,\Sigma))\leqslant2\sum_{j:\,\sigma_{j}>0}
\min\left\{  \sigma_{j}^{2},O\left(  \frac{1}{n^{2}}\right)  \right\}  ,
\label{eq:aggregate-CLT}
\end{equation}
where $\sigma_{j}^{2}$ are eigenvalues of $\Sigma.$ They use the martingale
property of the drift of a F\"{o}llmer process to bound the transport cost by
an integral of a squared conditional score. Then combine this with resampling
and Bessel's inequality.

In several dimensions, the bound $\left(  \ref{eq:aggregate-CLT}\right)  $
improves several early results including \cite[Theorem~1]{bonis2020stein}
\cite{eldan2020clt}, but separately coupling each coordinate to a normal
variable leaves open whether the normal variables have the required joint
Gaussian distribution. In other words, we asks for the corresponding bounds on
the individual errors within one coupling. This is the \textbf{conjecture} in
\cite[Theorem~2.1]{fang2026sharp}.

The current article proves the conjecture in Theorem \ref{thm:main}, thus
yielding a stronger multivariate central limit theorem. The proof uses
stochastic localization. To obtain a tractable bound we couple the original
process with one that starts with a Gaussian in Sections~\ref{sec:observation}
-\ref{sec:separation}. Exchangeability \cite[Theorem~1]{bonis2020stein} allows
us to compute the bound via Hermite expansion. The stochastic localization is
initiated in \cite{eldan2013thin}, see also \cite{el2022information} for a review.

We also deal with independent random objects taking values in real separable
Hilbert space. The coordinate-wise central limit theorem is established and in
an i.i.d. Hilbert space example we demonstrate the strength of the
coordinate-wise bound with a sharp order of $\left(  n^{-1/2}\right)  .$

\subsection{notations}

$\langle\cdot,\cdot\rangle$ and $|\cdot|$ are the Euclidean inner product and
norm. For matrices, $\Vert A\Vert^{2}=\sum_{j,\ell}A_{j,\ell}^{2}$ is the
Frobenius Norm, $\operatorname{tr}A=\sum_{j}A_{j,j}$ is the trace, $A^{\top}$
is the transpose of $A$, $I_{d}$ is the identity matrix.

In a chosen orthonormal basis $(e_{j})_{j=1}^{d}$ of $\mathbb{R}^{d}$,
$v_{j}=\langle v,e_{j}\rangle$ for a vector $v$ and $A_{j,\ell}=\langle
e_{j},Ae_{\ell}\rangle$ for a matrix $A$. The chosen eigenbasis of $\Sigma$
that diagonalizes the covariance is always specified in Theorem~\ref{thm:main}
and Section~\ref{sec:separation}. All vectors and matrices in a calculation
use that same basis.

For a subset $J$, write $\mathbb{R}^{J}$ for real vectors indexed by $J$,
write $x_{J}=(x_{j})_{j\in J}$. We write $1_{E}$ for the indicator of an event
$E$. Write $\mathcal{L}(V)$ for the distributioin of $V$. $N(0,L)$ denotes the
centered Gaussian distribution with covariance $L$, including singular $L$.
Write $\Pi(P,Q)$ for the couplings with marginals $P,Q$ of finite second
moments. Write $\mathbb{E}$ for expectation and $\Vert V\Vert_{L^{p}
}=(\mathbb{E}|V|^{p})^{1/p}$ for $1\leqslant p<\infty$.

\subsection{main result}

Throughout the article,
\begin{equation}
C_{0}=\frac{197}{208}. \label{eq:C0}
\end{equation}

\begin{theorem}
\label{thm:main}Let $n\geqslant1$ be an integer and $X_{1},\ldots,X_{n}$ be
independent random vectors in $\mathbb{R}^{d}$ with $\mathbb{E}X_{i}=0,$
$\mathbb{E}\left\vert X_{i}\right\vert ^{2}<\infty,$ $i=1,...,n$. Let $\left(
X_{i}^{\prime}\right)  _{i=1}^{n}$ be an independent copy of $\left(
X_{i}\right)  _{i=1}^{n}$. Define
\[
W_{n}=\frac{1}{\sqrt{n}}\sum_{i=1}^{n}X_{i}.
\]
Let $\Delta_{i}=X_{i}^{\prime}-X_{i}$ and $\Sigma=\operatorname{Cov}\left(
W_{n}\right)  .$ Fix an orthonormal eigenbasis $\left(  e_{j}\right)
_{j=1}^{d}$ of $\Sigma$ and write $\Sigma e_{j}=\sigma_{j}^{2}e_{j}$ with
$\sigma_{j}\geqslant0$. In particular $\Delta_{i,j}=\left\langle X_{i}
^{\prime}-X_{i},e_{j}\right\rangle $ and $W_{n,j}=\left\langle W_{n}
,e_{j}\right\rangle $. For $\sigma_{j}>0$, define
\begin{equation}
b_{j}=\frac{1}{n^{2}\sigma_{j}^{2}}\sum_{i=1}^{n}\mathbb{E}\left(  \left\vert
\Delta_{i}\right\vert ^{2}\Delta_{i,j}^{2}\right)  \in\left[  0,\infty\right]
. \label{eq:budget}
\end{equation}
There exists a coupling between $W_{n}$ and $G\sim N\left(  0,\Sigma\right)  $
such that, simultaneously for every $j$ with $\sigma_{j}>0$,
\begin{equation}
W_{2}^{2}\left[  \mathcal{L}\left(  W_{n,j}\right)  ,N\left(  0,\sigma_{j}
^{2}\right)  \right]  \leqslant\mathbb{E}\left[  \left(  W_{n,j}-G_{j}\right)
^{2}\right]  \leqslant\min\left\{  2\sigma_{j}^{2},C_{0}b_{j}\right\}  .
\label{eq:main}
\end{equation}
For $\sigma_{j}=0$ both coupled coordinates vanish almost surely.
\end{theorem}

\begin{remark}
If $\mathbb{E}\left\vert X_{i}\right\vert ^{4}<\infty$ for all $i=1,...,n$,
then $b_{j}<\infty$ for all $j$ with $\sigma_{j}>0.$ Under the finite second
moment assumption alone, $b_{j}$ may be infinite, in which case $\left(
\ref{eq:main}\right)  $ gives
\[
\mathbb{E}\left[  \left(  W_{n,j}-G_{j}\right)  ^{2}\right]  \leqslant
2\sigma_{j}^{2}.
\]

\end{remark}

\section{Estimates of covariance}

The main estimate uses the covariance of the sum to define a quadratic cost.
In its statement the normalization by $\sqrt{n}$ is absorbed into the summands.

\begin{theorem}
\label{thm:weighted} Let $U_{1},\ldots,U_{n}$ be independent random vectors in
$\mathbb{R}^{d}$ with $\mathbb{E}\left(  U_{i}\right)  =0,$ $\mathbb{E}
\left\vert U_{i}\right\vert ^{4}<\infty,$ $i=1,...,n$.. Let $S=\sum_{i=1}
^{n}U_{i}$, and assume $L=\operatorname{Cov}\left(  S\right)  $ is positive
definite. Let $\left(  U_{i}^{\prime}\right)  _{i=1}^{n}$ be an independent
copy of $\left(  U_{i}\right)  _{i=1}^{n}$ and $d_{i}=U_{i}^{\prime}-U_{i}$.
Then there exists a coupling between $S$ and $G^{\ast}\sim N\left(
0,L\right)  $ such that
\begin{equation}
\mathbb{E}\left[  \left(  S-G^{\ast}\right)  ^{\top}L\left(  S-G^{\ast
}\right)  \right]  \leqslant C_{0}\sum_{i=1}^{n}\mathbb{E}\left\vert
d_{i}\right\vert ^{4}. \label{eq:weighted}
\end{equation}

\end{theorem}

The proof is divided into several steps.

\section{Gaussian observation and the weighted cost}

\label{sec:observation}

For the proof of Theorem~\ref{thm:weighted}, we can assume without loss of
generality that the covariance $L=\operatorname{Cov}\left(  S\right)  $ is
diagonal. Indeed, since $L$ is symmetric, there exist $\left(  \lambda
_{1},\ldots,\lambda_{d}\right)  $ and an orthogonal matrix $Q$ with $Q^{\top
}LQ=\Lambda=\operatorname{diag}\left(  \lambda_{1},\ldots,\lambda_{d}\right)
$. The rotated random vectors $\bar{U}_{i}=Q^{\top}U_{i}$ are independent and
centered. Their sum is $\bar{S}=Q^{\top}S$ has $\operatorname{Cov}\left(
\bar{S}\right)  =Q^{\top}LQ=\Lambda$. The rotated differences $Q^{\top}d_{i}$
have norm $\left\vert d_{i}\right\vert $. Given a coupling between $S$ and
$\bar{G}=Q^{\top}G\sim N\left(  0,\Lambda\right)  $, thus $G_{\ast}\triangleq
Q\bar{G}\sim N(0,L)$ and the left side of $\left(  \ref{eq:weighted}\right)  $
is preserved by the identity
\[
\left(  \bar{S}-\bar{G}\right)  ^{\top}\Lambda\left(  \bar{S}-\bar{G}\right)
=\left(  S-G_{\ast}\right)  ^{\top}L\left(  S-G_{\ast}\right)  .
\]
Consequently, throughout Sections~\ref{sec:observation}-\ref{sec:integrate},
we use the eigenbasis of the covariance and assume that the covariance of $S$
has already been in a diagonal form $\operatorname{Cov}(S)=\operatorname{diag}
(\lambda_{1},\ldots,\lambda_{d}).$

\begin{proposition}
\label{prop:observation}Let $S\in\mathbb{R}^{d}$ be a random vector with
$\mathbb{E}S=0$, $\left\vert S\right\vert \leqslant M_{0}$ almost surely for
some $M_{0}<\infty$, and $L=\operatorname{Cov}\left(  S\right)
=\operatorname{diag}\left(  \lambda_{1},\ldots,\lambda_{d}\right)  $ with
every $\lambda_{j}>0$. Let $Z\sim N\left(  0,I_{d}\right)  $ be independent of
$S$. For $s>0$ define
\begin{equation}
Y_{s}=S+\sqrt{s}Z,\text{ }R_{s}=\mathbb{E}\left[  S-LZ/\sqrt{s}\mid
Y_{s}\right]  ,\text{ }w_{j}\left(  s\right)  =\frac{2\lambda_{j}^{2}
s}{\left(  s+\lambda_{j}\right)  ^{3}}. \label{eq:score}
\end{equation}
There exists $G^{\ast}\sim N\left(  0,L\right)  $ on a common probability
space with $S$ such that
\begin{equation}
\mathbb{E}\left[  \left(  S-G^{\ast}\right)  ^{\top}L\left(  S-G^{\ast
}\right)  \right]  =\sum_{j=1}^{d}\int_{0}^{\infty}w_{j}\left(  s\right)
\mathbb{E}R_{s,j}^{2}ds. \label{eq:identity}
\end{equation}

\end{proposition}

\begin{proof}
\textbf{1}. Suppose that $S$ is defined on the filtered probability space
$\left(  \Omega,\mathcal{A},\mathbb{P}\right)  $. Let $\mu=\mathcal{L}\left(
S\right)  $ be the distribution of $S$ $,$ $\mathbb{W}$ be one-dimensional
Wiener measure on the space $\mathcal{C}$ of continuous real paths starting at
$0$, with the Borel sigma-algebra for locally uniform convergence
\cite[Theorem~14.5]{kallenberg2021foundations}. We complete the extented
product space
\[
\left(  \Omega\times\mathcal{C}^{d},\ \mathcal{A}\otimes\mathcal{B}\left(
\mathcal{C}^{d}\right)  ,\ \mathbb{P}\otimes\mathbb{W}^{\otimes d}\right)
\]
by augmenting its filtration with the collection $\mathcal{N}$ of all null
sets. In the following we assume that the product space has already been
completed and use $\mathbb{P},$ $\mathbb{E}$ to denote the probability and
expectation on it. On the completed product space, the extended process
$S\left(  \omega,b\right)  =S(\omega)$ and the coordinate process
$B_{t}\left(  \omega,b\right)  =\left(  b_{1}\left(  t\right)  ,\ldots
,b_{d}\left(  t\right)  \right)  $ are independent. Define the observation
process
\begin{equation}
O_{t}=tS+B_{t}. \label{eq:observation}
\end{equation}
The raw observation and its completed filtration are given by
\begin{equation}
\mathcal{F}_{t}^{0}=\sigma\left(  O_{v}:0\leqslant v\leqslant t\right)
,\text{ }\mathcal{F}_{t}^{c}=\sigma\left(  \mathcal{F}_{t}^{0}\cup
\mathcal{N}\right)  . \label{eq:filtrations}
\end{equation}
Write $\mathcal{F}_{t}=\bigcap_{u>t}\mathcal{F}_{u}^{c}.$

\textbf{2}. Fix $t>0$ and a partition $\mathcal{P}:0=t_{0}<t_{1}<\cdots
<t_{m}=t$. For $1\leqslant\ell\leqslant m,$ write $h_{\ell}=t_{\ell}
-t_{\ell-1}>0$, $\Delta O_{\ell}=O_{t_{\ell}}-O_{t_{\ell-1}}$, and
$\mathbf{Z}_{\mathcal{P}}=\left(  \Delta O_{1},\ldots,\Delta O_{m}\right)  $.
Note $O_{0}=0$ by definition. For an initial observation $x\in\mathbb{R}^{d}$
and increment values $\mathbf{z}=\left(  z_{1},\ldots,z_{m}\right)  \in\left(
\mathbb{R}^{d}\right)  ^{m}$, consider
\[
\Delta O_{\ell}=h_{\ell}S+\left(  B_{t_{\ell}}-B_{t_{\ell-1}}\right)  .
\]
The Brownian increments are independent, they are also independent of $S$.
Writing the product Lebesgue measure as $d\mathbf{z}=dz_{1}\cdots dz_{m}.$ We
have the conditional densities of the the increment vector w.r.t. to the
Lebesgue measure,
\[
f_{\mathcal{P}}\left(  \mathbf{z}\mid x\right)  =\prod_{\ell=1}^{m}\left(
2\pi h_{\ell}\right)  ^{-d/2}\exp\left(  -\frac{\left\vert z_{\ell}-h_{\ell
}x\right\vert ^{2}}{2h_{\ell}}\right)  ,
\]
and
\[
f_{\mathcal{P}}\left(  \mathbf{z}\mid0\right)  =\prod_{\ell=1}^{m}\left(  2\pi
h_{\ell}\right)  ^{-d/2}\exp\left(  -\frac{\left\vert z_{\ell}\right\vert
^{2}}{2h_{\ell}}\right)  .
\]
These are likelihoods of the observation increments. As functions of $\left(
\mathbf{z},x\right)  ,$ they are Borel and strictly positive. The density
$f_{\mathcal{P}}\left(  \mathbf{z}\mid0\right)  $ is the driftless reference.
The increments has the likelihood ratio
\begin{equation}
\frac{f_{\mathcal{P}}(\mathbf{z}\mid x)}{f_{\mathcal{P}}(\mathbf{z}\mid
0)}=\prod_{\ell=1}^{m}e^{\left\langle x,z_{\ell}\right\rangle -h_{\ell}
|x|^{2}/2}=e^{\left\langle x,\sum_{i=1}^{m}z_{\ell}\right\rangle -t\left\vert
x\right\vert ^{2}/2}. \label{eq:likelihood}
\end{equation}

\textbf{3}. For any bounded Borel function $g:\mathbb{R}^{d}\mapsto\mathbb{R}$
and $y\in\mathbb{R}^{d}$, define
\begin{equation}
q_{g}\left(  t,y\right)  =\int g\left(  x\right)  e^{\left\langle
x,y\right\rangle -t\left\vert x\right\vert ^{2}/2}\mu\left(  dx\right)  .
\label{eq:q-definition}
\end{equation}
Write $q_{1}=q_{1}\left(  t,y\right)  $ for $g\equiv1$. Note that
$q_{1}\left(  t,y\right)  >0$ for all $\left(  t,y\right)  .$ To employ the
Bayes' theorem, we interpret $S$ as the parameter with $\mu$ as its prior and
$\mathbf{Z}_{\mathcal{P}}$ as the data. The joint density is $d\mathbf{z}
f_{\mathcal{P}}\left(  \mathbf{z}\mid x\right)  \mu(dx),$ thus the marginal
data density has
\[
p_{\mathcal{P}}(\mathbf{z})=\int f_{\mathcal{P}}\left(  \mathbf{z}\mid
u\right)  \mu\left(  du\right)  .
\]
Since $f_{\mathcal{P}}\left(  \mathbf{z}\mid u\right)  \leqslant\prod_{\ell
}(2\pi h_{\ell})^{-d/2}<\infty$, we have $0<p_{\mathcal{P}}\left(
\mathbf{z}\right)  <\infty$ for all $\mathbf{z}$. Using the likelihood ratio
$\left(  \ref{eq:likelihood}\right)  $ and the definition $\left(
\ref{eq:q-definition}\right)  ,$ we have
\[
p_{\mathcal{P}}\left(  \mathbf{z}\right)  =f_{\mathcal{P}}(\mathbf{z}
\mid0)\int e^{\langle x,y\rangle-t|x|^{2}/2}\mu\left(  dx\right)
=f_{\mathcal{P}}\left(  \mathbf{z}\mid0\right)  q_{1}\left(  t,y\right)  .
\]
The posterior kernel is
\begin{equation}
\kappa_{\mathcal{P}}\left(  \mathbf{z},dx\right)  =\frac{f_{\mathcal{P}
}\left(  \mathbf{z}\mid x\right)  \mu\left(  dx\right)  }{p_{\mathcal{P}
}\left(  \mathbf{z}\right)  }=\frac{e^{\left\langle x,y\right\rangle
-t|x|^{2}/2}\mu\left(  dx\right)  }{q_{1}\left(  t,y\right)  }.
\label{eq:posterior-kernel}
\end{equation}
For any $g$ as defined above, we have
\[
\int g\left(  x\right)  \kappa_{\mathcal{P}}\left(  \mathbf{z},dx\right)
=\frac{\int g\left(  x\right)  e^{\langle x,y\rangle-t|x|^{2}/2}\mu\left(
dx\right)  }{\int e^{\langle x,y\rangle-t|x|^{2}/2}\mu\left(  dx\right)
}=\frac{q_{g}\left(  t,y\right)  }{q_{1}\left(  t,y\right)  }.
\]
We now define the posterior process
\[
M_{t}^{g}\triangleq\frac{q_{g}\left(  t,O_{t}\right)  }{q_{1}\left(
t,O_{t}\right)  }.
\]
Hence for any Borel set $D\subset\left(  \mathbb{R}^{d}\right)  ^{m}$, we have
\begin{align*}
\mathbb{E}\left[  g\left(  S\right)  1_{\left\{  \mathbf{Z}_{\mathcal{P}}\in
D\right\}  }\right]   &  =\int_{D}\int g\left(  x\right)  f_{\mathcal{P}
}\left(  \mathbf{z}\mid x\right)  \mu\left(  dx\right)  d\mathbf{z=}\int
_{D}\int g\left(  x\right)  \kappa_{\mathcal{P}}\left(  \mathbf{z},dx\right)
p_{\mathcal{P}}\left(  \mathbf{z}\right)  d\mathbf{z}\\
&  \mathbf{=}\int_{D}\frac{q_{g}\left(  t,\sum_{\ell}z_{\ell}\right)  }
{q_{1}\left(  t,\sum_{\ell}z_{\ell}\right)  }p_{\mathcal{P}}\left(
\mathbf{z}\right)  d\mathbf{z}=\mathbb{E}\left[  M_{t}^{g}1_{\left\{
\mathbf{Z}_{\mathcal{P}}\in D\right\}  }\right]  .
\end{align*}
Here we used that $O_{t}=\sum_{\ell}\Delta O_{\ell}$. With the usual
finite-dimensinal cylinder set approximation \cite[Chapter 2]
{shreve1991brownian}\ and the monotone-class theorem \cite[Theorem
1.1]{kallenberg2021foundations}, we see that $M_{t}^{g}\left(  t,O_{t}\right)
$ is the version of the conditional expection under $\mathcal{F}_{t}^{0},$
\[
M_{t}^{g}=\mathbb{E}\left[  g\left(  S\right)  \mid\mathcal{F}_{t}^{0}\right]
.
\]
For $t=0,$ $O_{0}=0,$ $q_{1}\left(  0,O_{0}\right)  =1$ and $M_{0}^{g}=\int
g\left(  x\right)  \mu\left(  dx\right)  $. To work with Ito's theory, we want
$M_{t}^{g}$ to be adapted to a filtration that satisfies the usual condition
\cite[1.2.25]{shreve1991brownian}. $O_{t}$ has continuous path by definition.
By the form of $\left(  \ref{eq:q-definition}\right)  $ and the dominated
convergence theorem, we get that $M_{t}^{g}$ has continuous path. Fix
$t\geqslant0$ and let $u_{m}=t+1/m$, $m\geqslant1.$ Since every event of
$\mathcal{F}_{u_{m}}^{c}$ agrees almost surely with an event in $\mathcal{F}
_{t}^{0}$, so $\mathbb{E}\left[  g\left(  S\right)  \mid\mathcal{F}_{u_{m}
}^{c}\right]  =M_{u_{m}}^{g}.$ The sigma-algebras $\mathcal{F}_{u_{m}}^{c}$
decreases to $\mathcal{F}_{t}$ by the definition $\left(  \ref{eq:filtrations}
\right)  .$ Applying the backward case of \cite[Theorem 9.24]
{kallenberg2021foundations} to the integrable variable $g(S)$ and
$\mathcal{F}_{u_{m}}^{c}$, we obtain that
\[
\mathbb{E}\left[  g\left(  S\right)  \mid\mathcal{F}_{t}\right]
=\lim_{m\rightarrow\infty}\mathbb{E}\left[  g\left(  S\right)  \mid
\mathcal{F}_{u_{m}}^{c}\right]  =\lim_{m\rightarrow\infty}M_{u_{m}}^{g}
=M_{t}^{g}.
\]
The last equality uses the path continuity of $M_{t}^{g}$. Therefore, we have
proved for all $t\geqslant0,$
\begin{equation}
\mathbb{E}\left[  g\left(  S\right)  \mid\mathcal{F}_{t}\right]
=\mathbb{E}\left[  g\left(  S\right)  \mid O_{t}\right]  =\frac{q_{g}\left(
t,O_{t}\right)  }{q_{1}\left(  t,O_{t}\right)  }. \label{eq:posterior}
\end{equation}
Setting $g\left(  x\right)  =x$ yields the posterior mean vector
\begin{equation}
m_{t}=\mathbb{E}\left[  S\mid\mathcal{F}_{t}\right]  .
\label{eq:posterior-mean}
\end{equation}
We define the posterior covariance matrix
\begin{equation}
K_{t}=\operatorname{Cov}\left(  S\mid\mathcal{F}_{t}\right)  =\mathbb{E}
\left[  \left(  S-m_{t}\right)  \left(  S-m_{t}\right)  ^{\top}\mid
\mathcal{F}_{t}\right]  . \label{eq:posterior-variance}
\end{equation}
By the form of $\left(  \ref{eq:q-definition}\right)  $ and the dominated
convergence theorem, the kernel $q_{g}\left(  t,y\right)  $ is smooth w.r.t.
$\left(  t,y\right)  .$ So by differentiation we have
\[
\partial_{y_{\ell}}q_{g}=q_{gx_{\ell}},\text{ }\partial_{t}q_{g}
=-\frac{q_{g|x|^{2}}}{2}\text{ and }\sum_{\ell}\partial_{y_{\ell}}^{2}
q_{g}=q_{g|x|^{2}}.
\]
Here $q_{gx_{\ell}},$ $q_{g|x|^{2}}$ etc. mean replacing respectively $g$ by
$gx_{\ell},$ $g|x|^{2}$ in the definition $\left(  \ref{eq:q-definition}
\right)  .$ The last two equalities show that $q_{g}$ solves the backward heat
equation
\begin{equation}
\left(  \partial_{t}+\tfrac{1}{2}\Delta_{y}\right)  q_{g}=0.
\label{eq:backward-heat}
\end{equation}
Analogously $q_{x_{j}}$ means $q_{g}$ with $g\left(  x\right)  =x_{j}$. Recall
that $q_{1}\left(  t,y\right)  $ is positive for all $\left(  t,y\right)  ,$
so we can define
\[
p_{j}=\frac{q_{x_{j}}}{q_{1}}.
\]
With this notation, for each $j,$ the $j$-th component of $m_{t}$ can be
written as
\begin{equation}
m_{t,j}=p_{j}\left(  t,O_{t}\right)  \label{eq:mt-function}
\end{equation}
By differentiating $p_{j}$ through its definition, we have
\[
\partial_{y_{\ell}}p_{j}=\frac{q_{x_{j}x_{\ell}}}{q_{1}}-\frac{q_{x_{j}
}q_{x_{\ell}}}{q_{1}^{2}},\text{ }\partial_{y_{\ell}}\log q_{1}=p_{\ell}.
\]
and
\[
\partial_{y_{\ell}}p_{j}\left(  t,O_{t}\right)  =\mathbb{E}\left[
S_{j}S_{\ell}\mid\mathcal{F}_{t}\right]  -\mathbb{E}\left[  S_{j}
\mid\mathcal{F}_{t}\right]  \mathbb{E}\left[  S_{\ell}\mid\mathcal{F}
_{t}\right]  =\left(  K_{t}\right)  _{j,\ell}.
\]
The backward heat equation is applicable to $g\left(  x\right)  =x_{j}$ and
$g\left(  x\right)  =1$. So applying the oeprator $\partial_{t}+\tfrac{1}
{2}\Delta_{y}$ to $q_{x_{j}}=p_{j}q_{1}$ yields that
\[
0=q_{1}\left(  \partial_{t}+\tfrac{1}{2}\Delta_{y}\right)  p_{j}+\left\langle
\nabla_{y}p_{j},\nabla_{y}q_{1}\right\rangle .
\]
Combing these computations, we get
\begin{equation}
\left(  \partial_{t}+\tfrac{1}{2}\Delta_{y}\right)  p_{j}=-\left\langle
\nabla_{y}p_{j},\nabla_{y}\log q_{1}\right\rangle . \label{eq:pde}
\end{equation}

\textbf{4}. Define the innovation process which gives the observation after
subtracting the accumulated predictable drift in the observation filtration,
\[
\mathcal{I}_{0}=0,\text{ }\mathcal{I}_{t}=O_{t}-\int_{0}^{t}m_{s}ds.
\]
Thus $\mathcal{I}_{t}$ is adapted to $\mathcal{F}_{t}^{0}\subset
\mathcal{F}_{t}.$ Recall the definition of $O_{t}$ in $\left(
\ref{eq:observation}\right)  .$ Let $s<t$. Since $S$ is independent of the
whole path filtration of $B,$ and $B$ has independent increments, thus the
increment $B_{t}-B_{s}$ is independent of $\sigma\left(  S,B_{r}:r\leqslant
s\right)  \supset\mathcal{F}_{s}^{0},$ hence
\[
\mathbb{E}\left[  O_{t}-O_{s}\mid\mathcal{F}_{s}^{0}\right]  =\left(
t-s\right)  \mathbb{E}\left[  S\mid\mathcal{F}_{s}^{0}\right]  +\mathbb{E}
\left[  B_{t}-S_{s}\mid\mathcal{F}_{s}^{0}\right]  =(t-s)m_{s}.
\]
Also since $\left\vert m_{u}\right\vert \leqslant M_{0},$ so Fubini's theorem
applies and gives,
\[
\mathbb{E}\left[  \int_{s}^{t}m_{u}du\mid\mathcal{F}_{s}^{0}\right]  =\int
_{s}^{t}\mathbb{E}\left[  \mathbb{E}\left[  S\mid\mathcal{F}_{u}^{0}\right]
\mid\mathcal{F}_{s}^{0}\right]  du=\left(  t-s\right)  m_{s}.
\]
Each component $\mathcal{I}_{t,j}$ of $\mathcal{I}_{t}$ has
\[
\mathbb{E}\left\vert \mathcal{I}_{t,j}\right\vert \leqslant2tM_{0}
+\mathbb{E}\left\vert B_{t,j}\right\vert \leqslant2tM_{0}+t^{1/2}<\infty.
\]
Hence $\mathcal{I}_{\cdot,j}$ is a $\mathcal{F}_{t}^{0}$-martingale. Here we
write $\mathcal{I}_{\cdot,j}$ for the $j$-th component of the process
$\mathcal{I}_{\cdot}$. Using \cite[Lemma 9.8]{kallenberg2021foundations} and
applying \cite[Theorem~2.8]{revuz2013continuous} to $\mathcal{I}_{\cdot},$
$-\mathcal{I}_{\cdot}$, we see that $\mathcal{I}_{\cdot}$ is also a
$\mathcal{F}_{t}$-martingale. Since
\[
\mathcal{I}_{t}-B_{t}=\int_{0}^{t}\left(  S-m_{u}\right)  du,\text{ }\forall
t\geqslant0,
\]
and the integral has zero quadratic variation. So the quadratic variation of
$\mathcal{I}_{\cdot}$ satisfies
\[
\left[  \mathcal{I}_{\cdot,j},\mathcal{I}_{\cdot,\ell}\right]  _{t}=\left[
B_{\cdot,j},B_{\cdot,\ell}\right]  _{t}=t1_{\{j=\ell\}}.
\]
This proves, by L\'{e}vy's criterion \cite[Theorem 19.3]
{kallenberg2021foundations}, that $\mathcal{I}_{t}$ is also a Brownian motion
w.r.t. the complete right-continuous filtration $\mathcal{F}_{t}.$ Note that
$p_{j}\left(  t,y\right)  $ is a smooth on any bounded open set of
$\mathbb{R\times R}^{d}$ because in its definition $q_{x_{j}}$ is smooth and
$q_{1}$ is positive everywhere$.$ In particularly $p_{j}$ is smooth in any
bounded open set containing $\left(  0,0\right)  .$ So we can invoke Ito's
lemma \cite[Theorem 18.18]{kallenberg2021foundations} to $\left(
\ref{eq:mt-function}\right)  $ through the smooth function $p_{j}$ starting
with $t=0.$ Since
\[
dO_{t,\ell}=m_{t,\ell}dt+d\mathcal{I}_{t,\ell}\text{ and }d\left[
O_{\cdot,\ell},O_{\cdot,h}\right]  _{t}=1_{\left\{  \ell=h\right\}  }dt,
\]
we therefore obtain that
\begin{align*}
dm_{t,j}  &  =\left(  \partial_{t}p_{j}+\tfrac{1}{2}\Delta_{y}p_{j}+\sum
_{\ell=1}^{d}m_{t,\ell}\partial_{y_{\ell}}p_{j}\right)  \left(  t,O_{t}
\right)  \,dt+\sum_{\ell=1}^{d}\partial_{y_{\ell}}p_{j}\left(  t,O_{t}\right)
d\mathcal{I}_{t,\ell}\\
&  =\left(  -\sum_{\ell=1}^{d}\left(  K_{t}\right)  _{j,\ell}m_{t,\ell}
+\sum_{\ell=1}^{d}m_{t,\ell}\left(  K_{t}\right)  _{j,\ell}\right)
dt+\sum_{\ell=1}^{d}\left(  K_{t}\right)  _{j,\ell}d\mathcal{I}_{t,\ell}
=\sum_{\ell=1}^{d}\left(  K_{t}\right)  _{j,\ell}d\mathcal{I}_{t,\ell}.
\end{align*}
To interpret it in the Ito's sense of integral, we need to show that $K_{t}$
is a square integrable process on any finite time interval. Indeed, since $S$
is positive semi-definite, for any vector $v$, we have
\[
v^{\top}K_{t}v=\operatorname{Var}\left[  \left\langle v,S\right\rangle
\mid\mathcal{F}_{t}\right]  \geqslant0.
\]
So $K_{t}$ is positive semi-definite. Using the trace property and $\left(
\ref{eq:posterior}\right)  $ yields
\[
\left\Vert K_{t}\right\Vert \leqslant\operatorname{tr}K_{t}=\mathbb{E}\left[
\left\vert S\right\vert ^{2}\mid\mathcal{F}_{t}\right]  -\left\vert
m_{t}\right\vert ^{2}\leqslant M_{0}^{2}.
\]
This shows that $K_{t}$ is square integrable, i.e. for any $T>0,$
\[
\mathbb{E}\int_{0}^{T}\left\Vert K_{t}\right\Vert ^{2}dt=\sum_{j,\ell=1}
^{d}\int_{0}^{T}\mathbb{E}\left\vert K_{u}\right\vert _{j,\ell}^{2}du\leqslant
TM_{0}^{4}<\infty.
\]
Thus we have proved the Ito's differential equation,
\begin{equation}
dO_{t}=m_{t}dt+d\mathcal{I}_{t}\text{ and }dm_{t}=K_{t}d\mathcal{I}_{t}\text{
with }m_{0}=0. \label{eq:innovation}
\end{equation}

\textbf{5}. We next compute a covariance that results from a Gaussian prior,
i.e. keeping the same observation model as above but with a Gassian prior
$N(0,L)$. Recall that $L$ is positive definite. Replacing the previous prior
$\mu$ by the Gaussian in the posterior kernel $\left(
\ref{eq:posterior-kernel}\right)  $ after observing $O_{t}=y$ at the time $t,$
we have the unnormalized posterior density,
\begin{align*}
&  \left(  2\pi\right)  ^{-d/2}\left(  \det L\right)  ^{-1/2}\exp\left(
-\tfrac{1}{2}\left\langle x,L^{-1}x\right\rangle \right)  \exp\left(
\left\langle x,y\right\rangle -t\tfrac{1}{2}|x|^{2}\right) \\
&  =\left(  2\pi\right)  ^{-d/2}\left(  \det L\right)  ^{-1/2}\exp\left(
\left\langle x,y\right\rangle -\tfrac{1}{2}\left\langle x,P_{t}x\right\rangle
\right)  ,
\end{align*}
where $P_{t}=L^{-1}+tI_{d}$. We now identify the exact distribution of this
posterior. Denote the inverse of $P_{t}$ by $D_{t}=P_{t}^{-1}$. Since
$P_{t}D_{t}=I_{d}$ and both matrices are symmetric, we have
\begin{align*}
\left\langle x-D_{t}y,P_{t}\left(  x-D_{t}y\right)  \right\rangle  &
=\left\langle x,P_{t}x\right\rangle -\left\langle x,P_{t}D_{t}y\right\rangle
-\left\langle y,D_{t}P_{t}x\right\rangle +\left\langle y,D_{t}P_{t}
D_{t}y\right\rangle \\
&  =\left\langle x,P_{t}x\right\rangle -2\left\langle x,y\right\rangle
+\left\langle y,D_{t}y\right\rangle .
\end{align*}
Thus the exponent of the posterior density is
\[
\left(  2\pi\right)  ^{-d/2}\left(  \det L\right)  ^{-1/2}\exp\left[
-\tfrac{1}{2}\left\langle x-D_{t}y,P_{t}\left(  x-D_{t}y\right)  \right\rangle
+\tfrac{1}{2}\left\langle y,D_{t}y\right\rangle \right]
\]
With a change of variable $z=D_{t}^{-1/2}\left(  x-D_{t}y\right)  $ we have
\[
\int e^{-\left\langle x-D_{t}y,P_{t}\left(  x-D_{t}y\right)  \right\rangle
/2}dx=\left(  \det D_{t}\right)  ^{1/2}\int e^{-\left\vert z\right\vert
^{2}/2}dz=\left(  2\pi\right)  ^{d/2}\left(  \det D_{t}\right)  ^{1/2}.
\]
Therefore after normalization we see that the posterior distribution is
$N\left(  D_{t}y,D_{t}\right)  $. So when $O_{t}$ is observed, the posterior
mean and covariance will be $D_{t}O_{t},$ $D_{t}.$ They respectively play the
roles of $m_{t},$ $K_{t}$ in our early setup. At $t=0$, $P_{0}=L^{-1}$ and
$D_{0}=L$, recovering the prior.

Rewriting $P_{t}=L^{-1}(I_{d}+tL)$, we obtain that $D_{t}=(I_{d}+tL)^{-1}L$.
In the chosen eigenbasis, its $j$-th diagonal entry is $\left(  \lambda
_{j}^{-1}+t\right)  ^{-1}=\lambda_{j}/\left(  1+t\lambda_{j}\right)  $. Hence,
for $0\leqslant a<b<\infty$, we have
\[
\int_{a}^{b}\left(  D_{t}\right)  _{j,j}^{2}dt=\left.  -\frac{\lambda_{j}
}{1+t\lambda_{j}}\right\vert _{a}^{b}=\frac{\lambda_{j}}{1+a\lambda_{j}}
-\frac{\lambda_{j}}{1+b\lambda_{j}}.
\]
Define, for $T>0,$
\[
G_{T}^{\ast}=\int_{0}^{T}D_{t}d\mathcal{I}_{t}.
\]
Using Ito's isometry we get, for $0<T<T^{\prime},$
\begin{align*}
\mathbb{E}\left\vert G_{T^{\prime}}^{\ast}-G_{T}^{\ast}\right\vert ^{2}  &
=\mathbb{E}\left\vert \int_{T}^{T^{\prime}}D_{t}d\mathcal{I}_{t}\right\vert
^{2}=\mathbb{E}\int_{T}^{T^{\prime}}\left\Vert D_{t}\right\Vert ^{2}dt\\
&  =\sum_{j}\int_{T}^{T^{\prime}}\frac{\lambda_{j}^{2}}{\left(  1+t\lambda
_{j}\right)  ^{2}}\,dt\leqslant\sum_{j}\frac{\lambda_{j}}{1+T\lambda_{j}
}\rightarrow0\text{ as }T\rightarrow\infty.
\end{align*}
Thus we can define$,$
\[
G^{\ast}=\lim_{T\rightarrow\infty}G_{T}^{\ast}=\int_{0}^{\infty}
D_{t}d\mathcal{I}_{t}\text{ in }L^{2}\left(  d\mathbb{P}\right)  .
\]
For each $\theta\in\mathbb{R}^{d}$, $\theta^{\top}G^{\ast}$ is a centered
Gaussian. Indeed,
\[
\int_{0}^{\infty}\left\vert D_{t}^{\top}\theta\right\vert ^{2}dt=\sum
_{j}\theta_{j}^{2}\int_{0}^{\infty}\left(  D_{t}\right)  _{j,j}^{2}
dt=\theta^{\top}L\theta.
\]
Since $\left\vert e^{\mathrm{i}a}-e^{\mathrm{i}b}\right\vert \leqslant|a-b|,$
the convergence in $L^{2}$ guarantees the convergence of the characteristic
functions. Hence $\mathbb{E}e^{\mathrm{i}\theta^{\top}G^{\ast}}=e^{-\theta
^{\top}L\theta/2}$, which identifies the required law. Hence $G^{\ast}\sim
N\left(  0,L\right)  $.

\textbf{6}. For any square-integrable $\mathcal{F}_{t}$-measurable vector
$V$,
\[
\mathbb{E}\left\langle S-m_{t},m_{t}-V\right\rangle =\mathbb{E}\left\langle
\mathbb{E}[S-m_{t}\mid\mathcal{F}_{t}],m_{t}-V\right\rangle =0.
\]
It follows that
\[
\mathbb{E}\left\vert S-V\right\vert ^{2}=\mathbb{E}\left\vert S-m_{t}
\right\vert ^{2}+\mathbb{E}\left\vert m_{t}-V\right\vert ^{2}.
\]
Taking $V=O_{t}/t$ and using $\mathbb{E}\left\vert B_{t}\right\vert ^{2}=dt$
we obtain that
\begin{equation}
\mathbb{E}\left\vert S-m_{t}\right\vert ^{2}\leqslant\mathbb{E}\left\vert
S-O_{t}/t\right\vert ^{2}=t^{-2}\mathbb{E}\left\vert B_{t}\right\vert
^{2}=d/t. \label{eq:posterior-mse}
\end{equation}
Therefore
\[
\lim_{t\rightarrow\infty}m_{t}=S\text{ in }L^{2}\left(  d\mathbb{P}\right)
\text{ as }t\rightarrow\infty.
\]
By $\left(  \ref{eq:innovation}\right)  $ and Ito's isometry,
\[
\int_{0}^{T}\mathbb{E}\left\Vert K_{t}\right\Vert ^{2}dt=\mathbb{E}\left\vert
m_{T}\right\vert ^{2}\leqslant\mathbb{E}\left\vert S\right\vert ^{2}
=\operatorname{tr}L<\infty.
\]
Thus by $\left(  \ref{eq:posterior-mse}\right)  $ we obtain
\begin{equation}
S-G^{\ast}=\lim_{T\rightarrow\infty}(m_{T}-G_{T}^{\ast})=\int_{0}^{\infty
}\left(  K_{t}-D_{t}\right)  d\mathcal{I}_{t}\text{ in }L^{2}\left(
d\mathbb{P}\right)  . \label{eq:error-integral}
\end{equation}
\textbf{7}. Let
\begin{equation}
A_{t}=\left(  I_{d}+tL\right)  m_{t}-LO_{t}. \label{eq:At-time-change}
\end{equation}
By $\left(  \ref{eq:innovation}\right)  $ and $\left(  I_{d}+tL\right)
D_{t}=L$ we have
\[
dA_{t}=Lm_{t}\,dt+\left(  I_{d}+tL\right)  K_{t}d\mathcal{I}_{t}-L\left(
m_{t}dt+d\mathcal{I}_{t}\right)  =\left(  I_{d}+tL\right)  \left(  K_{t}
-D_{t}\right)  d\mathcal{I}_{t}.
\]
Here $A_{0}=0$ and $A$ is square integrable on each finite time interval.
Denote the $j$-th row of $A_{t}$ by $A_{t,j}$ and define
\[
h_{j}(t)=\mathbb{E}A_{t,j}^{2},
\]
By Ito's isometry
\[
h_{j}(t)=\int_{0}^{t}\left(  1+u\lambda_{j}\right)  ^{2}\sum_{\ell}
\mathbb{E}\left\vert \left(  K_{u}-D_{u}\right)  _{j,\ell}\right\vert ^{2}du.
\]
Thus $h_{j}$ is absolutely continuous and its derivative equals the integrand
almost everywhere, i.e.
\[
\frac{dh_{j}\left(  t\right)  }{\left(  1+t\lambda_{j}\right)  ^{2}}
=\sum_{\ell}\mathbb{E}|(K_{t}-D_{t})_{j,\ell}|^{2}dt.
\]
On the other hand, applying the same isometry to row $j$ of $\left(
\ref{eq:error-integral}\right)  $ yields
\[
\mathbb{E}\left(  S_{j}-G_{j}^{\ast}\right)  ^{2}=\int_{0}^{\infty}\sum_{\ell
}\mathbb{E}\left\vert \left(  K_{t}-D_{t}\right)  _{j,\ell}\right\vert
^{2}dt.
\]
Since
\[
\frac{\lambda_{j}}{\left(  1+t\lambda_{j}\right)  ^{2}}=\int_{t}^{\infty}
\frac{2\lambda_{j}^{2}}{\left(  1+u\lambda_{j}\right)  ^{3}}du,
\]
we get by Tonelli's theorem that
\begin{align*}
\lambda_{j}\mathbb{E}\left(  S_{j}-G_{j}^{\ast}\right)  ^{2}  &  =\int
_{0}^{\infty}\frac{\lambda_{j}}{\left(  1+t\lambda_{j}\right)  ^{2}}
dh_{j}\left(  t\right) \\
&  =\int_{0}^{\infty}\int_{t}^{\infty}\frac{2\lambda_{j}^{2}}{\left(
1+u\lambda_{j}\right)  ^{3}}dudh_{j}\left(  t\right)  =\int_{0}^{\infty}
\frac{2\lambda_{j}^{2}h_{j}\left(  u\right)  }{\left(  1+u\lambda_{j}\right)
^{3}}du.
\end{align*}
The interchange uses $\int_{[0,u]}dh_{j}(t)=h_{j}(u)$, since $h_{j}(0)=0$.
Substituing $\left(  \ref{eq:posterior-mean}\right)  $ into $\left(
\ref{eq:At-time-change}\right)  $ and using that $B_{t}=O_{t}-tS$ we get
\[
A_{t}=\left(  I_{d}+tL\right)  \mathbb{E}\left[  S\mid O_{t}\right]
-LO_{t}=\mathbb{E}\left[  \left(  I_{d}+tL\right)  S-LO_{t}\mid O_{t}\right]
=\mathbb{E}\left[  S-LB_{t}\mid O_{t}\right]  .
\]
Fix $t>0$ and let $s=1/t$. The variables $O_{t}$ and $O_{t}/t$ generate the
same sigma-field because each is a measurable function of the other. Hence
\[
A_{t}=\mathbb{E}\left[  S-LB_{t}\mid O_{t}\right]  =\mathbb{E}\left[
S-LB_{t}\mid O_{t}/t\right]  .
\]
Rewrite it as
\[
A_{t}=\mathbb{E}\left[  S-L\left(  B_{t}/\sqrt{t}\right)  /\sqrt{s}\mid
O_{t}/t\right]  .
\]
Recall the definition of $R_{s}$ in $\left(  \ref{eq:score}\right)  ,$
\[
R_{s}=\mathbb{E}\left[  S-LZ/\sqrt{s}\mid Y_{s}\right]  \text{ with }
Y_{s}=S+\sqrt{s}Z.
\]
Since $Z^{(t)}\triangleq B_{t}/\sqrt{t}$ is a standard Gaussian independent of
$S$ and $O_{t}/t=S+\sqrt{s}Z^{(t)}.$ Thus $(S,O_{t}/t)$ and $(S,Y_{s})$ have
the same joint distribution at this fixed time. We see that
\[
h_{j}(t)=\mathbb{E}A_{t,j}^{2}=\mathbb{E}R_{s,j}^{2}.
\]
Therefore, under the change of variable $u=1/s$, we have
\[
\lambda_{j}\mathbb{E}\left(  S_{j}-G_{j}^{\ast}\right)  ^{2}=\int_{0}^{\infty
}\frac{2\lambda_{j}^{2}s^{-2}}{\left(  1+\lambda_{j}/s\right)  ^{3}}
\mathbb{E}R_{s,j}^{2}ds=\int_{0}^{\infty}w_{j}\left(  s\right)  \mathbb{E}
R_{s,j}^{2}ds.
\]
Summing over $j$ we get $\left(  \ref{eq:identity}\right)  $.
\end{proof}

The integrations in Section~\ref{sec:integrate} use three properties of the
weights in $\left(  \ref{eq:score}\right)  $.
\begin{equation}
\sup_{s>0}w_{j}(s)=\frac{8}{27},\text{ }\int_{0}^{\infty}w_{j}(s)\frac{ds}
{s}=1,\text{ }\frac{w_{j}(s)}{s}\leqslant\frac{2}{\lambda_{j}}.
\label{eq:weights}
\end{equation}
Indeed, for $z=s/\lambda_{j}$, the derivative of $w_{j}(s)=2z/(1+z)^{3}$
w.r.t. $z$ is $2(1-2z)/(1+z)^{4}$. It is positive for $0<z<1/2$ and negative
for $z>1/2$, hence
\[
\sup_{s>0}w_{j}(s)=w_{j}(\lambda_{j}/2)=\frac{2(1/2)}{(1+1/2)^{3}}=\frac
{8}{27}.
\]
Also
\[
\int_{0}^{\infty}2\lambda_{j}^{2}(s+\lambda_{j})^{-3}ds=\left.  -\lambda
_{j}^{2}(s+\lambda_{j})^{-2}\right\vert _{0}^{\infty}=1.
\]
Finally
\[
\frac{w_{j}(s)}{s}=\frac{2\lambda_{j}^{2}}{(s+\lambda_{j})^{3}}\leqslant
\frac{2\lambda_{j}^{2}}{\lambda_{j}^{3}}=\frac{2}{\lambda_{j}}.
\]

\section{Resampling and expansion at the noise scale}

\label{sec:hermite}

Recall that the symbol $\otimes$ denotes the tensor product. For vectors $u,v$
its entries are $(u\otimes v)_{\ell,h}=u_{\ell}v_{h}$. For $v\in\mathbb{R}
^{d}$ and an integer $k\geqslant1$, the tensor power $v^{\otimes k}$ has
entries
\begin{equation}
(v^{\otimes k})_{\ell_{1},\ldots,\ell_{k}}=\prod_{p=1}^{k}v_{\ell_{p}},\text{
}1\leqslant\ell_{1},\ldots,\ell_{k}\leqslant d. \label{eq:tensor-power}
\end{equation}
Set $v^{\otimes0}=1$. For two real tensors $T,V$ of order $k$ we use the inner
product $\langle T,V\rangle=\sum_{\ell_{1},\ldots,\ell_{k}}T_{\ell_{1}
,\ldots,\ell_{k}}V_{\ell_{1},\ldots,\ell_{k}}$ and its norm $\Vert
T\Vert=\langle T,T\rangle^{1/2}$. A tensor is \textit{symmetric} if permuting
its indices leaves each entry unchanged. In particular,
\[
\Vert v^{\otimes k}\Vert^{2}=(\sum_{\ell=1}^{d}v_{\ell}^{2})^{k}=|v|^{2k}.
\]
For independent pairs $(U_{i},U_{i}^{\prime})$, $1\leqslant i\leqslant n,$
write
\[
d_{i}=U_{i}^{\prime}-U_{i}\text{ and }r_{i}=|d_{i}|^{2}.
\]
In an eigen basis $\left(  e_{j}\right)  _{j=1}^{d}$, write
\[
d_{i,j}=\langle d_{i},e_{j}\rangle=U_{i,j}^{\prime}-U_{i,j},\text{
}j=1,...,d.
\]
For $s>0$ and $1\leqslant j\leqslant d$, let
\begin{equation}
b_{i,j}(s)=\mathbb{E}[d_{i,j}1_{\{r_{i}>s\}}\mid U_{i}],\text{ }
C_{i,j,k}(s)=\mathbb{E}[d_{i,j}d_{i}^{\otimes k}1_{\{r_{i}\leqslant s\}}\mid
U_{i}],\text{ }k\geqslant1. \label{eq:coeff}
\end{equation}
With $e_{j}$ denoting the $j$-th basis vector, define
\begin{equation}
T_{j,1}=\sum_{i}C_{i,j,1}(s)-2\lambda_{j}e_{j},\text{ }T_{j,k}=\sum
_{i}C_{i,j,k}(s)\text{ for }k\geqslant2. \label{eq:tensors}
\end{equation}
We suppress the argument $s$ on $T_{j,k}$.

The lemma gives the Hermite expansion used in the computation of the bounds of
Lemma \ref{lem:expansion}.

\begin{lemma}
\label{lm:hermite-expansion}For $v,$ $z\in\mathbb{R}^{d}$, let
\[
E_{v}\left(  z\right)  =\exp\left(  \left\langle z,v\right\rangle
-|v|^{2}/2\right)  .
\]
Let $H_{0}=1$. For $k\geqslant1$ define the symmetric tensor $H_{k}\left(
z\right)  $ by
\begin{equation}
\left(  H_{k}\left(  z\right)  \right)  _{\ell_{1},\ldots,\ell_{k}}=\left.
\partial_{v_{\ell_{1}}}\cdots\partial_{v_{\ell_{k}}}E_{v}\left(  z\right)
\right\vert _{v=0}. \label{eq:hermite-differential}
\end{equation}
Then the following hold.

(i) For symmetric tensors $T,V$ of orders $k,\ell\geqslant0$ respectively, we
have the isometry
\begin{equation}
\mathbb{E}\left[  \left\langle H_{k}\left(  Z\right)  ,T\right\rangle
\left\langle H_{\ell}\left(  Z\right)  ,V\right\rangle \right]  =\left\{
\begin{array}
[c]{ll}
k!\left\langle T,V\right\rangle , & k=\ell,\\
0, & k\neq\ell.
\end{array}
\right.  \label{eq:orthog}
\end{equation}

(ii) For $Z\sim\gamma_{d},$ the series converges in $L^{2}\left(  \gamma
_{d}\right)  $ uniformly in $\left\vert v\right\vert \leqslant1$,
\begin{equation}
E_{v}\left(  Z\right)  =\sum_{k=0}^{\infty}\frac{\left\langle H_{k}\left(
Z\right)  ,v^{\otimes k}\right\rangle }{k!}. \label{eq:hermite}
\end{equation}

\end{lemma}

\begin{proof}
\textbf{1}. Using the identity that $\partial_{v_{j}}E_{v}\left(  z\right)
=\left(  z_{j}-v_{j}\right)  E_{v}\left(  z\right)  ,$ we have $H_{1}\left(
z\right)  =z$ at $v=0$. Fix $u\in\mathbb{R}^{d}$. We have
\begin{equation}
\left.  \partial_{t}^{k}E_{tu}\left(  z\right)  \right\vert _{t=0}=\sum
_{a_{1},\ldots,a_{k}\in\left\{  1,...,d\right\}  }u_{a_{1}}\cdots u_{a_{k}
}\left(  H_{k}\left(  z\right)  \right)  _{a_{1},\ldots,a_{k}}=\left\langle
H_{k}\left(  z\right)  ,u^{\otimes k}\right\rangle . \label{eq:hermite-co}
\end{equation}
Recall that $\mathbb{E}e^{\left\langle a,Z\right\rangle }=e^{\left\vert
a\right\vert ^{2}/2}$ for $a\in\mathbb{R}^{d}$. For $u,v\in\mathbb{R}^{d}$ and
$t,h\in\mathbb{R}$, taking $a=tu+hv$ in the previous formula and using
\[
\tfrac{1}{2}\left\vert tu+hv\right\vert ^{2}-\tfrac{1}{2}t^{2}\left\vert
u\right\vert ^{2}-\tfrac{1}{2}h^{2}\left\vert v\right\vert ^{2}=th\left\langle
u,v\right\rangle
\]
we get
\begin{equation}
\mathbb{E}\left[  E_{tu}\left(  Z\right)  E_{hv}\left(  Z\right)  \right]
=e^{th\left\langle u,v\right\rangle }. \label{eq:gen-hermite}
\end{equation}
Write $b=\left\langle u,v\right\rangle .$ Since Gaussian density decays faster
than any polynomial, differentiations pass through the above expectation.
Observe that the coefficient of $t^{k}h^{\ell}$ in the expansion
\[
e^{thb}=\sum_{r\geqslant0}t^{r}h^{r}\frac{b^{r}}{r!},
\]
vanishes for $k\neq\ell$ and equals $b^{k}/k!$ for $k=\ell$. The mixed
derivative at $t=h=0$ of the right side of $\left(  \ref{eq:gen-hermite}
\right)  $ is this coefficient times $k!\ell!$. Using this observation and
that the corresponding coefficients of the two sides of $\left(
\ref{eq:gen-hermite}\right)  $ after differentiation are equal, we therefore
get
\begin{equation}
\mathbb{E}\left[  \left\langle H_{k}\left(  Z\right)  ,u^{\otimes
k}\right\rangle \left\langle H_{\ell}\left(  Z\right)  ,v^{\otimes\ell
}\right\rangle \right]  =\left.  \partial_{t}^{k}\partial_{h}^{\ell}
e^{thb}\right\vert _{t=h=0}=\left\{
\begin{array}
[c]{ll}
(k!)^{2}b^{k}/k!=k!b^{k}, & k=\ell,\\
0, & k\neq\ell.
\end{array}
\right.  \label{eq:ts}
\end{equation}
For $k=\ell$, using $\left(  \ref{eq:tensor-power}\right)  $ we can write the
tensor inner product as
\begin{equation}
\left\langle u^{\otimes k},v^{\otimes k}\right\rangle =\sum_{a_{1}
,\ldots,a_{k}\in\left\{  1,...,d\right\}  }\prod_{p=1}^{k}u_{a_{p}}v_{a_{p}
}=\left(  \sum_{a=1}^{d}u_{a}v_{a}\right)  ^{k}=b^{k}. \label{eq:ts-power}
\end{equation}

\textbf{2}. Consider symmetric tensors $T,V$ of positive orders $k,\ell
\geqslant1$, there exist real representations \cite[Section~8.2,
Eq.~(8.3)]{comon2008symmetric} that $T=\sum_{a}c_{a}u_{a}^{\otimes k}$ and
$V=\sum_{b}\beta_{b}v_{b}^{\otimes\ell}$ with real coefficients $c_{a}
,\beta_{b}$. Using this representations together with $\left(  \ref{eq:ts}
\right)  $ and substituting $\left(  \ref{eq:ts-power}\right)  $ we have
\begin{align*}
&  \mathbb{E}\left[  \left\langle H_{k}\left(  Z\right)  ,T\right\rangle
\left\langle H_{\ell}\left(  Z\right)  ,V\right\rangle \right] \\
&  =\sum_{a,b}c_{a}\beta_{b}\mathbb{E}\left[  \left\langle H_{k}\left(
Z\right)  ,u_{a}^{\otimes k}\right\rangle \left\langle H_{\ell}(Z),v_{b}
^{\otimes\ell}\right\rangle \right]  =\left\{
\begin{array}
[c]{ll}
k!\sum_{a,b}c_{a}\beta_{b}\left\langle u_{a}^{\otimes k},v_{b}^{\otimes\ell
}\right\rangle , & k=\ell,\\
0, & k\neq\ell.
\end{array}
\right.  .
\end{align*}
At unequal orders each summand vanishes due to $\left(  \ref{eq:ts}\right)  $.
The sum of equal order terms is
\[
k!\sum_{a,b}c_{a}\beta_{b}\langle u_{a}^{\otimes k},v_{b}^{\otimes k}
\rangle=k!\langle T,V\rangle.
\]
If both $T,V$ are both zero order tensor, $k=0,$ $\ell=0$, then $T,$ $V$ are
constant. Also $H_{0}\left(  Z\right)  =1.$ Thus
\[
\mathbb{E}\left[  \left\langle H_{k}\left(  Z\right)  ,T\right\rangle
\left\langle H_{\ell}\left(  Z\right)  ,V\right\rangle \right]  =\mathbb{E}
\left[  \left\langle H_{0}\left(  Z\right)  ,T\right\rangle \left\langle
H_{0}\left(  Z\right)  ,V\right\rangle \right]  =TV=\left\langle
T,V\right\rangle .
\]
If either one of $T,V$ has order zero, say $k=0,$ $\ell\geqslant1$, then
\[
\mathbb{E}\left[  \left\langle H_{k}\left(  Z\right)  ,T\right\rangle
\left\langle H_{\ell}\left(  Z\right)  ,V\right\rangle \right]  =T\mathbb{E}
\left[  \left\langle H_{\ell}\left(  Z\right)  ,V\right\rangle \right]
\]
Setting $t=0$ in $\left(  \ref{eq:gen-hermite}\right)  $ and using
$E_{0}\left(  Z\right)  =1$, we obtain $\mathbb{E}E_{hv}\left(  Z\right)  =1$
for every $h$. Since $\ell\geqslant1,$ differentiation then yields
\[
\mathbb{E}\left\langle H_{\ell}\left(  Z\right)  ,v^{\otimes\ell}\right\rangle
=\left.  \frac{d^{\ell}}{dh^{\ell}}\mathbb{E}E_{hv}\left(  Z\right)
\right\vert _{h=0}=\left.  \frac{d^{\ell}}{dh^{\ell}}1\right\vert _{h=0}=0.
\]
Thus $\mathbb{E}\left\langle H_{\ell}\left(  Z\right)  ,V\right\rangle =0.$
Combining the above computations we obtain
\[
\mathbb{E}\left[  \left\langle H_{k}\left(  Z\right)  ,T\right\rangle
\left\langle H_{\ell}\left(  Z\right)  ,V\right\rangle \right]  =\left\{
\begin{array}
[c]{ll}
k!\left\langle T,V\right\rangle , & k=\ell,\\
0, & k\neq\ell.
\end{array}
\right.
\]
This proves $\left(  \ref{eq:orthog}\right)  .$

\textbf{3}\emph{. }Fix $v\in\mathbb{R}^{d}$. The function $t\mapsto
E_{tv}\left(  z\right)  =e^{t\left\langle z,v\right\rangle -t^{2}\left\vert
v\right\vert ^{2}/2}$ has a Taylor series convergent for all $t\in\mathbb{R}$.
So the computation of $\left(  \ref{eq:hermite-co}\right)  $ shows, at $t=1$,
that
\[
E_{v}\left(  z\right)  =\sum_{k=0}^{\infty}\frac{\left\langle H_{k}\left(
z\right)  ,v^{\otimes k}\right\rangle }{k!}.
\]
For $Z\sim\gamma_{d},$ we verify that the series is convergent in
$L^{2}\left(  \gamma_{d}\right)  $ uniformly in $\left\vert v\right\vert
\leqslant1$. For integers $k,$ $K\geqslant0$ define
\begin{equation}
F_{k}\left(  v,Z\right)  =\frac{\left\langle H_{k}\left(  Z\right)
,v^{\otimes k}\right\rangle }{k!}\text{ and }P_{K}\left(  v,Z\right)
=\sum_{k=0}^{K}F_{k}\left(  v,Z\right)  . \label{eq:finite-terms}
\end{equation}
Taking $u=v,t=1$ in $\left(  \ref{eq:gen-hermite}\right)  $ and
differentiating in $h$ yield
\[
\mathbb{E}\left[  E_{v}\left(  Z\right)  F_{k}\left(  v,Z\right)  \right]
=\frac{1}{k!}\left.  \frac{d^{k}}{dh^{k}}\mathbb{E}\left[  E_{v}\left(
Z\right)  E_{hv}\left(  Z\right)  \right]  \right\vert _{h=0}=\frac{1}
{k!}\left.  \frac{d^{k}}{dh^{k}}e^{h\left\vert v\right\vert ^{2}}\right\vert
_{h=0}=\frac{\left\vert v\right\vert ^{2k}}{k!}.
\]
The differentiation is justified by dominated convergence. By the already
proved formula $\left(  \ref{eq:orthog}\right)  $,
\[
\mathbb{E}\left[  F_{k}\left(  v,Z\right)  F_{\ell}\left(  v,Z\right)
\right]  =1_{\left\{  k=\ell\right\}  }\frac{\left\vert v\right\vert ^{2k}
}{k!}.
\]
We also compute
\[
\mathbb{E}\left[  E_{v}\left(  Z\right)  ^{2}\right]  =e^{-\left\vert
v\right\vert ^{2}}\mathbb{E}e^{\left\langle 2v,Z\right\rangle }=e^{-\left\vert
v\right\vert ^{2}}e^{2\left\vert v\right\vert ^{2}}=e^{\left\vert v\right\vert
^{2}}.
\]
Hence, expanding the finite squared error we obtain
\begin{align*}
\mathbb{E}|E_{v}(Z)-P_{K}(v,Z)|^{2}  &  =\mathbb{E}E_{v}(Z)^{2}-2\sum
_{k=0}^{K}\mathbb{E}[E_{v}(Z)F_{k}(v,Z)]+\sum_{k,\ell=0}^{K}\mathbb{E}
[F_{k}(v,Z)F_{\ell}(v,Z)]\\
&  =e^{|v|^{2}}-2\sum_{k=0}^{K}\frac{|v|^{2k}}{k!}+\sum_{k=0}^{K}
\frac{|v|^{2k}}{k!}.
\end{align*}
Thus
\begin{equation}
\mathbb{E}\left\vert E_{v}\left(  Z\right)  -\sum_{k=0}^{K}\frac{\left\langle
H_{k}\left(  Z\right)  ,v^{\otimes k}\right\rangle }{k!}\right\vert
^{2}=e^{\left\vert v\right\vert ^{2}}-\sum_{k=0}^{K}\frac{\left\vert
v\right\vert ^{2k}}{k!}=\sum_{k>K}\frac{\left\vert v\right\vert ^{2k}}{k!}.
\label{eq:tail}
\end{equation}
For $\left\vert v\right\vert \leqslant1$, the right side is bounded by the
tail series
\[
\sum_{k>K}\frac{1}{k!}\leqslant\frac{1}{\left(  K+1\right)  !}\sum
_{r=0}^{\infty}\left(  K+2\right)  ^{-r}=\frac{K+2}{(K+1)(K+1)!}\rightarrow0.
\]
The inequality follows from $\left(  K+1+r\right)  !\geqslant\left(
K+1\right)  !\left(  K+2\right)  ^{r}$. For $v=0$, all positive degree terms
vanish and the remainder is zero. So uniformly for all $\left\vert
v\right\vert \leqslant1,$ the expansion $\left(  \ref{eq:hermite}\right)  $
holds in $L^{2}\left(  \gamma_{d}\right)  $.
\end{proof}

Equipped with the Hermite expansion we now proceed to the bound of
$\mathbb{E}R_{s,j}^{2}$.

\begin{lemma}
\label{lem:expansion}Let $n\geqslant1$ be an integer. Suppose that
$U_{1},\ldots,U_{n}\in\mathbb{R}^{d}$ are independent, centered and bounded.
Suppose that $S=\sum_{i}U_{i}$ has covariance matrix $L=\operatorname{Cov}
(S)=\operatorname{diag}(\lambda_{1},\ldots,\lambda_{d})$ with $\lambda_{i}>0$.
Let $(U_{i}^{\prime})_{i}$ be an independent copy of $(U_{i})_{i}$, and take
$Z\sim N(0,I_{d})$ independent of both families. With $Y_{s},$ $R_{s}$ from
$\left(  \ref{eq:score}\right)  $ and the coefficients from $\left(
\ref{eq:coeff}\right)  $-$\left(  \ref{eq:tensors}\right)  $, for each $s>0$
and $1\leqslant j\leqslant d$, it holds that
\begin{equation}
\mathbb{E}R_{s,j}^{2}\leqslant\mathbb{E}\left\vert \sum_{i}b_{i,j}\left(
s\right)  \right\vert ^{2}+\frac{1}{4}\sum_{k=1}^{\infty}\frac{\mathbb{E}
\left\Vert T_{j,k}\right\Vert ^{2}}{k!s^{k}}. \label{eq:hermite-bound}
\end{equation}

\end{lemma}

\begin{proof}
\textbf{1}. Fix $1\leqslant j\leqslant d$. Write $Y_{s}=S+\sqrt{s}Z.$ By the
definition $\left(  \ref{eq:score}\right)  $, the $j$-th component of $R_{s}$
in the eigenbasis of $L$ is $R_{s,j}=\mathbb{E}[S_{j}-\lambda_{j}Z_{j}
/\sqrt{s}\mid Y_{s}]$. Applying Jensen inequality we have
\[
\mathbb{E}R_{s,j}^{2}\leqslant\mathbb{E}\left(  S_{j}-\lambda_{j}Z_{j}
/\sqrt{s}\right)  ^{2}=\lambda_{j}+\lambda_{j}^{2}/s.
\]
The cross term vanishes since $S$ is centered and independent of $Z$. Thus
$R_{s,j}\in L^{2}$ for all $s>0$. The bound depends only on the covariance
$L$, whereas $\left(  \ref{eq:hermite-bound}\right)  $ must retain information
about each $d_{i}$. We therefore seek a some $V_{s,j}\in L^{2}$ so that it
preserves the desired information and
\begin{equation}
R_{s,j}=\mathbb{E}[V_{s,j}\mid Y_{s}]. \label{eq:equal-conditional}
\end{equation}
For any such choice, the contraction property will give $\mathbb{E}R_{s,j}
^{2}\leqslant\mathbb{E}V_{s,j}^{2}$. So in the following we will attempt to
find such $V_{s,j}$ with computable second moments. By the definition of
conditional expectation$,$ for $\left(  \ref{eq:equal-conditional}\right)  $
to hold, we must have, for any bounded Borel function $f$ on $\mathbb{R}^{d}$,
that
\[
\mathbb{E}[V_{s,j}f\left(  Y_{s}\right)  ]=\mathbb{E}[R_{s,j}f\left(
Y_{s}\right)  ]=\mathbb{E}[\left(  S_{j}-\lambda_{j}Z_{j}/\sqrt{s}\right)
f\left(  Y_{s}\right)  ].
\]
Therefore the focus is on $\mathbb{E}[S_{j}f\left(  Y_{s}\right)  ].$

\textbf{2}. Let
\[
\mathcal{U}=\sigma\left(  U_{1},\ldots,U_{n}\right)  ,\text{ }\mathcal{G}
=\sigma\left(  U_{1},\ldots,U_{n},Z\right)  ,\text{ }\mathcal{D}=\sigma\left(
U_{1},U_{1}^{\prime},\ldots,U_{n},U_{n}^{\prime}\right)  .
\]
Each $U_{i}^{\prime}$ is independent of $\mathcal{G}$, and $\mathbb{E}
U_{i}^{\prime}=0$. So
\[
\mathbb{E}\left[  d_{i,j}\mid\mathcal{G}\right]  =\mathbb{E}\left[
U_{i,j}^{\prime}\mid\mathcal{G}\right]  -U_{i,j}=\mathbb{E}U_{i,j}^{\prime
}-U_{i,j}=-U_{i,j}.
\]
Since $Y_{s}$ is $\mathcal{G}$-measurable and all integrands are bounded, we
get
\[
\mathbb{E}\left[  d_{i,j}f\left(  Y_{s}\right)  \right]  =\mathbb{E}\left[
f\left(  Y_{s}\right)  \mathbb{E}\left[  d_{i,j}\mid\mathcal{G}\right]
\right]  =-\mathbb{E}\left[  U_{i,j}f(Y_{s})\right]  .
\]
Using $S_{j}=\sum_{i}U_{i,j}$ yields
\begin{equation}
\mathbb{E}\left[  S_{j}f\left(  Y_{s}\right)  \right]  =-\sum_{i}
\mathbb{E}\left[  d_{i,j}f\left(  Y_{s}\right)  \right]  .
\label{eq:center-swap}
\end{equation}
It is sufficient to focus on $\mathbb{E}\left[  d_{i,j}f\left(  Y_{s}\right)
\right]  .$ We now split the expectation into two parts
\[
\mathbb{E}\left[  d_{i,j}f\left(  Y_{s}\right)  \right]  =\mathbb{E}\left[
d_{i,j}1_{\left\{  r_{i}\leqslant s\right\}  }f\left(  Y_{s}\right)  \right]
+\mathbb{E}\left[  d_{i,j}1_{\left\{  r_{i}>s\right\}  }f\left(  Y_{s}\right)
\right]  .
\]
The reason for the choice of the truncation threshold will become clear in the
computation. We address the two terms in the sum separately. The second term
has
\[
\mathbb{E}\left[  d_{i,j}1_{\left\{  r_{i}>s\right\}  }f\left(  Y_{s}\right)
\right]  =\mathbb{E}\left[  \mathbb{E}\left[  d_{i,j}1_{\left\{
r_{i}>s\right\}  }f\left(  Y_{s}\right)  \mid\mathcal{G}\right]  \right]
=\mathbb{E}\left[  \mathbb{E}\left[  d_{i,j}1_{\left\{  r_{i}>s\right\}  }
\mid\mathcal{G}\right]  f\left(  Y_{s}\right)  \right]  .
\]
The integrand $d_{i,j}1_{\left\{  r_{i}>s\right\}  }$ is a function of
$U_{i},U_{i}^{\prime}$. Except $U_{i},$ all other information of $\mathcal{G}$
are independent of $U_{i},U_{i}^{\prime}.$ Hence
\[
\mathbb{E}\left[  d_{i,j}1_{\left\{  r_{i}>s\right\}  }\mid\mathcal{G}\right]
=\mathbb{E}\left[  d_{i,j}1_{\left\{  r_{i}>s\right\}  }\mid U_{i}\right]
=b_{i,j}(s).
\]
Thus
\begin{equation}
\mathbb{E}\left[  d_{i,j}1_{\left\{  r_{i}>s\right\}  }f\left(  Y_{s}\right)
\right]  =\mathbb{E}\left[  b_{i,j}(s)f\left(  Y_{s}\right)  \right]  .
\label{eq:resample}
\end{equation}
We next turn to the first term $\mathbb{E}\left[  d_{i,j}1_{\left\{
r_{i}\leqslant s\right\}  }f\left(  Y_{s}\right)  \right]  .$ The expectation
is determined by the joint distribution of $U_{i}^{\prime}$, $U_{1}
,\ldots,U_{n},Z,\,$which is the product of individual distribution. In
particular, swapping $U_{i}^{\prime}$ and $U_{i}$ does not change the joint
distribution, and thus does not change the expectation. On the other hand, the
swapping smounts to replacing $d_{i,j},$ $r_{i},$ $Y_{s}$ respectively by
$-d_{i,j},$ $r_{i},$ $Y_{s}-U_{i}+U_{i}^{\prime}.$ Therefore we have
\[
\mathbb{E}\left[  d_{i,j}1_{\left\{  r_{i}\leqslant s\right\}  }f\left(
Y_{s}\right)  \right]  =-\mathbb{E}\left[  d_{i,j}1_{\left\{  r_{i}\leqslant
s\right\}  }f\left(  Y_{s}+d_{i}\right)  \right]  .
\]
It follows that
\begin{equation}
\mathbb{E}\left[  d_{i,j}1_{\left\{  r_{i}\leqslant s\right\}  }f\left(
Y_{s}\right)  \right]  =-\tfrac{1}{2}\mathbb{E}\left[  d_{i,j}1_{\left\{
r_{i}\leqslant s\right\}  }\left(  f\left(  Y_{s}+d_{i}\right)  -f\left(
Y_{s}\right)  \right)  \right]  . \label{eq:shift}
\end{equation}
The next step is compute the expectation on the right side.

\textbf{3}. Using the tower property and that $d_{i,j}1_{\left\{
r_{i}\leqslant s\right\}  }$ is $\mathcal{D}$-measurable, we write
\begin{equation}
\mathbb{E}\left[  d_{i,j}1_{\left\{  r_{i}\leqslant s\right\}  }\left(
f\left(  Y_{s}+d_{i}\right)  -f\left(  Y_{s}\right)  \right)  \right]
=\mathbb{E}\left[  d_{i,j}1_{\left\{  r_{i}\leqslant s\right\}  }
\mathbb{E}\left[  f\left(  Y_{s}+d_{i}\right)  -f\left(  Y_{s}\right)
\mid\mathcal{D}\right]  \right]  . \label{eq:copy}
\end{equation}
Recall that $Y_{s}=S+\sqrt{s}Z.$ Conditioned on $\mathcal{D}$, the random
vectors $S,d_{i}$ are fixed and the independent $Z$ still has distribution
$\gamma_{d}=N(0,I_{d})$. For $v,$ $z\in\mathbb{R}^{d}$, define
\[
E_{v}\left(  z\right)  =\exp\left(  \left\langle z,v\right\rangle
-|v|^{2}/2\right)  .
\]
Setting $v=d_{i}/\sqrt{s}$, by direct computation we have that
\begin{align*}
\left(  2\pi\right)  ^{-d/2}\int f(S+\sqrt{s}z+d_{i})e^{-\left\vert
z\right\vert ^{2}/2}dz  &  =\left(  2\pi\right)  ^{-d/2}\int f(S+\sqrt
{s}\left(  z+v\right)  )e^{-\left\vert z\right\vert ^{2}/2}dz\\
&  =\left(  2\pi\right)  ^{-d/2}\int f(S+\sqrt{s}w)e^{-\left\vert
w-v\right\vert ^{2}/2}dw\\
&  =\left(  2\pi\right)  ^{-d/2}\int f(S+\sqrt{s}w)E_{v}(w)e^{-\left\vert
w\right\vert ^{2}/2}dw.
\end{align*}
Here $w=z+v$ and the last equality uses the equation
\[
\exp\left(  -\frac{\left\vert w-v\right\vert ^{2}}{2}+\frac{\left\vert
w\right\vert ^{2}}{2}\right)  =\exp\left(  \left\langle w,v\right\rangle
-\left\vert v\right\vert ^{2}/2\right)  =E_{v}\left(  w\right)  .
\]
Hence we have
\[
\int f(S+\sqrt{s}z+d_{i})d\gamma_{d}\left(  z\right)  =\int f\left(
S+\sqrt{s}z\right)  E_{v}\left(  z\right)  d\gamma_{d}\left(  z\right)  .
\]
With these computations we get
\[
\mathbb{E}\left[  f\left(  Y_{s}+d_{i}\right)  -f\left(  Y_{s}\right)
\mid\mathcal{D}\right]  =\int f\left(  S+\sqrt{s}z\right)  \left[
E_{v}\left(  z\right)  -1\right]  d\gamma_{d}\left(  z\right)  .
\]
Substituting it into $\left(  \ref{eq:copy}\right)  $ yields
\begin{align}
\mathbb{E}\left[  d_{i,j}1_{\left\{  r_{i}\leqslant s\right\}  }\left(
f\left(  Y_{s}+d_{i}\right)  -f\left(  Y_{s}\right)  \right)  \right]   &
=\mathbb{E}\left[  d_{i,j}1_{\left\{  r_{i}\leqslant s\right\}  }
\mathbb{E}\left[  f\left(  Y_{s}\right)  \left(  E_{d_{i}/\sqrt{s}}\left(
Z\right)  -1\right)  \mid\mathcal{D}\right]  \right] \label{eq:shift-E}\\
&  =\mathbb{E}\left[  d_{i,j}1_{\left\{  r_{i}\leqslant s\right\}  }f\left(
Y_{s}\right)  \left(  E_{d_{i}/\sqrt{s}}\left(  Z\right)  -1\right)  \right]
.\nonumber
\end{align}
Note that the absolute value of the expectation on the right side is at most
$2\sqrt{s}\sup|f|$, because $\left\vert d_{i,j}\right\vert \leqslant\sqrt{s}$
on the truncated event $\left\{  r_{i}\leqslant s\right\}  $ and
$\mathbb{E}\left[  \left\vert E_{d_{i}/\sqrt{s}}\left(  Z\right)
-1\right\vert \mid\mathcal{D}\right]  \leqslant2$.

\textbf{4}. To proceed further with $\left(  \ref{eq:shift-E}\right)  ,$ we
invoke Hermite expanasion proved in Lemma \ref{lm:hermite-expansion} to write
\[
E_{v}\left(  Z\right)  =\sum_{k=0}^{\infty}\frac{\left\langle H_{k}\left(
Z\right)  ,v^{\otimes k}\right\rangle }{k!},
\]
where $H_{k}\left(  z\right)  $ is defined in $\left(
\ref{eq:hermite-differential}\right)  .$ The series converges in $L^{2}\left(
\gamma_{d}\right)  $ uniformly in $\left\vert v\right\vert \leqslant1$. Keep
the definitions of $F_{k}\left(  v,Z\right)  $ and $P_{K}\left(  v,Z\right)  $
from $\left(  \ref{eq:finite-terms}\right)  .$ For fixed $s,j$, define
\[
\xi_{i}=d_{i,j}1_{\left\{  r_{i}\leqslant s\right\}  }\left(  E_{d_{i}
/\sqrt{s}}\left(  Z\right)  -1\right)  ,\text{ }\xi_{i}^{\left(  K\right)
}=d_{i,j}1_{\left\{  r_{i}\leqslant s\right\}  }\sum_{k=1}^{K}F_{k}\left(
d_{i}/\sqrt{s},Z\right)  .
\]
Due to the convergent property of $E_{v}\left(  Z\right)  $ with $\left\vert
v\right\vert \leqslant1,$ $\xi_{i}$ is a well-defined element in $L^{2}\left(
\gamma_{d}\right)  $. The constant term of $E_{v}\left(  z\right)  $ is
$F_{0}=1$, so this series starts at $k=1$. Using the independence of $Z$ from
$\mathcal{D}$ and $\left(  \ref{eq:gen-hermite}\right)  $ we have
\[
\mathbb{E}\left[  \xi_{i}^{2}\mid\mathcal{D}\right]  =d_{i,j}^{2}1_{\left\{
r_{i}\leqslant s\right\}  }\left(  e^{r_{i}/s}-2+1\right)  \leqslant s\left(
e-1\right)  .
\]
Thus each $\xi_{i}$ is square integrable. $\tau_{K}=\sum_{k>K}1/k!$.
Conditional on $\mathcal{D}$, $d_{i}$ is fixed, so $\left(  \ref{eq:tail}
\right)  $ applies and yields
\[
\mathbb{E}\left[  \left.  \left\vert \xi_{i}-\xi_{i}^{(K)}\right\vert
^{2}\right\vert \mathcal{D}\right]  =d_{i,j}^{2}1_{\left\{  r_{i}\leqslant
s\right\}  }\sum_{k>K}\frac{\left(  r_{i}/s\right)  ^{k}}{k!}\leqslant
s\tau_{K}.
\]
Since $\left\vert \sum_{i}a_{i}\right\vert ^{2}\leqslant n\sum_{i}\left\vert
a_{i}\right\vert ^{2},$ we get
\[
\left\Vert \sum_{i}(\xi_{i}-\xi_{i}^{(K)})\right\Vert _{L^{2}}^{2}\leqslant
n\sum_{i}\mathbb{E}|\xi_{i}-\xi_{i}^{(K)}|^{2}\leqslant n^{2}s\tau
_{K}\rightarrow0,\text{ as }K\rightarrow\infty.
\]
Since conditional expectation is a contraction \cite[Theorem 4.1.11]
{durrett2019probability}, we obtain that
\begin{equation}
\left\Vert \mathbb{E}\left[  \left.  \sum_{i}\xi_{i}-\sum_{i}\xi_{i}
^{(K)}\right\vert \mathcal{G}\right]  \right\Vert _{L^{2}}^{2}\leqslant
n^{2}s\tau_{K}\rightarrow0,\text{ as }K\rightarrow\infty.
\label{eq:contraction}
\end{equation}

\textbf{5}. To expand the conditional expectation in $\left(
\ref{eq:contraction}\right)  $ and identify the coefficients, write $\nu_{i}$
for the distribution of $U_{i}$. By independence, the conditional distribution
of $U_{i}^{\prime}$ given $\mathcal{G}$ equals $\nu_{i}$. Hence
\[
\mathbb{E}\left[  d_{i,j}d_{i}^{\otimes k}1_{\left\{  r_{i}\leqslant
s\right\}  }\mid\mathcal{G}\right]  =\int\left(  u_{j}-U_{i,j}\right)  \left(
u-U_{i}\right)  ^{\otimes k}1_{\{\left\vert u-U_{i}\right\vert ^{2}\leqslant
s\}}\nu_{i}\left(  du\right)  =C_{i,j,k}\left(  s\right)  ,
\]
\[
\mathbb{E}[d_{i,j}1_{\{r_{i}>s\}}\mid\mathcal{G}]=\int\left(  u_{j}
-U_{i,j}\right)  1_{\{|u-U_{i}|^{2}>s\}}\nu_{i}\left(  du\right)
=b_{i,j}\left(  s\right)  .
\]
These are the $\sigma\left(  U_{i}\right)  $-measurable versions in $\left(
\ref{eq:coeff}\right)  $ due to independence of the copy $(U_{i}^{\prime}
)_{i}$ from $\mathcal{G}$. From the definition $\left(  \ref{eq:coeff}\right)
$, we see that permuting the $k$ indices of $C_{i,j,k}$ leaves it unchanged,
thus $C_{i,j,k}$ is symmetric in its $k$ tensor indices. Hence we can compute
the bounds,
\[
\left\Vert C_{i,j,k}\left(  s\right)  \right\Vert \leqslant\mathbb{E}\left[
\left\vert d_{i,j}\right\vert \left\vert d_{i}\right\vert ^{k}1_{\left\{
r_{i}\leqslant s\right\}  }\mid U_{i}\right]  \leqslant s^{\left(  k+1\right)
/2}.
\]
Since $\left(  d_{i}/\sqrt{s}\right)  ^{\otimes k}=s^{-k/2}d_{i}^{\otimes k}$,
pulling the $\mathcal{G}$-measurable entries of $H_{k}\left(  Z\right)  $ out
of the conditional sums we get
\begin{align*}
\mathbb{E}\left[  \left.  \sum_{i}\xi_{i}^{(K)}\right\vert \mathcal{G}
\,\right]   &  =\sum_{i}\mathbb{E}\left[  \left.  d_{i,j}1_{\left\{
r_{i}\leqslant s\right\}  }\sum_{k=1}^{K}\frac{\left\langle H_{k}\left(
Z\right)  ,\left(  d_{i}/\sqrt{s}\right)  ^{\otimes k}\right\rangle }
{k!}\right\vert \mathcal{G}\right] \\
&  =\sum_{k=1}^{K}\frac{1}{k!s^{k/2}}\sum_{\ell_{1},\ldots,\ell_{k}}\left(
H_{k}\left(  Z\right)  \right)  _{\ell_{1},\ldots,\ell_{k}}\sum_{i}\left(
C_{i,j,k}\left(  s\right)  \right)  _{\ell_{1},\ldots,\ell_{k}}=\sum_{k=1}
^{K}\frac{\left\langle H_{k}\left(  Z\right)  ,\sum_{i}C_{i,j,k}\left(
s\right)  \right\rangle }{k!s^{k/2}}.
\end{align*}
All terms are squared integrable for fixed $s>0$, because the coefficients
$\sum_{i}C_{i,j,k}\left(  s\right)  $ are bounded and the Hermite polynomials
$H_{k}\left(  Z\right)  $ have finite Gaussian moments. The contraction
inequality $\left(  \ref{eq:contraction}\right)  $ then justifies the limit in
$L^{2},$
\begin{equation}
\lim_{K\rightarrow\infty}\sum_{k=1}^{K}\frac{\left\langle H_{k}\left(
Z\right)  ,\sum_{i}C_{i,j,k}\left(  s\right)  \right\rangle }{k!s^{k/2}}
=\lim_{K\rightarrow\infty}\mathbb{E}\left[  \left.  \sum_{i}\xi_{i}
^{(K)}\right\vert \mathcal{G}\,\right]  =\mathbb{E}\left[  \left.  \sum_{i}
\xi_{i}\right\vert \mathcal{G}\,\right]  . \label{eq:K-limit}
\end{equation}
Therefore using $\left(  \ref{eq:shift-E}\right)  ,$
\begin{equation}
\mathbb{E}\left[  d_{i,j}1_{\left\{  r_{i}\leqslant s\right\}  }\left(
f\left(  Y_{s}+d_{i}\right)  -f\left(  Y_{s}\right)  \right)  \right]
=\mathbb{E}\left[  f\left(  Y_{s}\right)  \xi_{i}\right]  .
\label{eq:shift-Xi}
\end{equation}

\textbf{5}\emph{. }Combining the previous computations suggests the following
candidate for $V_{s,j}$,
\[
V_{s,j}=-\sum_{i}b_{i,j}\left(  s\right)  -\frac{\lambda_{j}Z_{j}}{\sqrt{s}
}+\frac{1}{2}\mathbb{E}\left[  \left.  \sum_{i}\xi_{i}\right\vert
\mathcal{G}\right]  .
\]
It is square integrable. Indeed, the last term is square integrable by the
early estimates and $\mathbb{E}\left[  b_{i,j}\left(  s\right)  ^{2}\right]
\leqslant\mathbb{E}\left[  d_{i,j}^{2}1_{\left\{  r_{i}>s\right\}  }\right]
<\infty$. For $K\geqslant1,$ define its finite $K$-approximants by
\[
V_{s,j}^{\left(  K\right)  }=-\sum_{i}b_{i,j}(s)-\frac{\lambda_{j}Z_{j}}
{\sqrt{s}}+\frac{1}{2}\mathbb{E}\!\left[  \left.  \sum_{i}\xi_{i}
^{(K)}\right\vert \mathcal{G}\right]  .
\]
Since $H_{1}(Z)=Z$, the sum of the terms of $Z$-order one is
\[
-\frac{\lambda_{j}Z_{j}}{\sqrt{s}}+\frac{\langle H_{1}(Z),\sum_{i}
C_{i,j,1}(s)\rangle}{2\sqrt{s}}=\frac{\langle Z,\sum_{i}C_{i,j,1}
(s)-2\lambda_{j}e_{j}\rangle}{2\sqrt{s}}=\frac{\langle H_{1}(Z),T_{j,1}
\rangle}{2\sqrt{s}},
\]
where $T_{j,1}$ is defined in $\left(  \ref{eq:tensors}\right)  .$ For
$k\geqslant2$, $T_{j,k}=\sum_{i}C_{i,j,k}(s)$. Consequently
\[
V_{s,j}^{\left(  K\right)  }=-\sum_{i}b_{i,j}\left(  s\right)  +\frac{1}
{2}\sum_{k=1}^{K}\frac{\left\langle H_{k}\left(  Z\right)  ,T_{j,k}
\right\rangle }{k!s^{k/2}},
\]
By $\left(  \ref{eq:contraction}\right)  $,
\[
\left\Vert V_{s,j}^{(K)}-V_{s,j}\right\Vert _{L^{2}}^{2}\leqslant\frac{n^{2}
s}{4}\tau_{K}\rightarrow0.
\]
Therefore we have derived, with convergence in $L^{2}$,
\begin{equation}
V_{s,j}=-\sum_{i}b_{i,j}\left(  s\right)  +\frac{1}{2}\sum_{k=1}^{\infty}
\frac{\left\langle H_{k}\left(  Z\right)  ,T_{j,k}\right\rangle }{k!s^{k/2}}.
\label{eq:score-proxy}
\end{equation}

\textbf{6}\emph{. }Combing $\left(  \ref{eq:center-swap}\right)  $,$\left(
\ref{eq:resample}\right)  $ and $\left(  \ref{eq:shift}\right)  $ with
$\left(  \ref{eq:shift-Xi}\right)  $ therefore yields
\begin{align*}
\mathbb{E}[S_{j}f(Y_{s})]  &  =-\sum_{i}\mathbb{E}[b_{i,j}(s)f(Y_{s}
)]+\frac{1}{2}\sum_{i}\mathbb{E}[d_{i,j}1_{\{r_{i}\leqslant s\}}(f(Y_{s}
+d_{i})-f(Y_{s}))]\\
&  =-\sum_{i}\mathbb{E}\left[  b_{i,j}\left(  s\right)  f\left(  Y_{s}\right)
\right]  +\frac{1}{2}\mathbb{E}\left[  f\left(  Y_{s}\right)  \sum_{i}\xi
_{i}\right] \\
&  =\mathbb{E}\left[  f\left(  Y_{s}\right)  \left(  -\sum_{i}b_{i,j}
(s)+\frac{1}{2}\mathbb{E}\left[  \left.  \sum_{i}\xi_{i}\right\vert
\mathcal{G}\right]  \right)  \right] \\
&  =\mathbb{E}\left[  f\left(  Y_{s}\right)  \left(  V_{s,j}+\frac{\lambda
_{j}Z_{j}}{\sqrt{s}}\right)  \right]  .
\end{align*}
Thus, using the definition of $R_{s}$ and the diagonal convariance $L$,
\[
\mathbb{E}\left[  V_{s,j}f\left(  Y_{s}\right)  \right]  =\mathbb{E}\left[
\left(  S_{j}-\lambda_{j}Z_{j}/\sqrt{s}\right)  f\left(  Y_{s}\right)
\right]  =\mathbb{E}\left[  \mathbb{E}\left[  S_{j}-\lambda_{j}Z_{j}/\sqrt
{s}\mid Y_{s}\right]  f\left(  Y_{s}\right)  \right]  =\mathbb{E}\left[
R_{s,j}f\left(  Y_{s}\right)  \right]  .
\]
Noting that both $V_{s,j}$, $R_{s,j}$ are integrable, $R_{s,j}$ is
$\sigma\left(  Y_{s}\right)  $-measurable and that $f$ is an arbitrary bounded
Borel function, we therefore have
\[
\mathbb{E}[V_{s,j}\mid Y_{s}]=R_{s,j},\text{ }\mathbb{P}\text{-}a.e.
\]

\textbf{7}. Conditioning on $\mathcal{U}$ fixes all coefficients in $\left(
\ref{eq:score-proxy}\right)  $ and leaves $Z$ standard Gaussian. Write
$\zeta_{j,k}=\left\langle H_{k}\left(  Z\right)  ,T_{j,k}\right\rangle
/\left(  k!s^{k/2}\right)  $. Conditioning on $\mathcal{U}$, $Z$ is a standard
Gaussian. Since $C_{i,j,k}$ is symmetric in its $k$ tensor indices, we see
from the definition $\left(  \ref{eq:tensors}\right)  $ that $T_{j,k}$ is
symmetric in its $k$ tensor indices. Hence using $\left(  \ref{eq:orthog}
\right)  $ we have, for $k,\ell\geqslant1$, that
\[
\mathbb{E}\left[  \zeta_{j,k}\mid\mathcal{U}\right]  =0,
\]
and
\[
\mathbb{E}\left[  \zeta_{j,k}\zeta_{j,\ell}\mid\mathcal{U}\right]  =
\begin{cases}
\frac{k!\left\Vert T_{j,k}\right\Vert ^{2}}{\left(  k!\right)  ^{2}s^{k}
}=\frac{\left\Vert T_{j,k}\right\Vert ^{2}}{k!s^{k}}, & k=\ell,\\
0, & k\neq\ell.
\end{cases}
\]
It follows that
\begin{align*}
\mathbb{E}\left[  \left.  \left(  V_{s,j}^{\left(  K\right)  }\right)
^{2}\right\vert \mathcal{U}\right]   &  =\left\vert \sum_{i}b_{i,j}\left(
s\right)  \right\vert ^{2}-\left(  \sum_{i}b_{i,j}\left(  s\right)  \right)
\sum_{k=1}^{K}\mathbb{E}\left[  \zeta_{j,k}\mid\mathcal{U}\right]  +\frac
{1}{4}\sum_{k,\ell=1}^{K}\mathbb{E}\left[  \zeta_{j,k}\zeta_{j,\ell}
\mid\mathcal{U}\right] \\
&  =\left\vert \sum_{i}b_{i,j}\left(  s\right)  \right\vert ^{2}+\frac{1}
{4}\sum_{k=1}^{K}\frac{\left\Vert T_{j,k}\right\Vert ^{2}}{k!s^{k}}.
\end{align*}
The $L^{2}$ convergence proved earlier implies convergence of second moments,
\[
\left\vert \mathbb{E}(V_{s,j}^{(K)})^{2}-\mathbb{E}V_{s,j}^{2}\right\vert
\leqslant\left\Vert V_{s,j}^{(K)}-V_{s,j}\right\Vert _{L^{2}}\left(
\left\Vert V_{s,j}^{(K)}\right\Vert _{L^{2}}+\left\Vert V_{s,j}\right\Vert
_{L^{2}}\right)  \rightarrow0.
\]
Therefore
\begin{align*}
\mathbb{E}R_{s,j}^{2}  &  =\mathbb{E}\left\vert \mathbb{E}\left[  V_{s,j}\mid
Y_{s}\right]  \right\vert ^{2}\leqslant\mathbb{E}V_{s,j}^{2}\\
&  =\lim_{K\rightarrow\infty}\left(  \mathbb{E}\left\vert \sum_{i}
b_{i,j}\left(  s\right)  \right\vert ^{2}+\frac{1}{4}\sum_{k=1}^{K}
\frac{\mathbb{E}\left\Vert T_{j,k}\right\Vert ^{2}}{k!s^{k}}\right)
=\mathbb{E}\left\vert \sum_{i}b_{i,j}(s)\right\vert ^{2}+\frac{1}{4}\sum
_{k=1}^{\infty}\frac{\mathbb{E}\Vert T_{j,k}\Vert^{2}}{k!s^{k}}.
\end{align*}
This completes the proof for $\left(  \ref{eq:hermite-bound}\right)  $.
\end{proof}

\section{Integration of the expansion}

\label{sec:integrate}

\begin{proof}
[Proof of Theorem~\ref{thm:weighted} (the bounded case)]\textbf{1}. First
suppose the summands are bounded. Rotate them as in
Section~\ref{sec:observation}. We use the diagonal form $\left(
\ref{eq:identity}\right)  $ of Proposition~\ref{prop:observation} in precisely
this eigenbasis, so that the factor $w_{j}(s)$ multiplies the squared residual
of the same coordinate $j$ in Lemma~\ref{lem:expansion}. Write
\[
M=\sum_{i}\mathbb{E}r_{i}^{2}=\sum_{i}\mathbb{E}|d_{i}|^{4}.
\]
Multiply $\left(  \ref{eq:hermite-bound}\right)  $ by $w_{j}(s)$, sum over
$j$, and integrate in $s$. By $\left(  \ref{eq:identity}\right)  $ this bounds
the desired cost. We separate the variance and the squared mean of each tensor
using
\begin{equation}
\mathbb{E}\left\Vert T_{j,k}\right\Vert ^{2}=\mathbb{E}\left\Vert
T_{j,k}-\mathbb{E}T_{j,k}\right\Vert ^{2}+\left\Vert \mathbb{E}T_{j,k}
\right\Vert ^{2}. \label{eq:split-variance}
\end{equation}
Where a formula contains a negative power of $r_{i}$, we restrict the
integration to $\{r_{i}>0\}$ and set the entire integrand to zero on
$\{r_{i}=0\}$, where the associated tensor vanishes.

\textbf{2}. Bounds for the terms $b_{i,j}(s).$ Since the distribution of
$d_{i}$ is symmetric, so
\[
\mathbb{E}\left[  b_{i,j}\left(  s\right)  \right]  =\mathbb{E}\left[
d_{i,j}1_{\{r_{i}>s\}}\right]  =0.
\]
The $b_{i,j}\left(  s\right)  $ for fixed $j,s$ are independent over $i$. So
for $i\neq h$,
\[
\mathbb{E}\left[  b_{i,j}\left(  s\right)  b_{h,j}\left(  s\right)  \right]
=\left(  \mathbb{E}b_{i,j}\left(  s\right)  \right)  \left(  \mathbb{E}
b_{h,j}\left(  s\right)  \right)  =0.
\]
It follows that
\[
\mathbb{E}\left\vert \sum_{i}b_{i,j}\left(  s\right)  \right\vert ^{2}
=\sum_{i}\mathbb{E}b_{i,j}\left(  s\right)  ^{2}.
\]
Conditional Jensen in $\left(  \ref{eq:coeff}\right)  $ and the supremum in
$\left(  \ref{eq:weights}\right)  $ give
\begin{align}
\int_{0}^{\infty}\sum_{j}w_{j}(s)\mathbb{E}\left\vert \sum_{i}b_{i,j}
(s)\right\vert ^{2}ds  &  \leqslant\sum_{i}\int_{0}^{\infty}\sum_{j}
w_{j}(s)\mathbb{E}[d_{i,j}^{2}1_{\{s<r_{i}\}}]ds\label{eq:zero}\\
&  \leqslant\frac{8}{27}\sum_{i}\mathbb{E}\left[  r_{i}\int_{0}^{r_{i}
}ds\right]  =\frac{8}{27}M.\nonumber
\end{align}

\textbf{3}\emph{. }Bounds for the terms $\mathbb{E}\Vert T_{j,k}
-\mathbb{E}T_{j,k}\Vert^{2}.$ For each $i,$ write $C_{i}^{\prime}
=C_{i,j,k}(s)-\mathbb{E}C_{i,j,k}(s)$. For $i\neq h$, the independence of
$U_{i},U_{h}$ and $\left(  \ref{eq:coeff}\right)  $ yield $\mathbb{E}
\left\langle C_{i}^{\prime},C_{h}^{\prime}\right\rangle =\langle
\mathbb{E}C_{i}^{\prime},\mathbb{E}C_{h}^{\prime}\rangle=0$. Consequently for
all $k\geqslant1$,
\begin{align}
\mathbb{E}\Vert T_{j,k}-\mathbb{E}T_{j,k}\Vert^{2}  &  =\sum_{i}
\mathbb{E}\Vert C_{i,j,k}(s)-\mathbb{E}C_{i,j,k}(s)\Vert^{2}
\label{eq:centered}\\
&  \leqslant\sum_{i}\mathbb{E}\Vert C_{i,j,k}(s)\Vert^{2}\leqslant\sum
_{i}\mathbb{E}[d_{i,j}^{2}r_{i}^{k}1_{\{r_{i}\leqslant s\}}].\nonumber
\end{align}
In the last inequality we used Jensen inequality and $\left(
\ref{eq:tensor-power}\right)  $ to write
\[
\left\Vert C_{i,j,k}(s)\right\Vert ^{2}\leqslant\mathbb{E}\left[  d_{i,j}
^{2}\left\vert d_{i}\right\vert ^{2k}1_{\{r_{i}\leqslant s\}}\mid
U_{i}\right]  .
\]
Define
\[
\mathcal{V}_{k}=\sum_{j}\int_{0}^{\infty}\frac{w_{j}(s)}{4k!s^{k}}
\mathbb{E}\Vert T_{j,k}-\mathbb{E}T_{j,k}\Vert^{2}ds,\text{ }k\geqslant1.
\]
For $k=1$, by $\left(  \ref{eq:weights}\right)  $ we have, for any value
$r\geqslant0$, that
\[
\int_{0}^{\infty}\frac{w_{j}(s)}{s}1_{\{r\leqslant s\}}ds\leqslant1.
\]
Taking $r=r_{i}$ in the inequality and that $\sum_{j}d_{i,j}^{2}r_{i}
=|d_{i}|^{2}r_{i}=r_{i}^{2}$ yield
\begin{equation}
\mathcal{V}_{1}\leqslant\frac{1}{4}\sum_{i,j}\mathbb{E}\left[  d_{i,j}
^{2}r_{i}\int_{r_{i}}^{\infty}\frac{w_{j}(s)}{s}\,ds\right]  \leqslant\frac
{1}{4}\sum_{i}\mathbb{E}r_{i}^{2}=\frac{1}{4}M. \label{eq:c1}
\end{equation}
The first inequality also uses Tonelli's theorem and the bound $\left(
\ref{eq:centered}\right)  .$ For $k\geqslant2$, using $w_{j}\leqslant8/27$
(see $\left(  \ref{eq:weights}\right)  $) and $\sum_{j}d_{i,j}^{2}=r_{i}$ in
the integral of $\left(  \ref{eq:centered}\right)  $ we obtain
\begin{equation}
\mathcal{V}_{k}\leqslant\frac{2}{27k!}\sum_{i}\mathbb{E}\left[  1_{\{r_{i}
>0\}}r_{i}^{k+1}\int_{r_{i}}^{\infty}s^{-k}\,ds\right]  =\frac{2M}{27k!(k-1)},
\label{eq:ck}
\end{equation}
since the integral equals $r_{i}^{1-k}/(k-1)$.

\textbf{4}\emph{. }Bounds for the terms $\left\Vert \mathbb{E}T_{j,k}
\right\Vert ^{2}.$ For tensors $V_{i}$ of a common finite shape that are
square integrable on the respective spaces of $(U_{i},U_{i}^{\prime})$, we
claim that
\begin{equation}
\sum_{j=1}^{d}\frac{1}{2\lambda_{j}}\left\Vert \sum_{i}\mathbb{E}[d_{i,j}
V_{i}]\right\Vert ^{2}\leqslant\sum_{i}\mathbb{E}\Vert V_{i}\Vert^{2}.
\label{eq:bessel}
\end{equation}
Indeed, let $\mathcal{H}=\bigoplus_{i}L^{2}(\mathcal{L}(U_{i},U_{i}^{\prime
}))$. By independence of the copies and centering we have
\[
\sum_{i}\mathbb{E}\left[  d_{i,j}d_{i,\ell}\right]  =2\sum_{i}
\operatorname{Cov}\left(  U_{i}\right)  _{j,\ell}=2\lambda_{j}1_{\{j=\ell\}}.
\]
Thus the $(d_{i,j}/\sqrt{2\lambda_{j}})_{i}$ are orthonormal. For each tensor
index $\alpha$, applying Bessel's inequality to $f^{\alpha}=((V_{i})_{\alpha
})_{i}$ and summing over $\alpha$ yield
\[
\sum_{\alpha,j}|\langle f^{\alpha},v_{j}\rangle|^{2}=\sum_{\alpha,j}
\frac{\left\vert \sum_{i}\mathbb{E}[d_{i,j}(V_{i})_{\alpha}]\right\vert ^{2}
}{2\lambda_{j}}\leqslant\sum_{\alpha}\Vert f^{\alpha}\Vert_{\mathcal{H}}
^{2}=\sum_{i}\mathbb{E}\Vert V_{i}\Vert^{2}.
\]
Write
\[
\mathcal{M}_{k}=\sum_{j}\int_{0}^{\infty}\frac{w_{j}(s)}{4k!s^{k}}
\Vert\mathbb{E}T_{j,k}\Vert^{2}ds.
\]
Since $\sum_{i}\mathbb{E}[d_{i,j}d_{i}]=2\lambda_{j}e_{j}$, by $\left(
\ref{eq:tensors}\right)  $ we have
\[
\mathbb{E}T_{j,1}=-\sum_{i}\mathbb{E}[d_{i,j}d_{i}1_{\{r_{i}>s\}}].
\]
Let $V_{i}=d_{i}1_{\{r_{i}>s\}}$ in $\left(  \ref{eq:bessel}\right)  $. These
vectors are square integrable since $\mathbb{E}\Vert V_{i}\Vert^{2}
=\mathbb{E}[r_{i}1_{\{r_{i}>s\}}]\leqslant\mathbb{E}r_{i}<\infty$. By the last
bound in $\left(  \ref{eq:weights}\right)  $,
\[
\sum_{j}\frac{w_{j}(s)}{4s}\Vert\mathbb{E}T_{j,1}\Vert^{2}\leqslant\sum
_{j}\frac{1}{2\lambda_{j}}\left\Vert \sum_{i}\mathbb{E}[d_{i,j}V_{i}
]\right\Vert ^{2}\leqslant\sum_{i}\mathbb{E}[r_{i}1_{\{s<r_{i}\}}].
\]
It follows that
\begin{equation}
\mathcal{M}_{1}\leqslant\sum_{i}\mathbb{E}\left[  r_{i}\int_{0}^{r_{i}
}ds\right]  =M. \label{eq:b1}
\end{equation}
If $k=2$, $d_{i,j}d_{i}^{\otimes2}$ is odd under $d_{i}\mapsto-d_{i}$ and
$1_{\{r_{i}\leqslant s\}}$ is even. By symmetry all $\mathbb{E}C_{i,j,2}(s)$
vanishes, so $\mathcal{M}_{2}=0$. For $k\geqslant3$, we use $V_{i}
=d_{i}^{\otimes k}1_{\{r_{i}\leqslant s\}}$ in $\left(  \ref{eq:bessel}
\right)  $. Their squared norms are $r_{i}^{k}1_{\{r_{i}\leqslant
s\}}\leqslant s^{k-1}r_{i}$, which verifies the square-integrability
hypothesis at each fixed $s>0$. Since $w_{j}(s)/(4s)\leqslant1/(2\lambda_{j}
)$,
\[
\sum_{j}\frac{w_{j}(s)}{4k!s^{k}}\Vert\mathbb{E}T_{j,k}\Vert^{2}\leqslant
\frac{1}{k!s^{k-1}}\sum_{i}\mathbb{E}\Vert V_{i}\Vert^{2}=\frac{1}{k!s^{k-1}
}\sum_{i}\mathbb{E}[r_{i}^{k}1_{\{r_{i}\leqslant s\}}].
\]
Fix $k\geqslant3.$ To compute $\mathcal{M}_{k}$, we keep $r_{i}>0$ fixed and
integrate the above inequality,
\[
r_{i}^{k}\int_{r_{i}}^{\infty}s^{1-k}\,ds=r_{i}^{k}\left[  \frac{s^{2-k}}
{2-k}\right]  _{r_{i}}^{\infty}=\frac{r_{i}^{k}r_{i}^{2-k}}{k-2}=\frac
{r_{i}^{2}}{k-2}.
\]
Using Tonelli's theorem followed by summation over $i$ therefore yields
\begin{equation}
\mathcal{M}_{k}\leqslant\frac{M}{k!(k-2)}\text{ for }k\geqslant3.
\label{eq:bk}
\end{equation}
On $\left\{  r_{i}=0\right\}  $ the tensor $T_{j,k}$ vanishes for every $s>0$,
thus the inequality is trivial.

\textbf{5}\emph{.} We now invoke Proposition \ref{prop:observation}.
Substituting $\left(  \ref{eq:split-variance}\right)  $ in $\left(
\ref{eq:hermite-bound}\right)  $ and adding $\left(  \ref{eq:zero}\right)  $,
$\left(  \ref{eq:c1}\right)  $, $\left(  \ref{eq:ck}\right)  $, $\left(
\ref{eq:b1}\right)  $, $\mathcal{M}_{2}=0$ and $\left(  \ref{eq:bk}\right)  $,
we obtain
\begin{align*}
\mathbb{E}[(S-G^{\ast})^{\top}L(S-G^{\ast})]  &  =\sum_{j}\int_{0}^{\infty
}w_{j}(s)\mathbb{E}R_{s,j}^{2}ds\leqslant\frac{8}{27}M+\sum_{k=1}^{\infty
}(\mathcal{V}_{k}+\mathcal{M}_{k})\\
&  \leqslant M\left(  \frac{8}{27}+\frac{1}{4}+1+\frac{2}{27}\sum
_{k=2}^{\infty}\frac{1}{k!(k-1)}+\sum_{k=3}^{\infty}\frac{1}{k!(k-2)}\right)
.
\end{align*}
The ratio of successive reciprocal factorials is $1/(k+1)$, which is at most
$1/3$ for $k\geqslant2$ and at most $1/4$ for $k\geqslant3$. In particular,
for each $m\geqslant0$,
\[
\frac{1}{(2+m)!}=\frac{1}{2!}\prod_{r=1}^{m}\frac{1}{2+r}\leqslant\tfrac{1}
{2}(1/3)^{m},\text{ }\frac{1}{(3+m)!}=\frac{1}{3!}\prod_{r=1}^{m}\frac{1}
{3+r}\leqslant\tfrac{1}{6}(1/4)^{m}.
\]
The products equal one for $m=0$. Since $k-1\geqslant1$ and $k-2\geqslant1$ on
the respective summation ranges, summing the two geometric series gives
\[
\sum_{k=2}^{\infty}\frac{1}{k!(k-1)}\leqslant\frac{1/2}{1-1/3}=\frac{3}
{4},\text{ }\sum_{k=3}^{\infty}\frac{1}{k!(k-2)}\leqslant\frac{1/6}
{1-1/4}=\frac{2}{9}.
\]
The resulting coefficient is
\[
\frac{8}{27}+\frac{1}{4}+1+\frac{2}{27}\frac{3}{4}+\frac{2}{9}=\frac
{32+27+108+6+24}{108}=\frac{197}{108}.
\]
Therefore
\[
\mathbb{E}\left[  (S-G^{\ast})^{\top}L(S-G^{\ast})\right]  \leqslant\frac
{197}{108}M=C_{0}\sum_{i}\mathbb{E}|U_{i}^{\prime}-U_{i}|^{4},
\]
which proves $\left(  \ref{eq:weighted}\right)  $ for bounded summands.
\end{proof}

We next try to remove the boundedness assumptions.

\begin{proof}
[Proof of Theorem~\ref{thm:weighted} (the general case)]\textbf{1}. For each
positive integer $m$ set
\[
U_{i}^{(m)}=U_{i}1_{\{\left\vert U_{i}\right\vert \leqslant m\}}
-\mathbb{E}\left[  U_{i}1_{\{|U_{i}|\leqslant m\}}\right]  .
\]
Write $H_{i,m}=U_{i}1_{\{|U_{i}|>m\}}$. Since $\mathbb{E}U_{i}=0$,
\[
U_{i}-U_{i}^{(m)}=H_{i,m}-\mathbb{E}H_{i,m}.
\]
By the triangle inequality in $L^{4}$ norm and Jensen's inequality we have
\[
\left\Vert U_{i}-U_{i}^{(m)}\right\Vert _{L^{4}}\leqslant\left\Vert
H_{i,m}\right\Vert _{L^{4}}+\left\vert \mathbb{E}H_{i,m}\right\vert
\leqslant2\left[  \mathbb{E}[|U_{i}|^{4}1_{\{|U_{i}|>m\}}\right]
^{1/4}\rightarrow0.
\]
Truncating each independent copy with the same deterministic map preserves its
law relative to the truncated original and preserves independence. Let
$S_{m}=\sum_{i}U_{i}^{(m)}$ and $L_{m}=\operatorname{Cov}(S_{m})$. Then
\[
\Vert S_{m}-S\Vert_{L^{4}}\leqslant\sum_{i}\Vert U_{i}^{(m)}-U_{i}\Vert
_{L^{4}}.
\]
Thus $S_{m}\rightarrow S$ in $L^{4}$. Since
\[
\Vert S_{m}-S\Vert_{L^{2}}\leqslant\Vert S_{m}-S\Vert_{L^{4}},
\]
there is also $L^{2}$ convergence and $\sup_{m}\Vert S_{m}\Vert_{L^{2}}
<\infty$.

\textbf{2}. All sums are centered, hence
\[
L_{m}-L=\mathbb{E}\left[  S_{m}S_{m}^{\top}-SS^{\top}\right]  .
\]
For column vectors $a,b$, $\Vert ab^{\top}\Vert^{2}=\sum_{j,k}a_{j}^{2}
b_{k}^{2}=|a|^{2}|b|^{2}$. Adding and subtracting $SS_{m}^{\top}$ inside the
expectation gives
\[
\Vert L_{m}-L\Vert\leqslant\mathbb{E}\Vert(S_{m}-S)S_{m}^{\top}\Vert
+\mathbb{E}\Vert S(S_{m}-S)^{\top}\Vert=\mathbb{E}[|S_{m}-S|\,|S_{m}
|]+\mathbb{E}[|S|\,|S_{m}-S|].
\]
Cauchy-Schwarz therefore yields
\[
\left\Vert L_{m}-L\right\Vert \leqslant\left(  \mathbb{E}\left\vert
S_{m}-S\right\vert ^{2}\right)  ^{1/2}\left[  (\mathbb{E}|S_{m}|^{2}
)^{1/2}+(\mathbb{E}|S|^{2})^{1/2}\right]  \rightarrow0
\]
as claimed.

\textbf{3}. Let $\lambda_{\min}(L)$ be the smallest eigenvalue of $L$, and use
$L^{\alpha}$ for the symmetric matrix obtained by raising each eigenvalue to
the power $\alpha$ in its spectral decomposition. For a unit vector
$u\in\mathbb{R}^{d}$,
\[
u^{\top}L_{m}u\geqslant\lambda_{\min}(L)-\Vert L_{m}-L\Vert,
\]
because
\[
\left\vert u^{\top}(L_{m}-L)u\right\vert \leqslant\Vert uu^{\top}\Vert\Vert
L_{m}-L\Vert=\Vert L_{m}-L\Vert.
\]
Once the last norm is less than $\lambda_{\min}(L)/2$, this lower bound is
positive for every unit $u$, so $L_{m}$ is positive definite. Write
$F_{m}=L_{m}^{1/2}-L^{1/2}$. For positive semi-definite real symmetric $d$ by
$d$ matrices $A,B$, the Frobenius bound of \cite[(3.2)]{Wihler2009Holder}
gives
\[
\left\Vert A^{1/2}-B^{1/2}\right\Vert ^{2}\leqslant\sqrt{d}\left\Vert
A-B\right\Vert .
\]
Here the cited Frobenius norm is exactly $\Vert\cdot\Vert$.

Both $L_{m}$ and $L$ are covariances in the same dimension $d$. Taking
$A=L_{m},B=L$ yields
\[
\Vert F_{m}\Vert\leqslant d^{1/4}\Vert L_{m}-L\Vert^{1/2}\rightarrow0,\text{
as }m\rightarrow\infty.
\]
\textbf{4}. For a centered square integrable vector $X$ with positive definite
covariance $A$, the invertible map $F:$ $\left(  x,y\right)  \mapsto\left(
A^{1/2}x,A^{1/2}y\right)  $ has
\[
|L^{1/2}x-L^{1/2}y|^{2}=(x-y)^{\top}A(x-y).
\]
Also $\operatorname{Cov}(A^{1/2}G)=A^{1/2}AA^{1/2}=A^{2}$ for $G\sim N(0,A)$.
Thus for any coupling $\pi$ of $\mathcal{L}\left(  X\right)  $ and $N\left(
0,A\right)  $, the first and second marginals of the transformed coupling
$F_{\#}\pi$ are $\mathcal{L}\left(  L^{1/2}X\right)  $ and $N\left(
0,A^{2}\right)  $. Moreover
\begin{equation}
\inf_{\pi\in\Pi(\mathcal{L}(X),N(0,A))}\int(x-y)^{\top}A(x-y)d\pi=W_{2}
^{2}(\mathcal{L}(A^{1/2}X),N(0,A^{2})). \label{eq:weighted-infimum}
\end{equation}
\textbf{5}. To apply $\left(  \ref{eq:weighted-infimum}\right)  $ with
$S_{m},L_{m}$, we first observe that
\[
\left\Vert L_{m}^{1/2}S_{m}-L^{1/2}S\right\Vert _{L^{2}}=\left\Vert
L_{m}^{1/2}\left(  S_{m}-S\right)  +F_{m}S\right\Vert _{L^{2}}\leqslant
\left\Vert L_{m}^{1/2}\right\Vert \left\Vert S_{m}-S\right\Vert _{L^{2}
}+\left\Vert F_{m}\right\Vert \left\Vert S\right\Vert _{L^{2}}\rightarrow0.
\]
Here we used $\left\vert Ax\right\vert \leqslant\left\Vert A\right\Vert
\left\vert x\right\vert $ for the Frobenius norm of any matrix $A$ and
$\left\Vert L_{m}^{1/2}\right\Vert \leqslant\left\Vert L^{1/2}\right\Vert
+\left\Vert F_{m}\right\Vert $. For $Z\sim N(0,I_{d})$, $L_{m}Z\sim
N(0,L_{m}^{2})$ and $LZ\sim N(0,L^{2})$. Moreover,
\[
\mathbb{E}\left\vert \left(  L_{m}-L\right)  Z\right\vert ^{2}
=\operatorname{tr}\left[  \left(  L_{m}-L\right)  ^{\top}\left(
L_{m}-L\right)  \mathbb{E}\left(  ZZ^{\top}\right)  \right]  =\left\Vert
L_{m}-L\right\Vert ^{2}.
\]
Let $P_{m}=\mathcal{L}(L_{m}^{1/2}S_{m})$, $P=\mathcal{L}(L^{1/2}S)$,
$Q_{m}=N(0,L_{m}^{2})$, and $Q=N(0,L^{2})$. These four distributions have
finite second moments. By the triangle inequality, we obtain, as
$m\rightarrow\infty,$ that
\[
\left\vert W_{2}(P_{m},Q_{m})-W_{2}(P,Q)\right\vert \leqslant W_{2}
(P_{m},P)+W_{2}(Q_{m},Q)\leqslant\Vert L_{m}^{1/2}S_{m}-L^{1/2}S\Vert_{L^{2}
}+\Vert L_{m}-L\Vert\rightarrow0.
\]
Let $d_{i}^{(m)}=\left(  U_{i}^{\prime}\right)  ^{(m)}-U_{i}^{(m)}$, we have
\[
\Vert d_{i}^{(m)}-d_{i}\Vert_{L^{4}}\leqslant\Vert\left(  U_{i}^{\prime
}\right)  ^{(m)}-U_{i}^{\prime}\Vert_{L^{4}}+\Vert U_{i}^{(m)}-U_{i}
\Vert_{L^{4}}=2\Vert U_{i}^{(m)}-U_{i}\Vert_{L^{4}}.
\]
It follows that
\[
\left\vert \Vert d_{i}^{(m)}\Vert_{L^{4}}-\Vert d_{i}\Vert_{L^{4}}\right\vert
\leqslant\Vert d_{i}^{(m)}-d_{i}\Vert_{L^{4}}\rightarrow0,\text{ as
}m\rightarrow\infty.
\]
Hence
\[
\sum_{i}\mathbb{E}|d_{i}^{(m)}|^{4}\rightarrow\sum_{i}\mathbb{E}|d_{i}
|^{4},\text{ as }m\rightarrow\infty.
\]
\textbf{6}. We have proved that For each sufficiently large $m$, the truncated
summands $U_{i}^{(m)}$ are bounded, independent, centered, and with
$S_{m}=\sum_{i}U_{i}^{(m)}$, $L_{m}=\operatorname{Cov}(S_{m})$ is positive
definite. So the already proved bounded case applies to $S_{m},$ that is,
there exists a minimizing coupling of $P_{m}$ and $Q_{m}$ such that $\left(
\ref{eq:weighted-infimum}\right)  $ holds (with $S_{m}$, $L_{m}$ replacing
$X$, $A$ there) and thus
\[
W_{2}^{2}(P_{m},Q_{m})\leqslant C_{0}\sum_{i}\mathbb{E}|d_{i}^{(m)}|^{4}.
\]
It follows that
\[
W_{2}^{2}(P,Q)=\lim_{m\rightarrow\infty}W_{2}^{2}(P_{m},Q_{m})\leqslant
C_{0}\lim_{m\rightarrow\infty}\sum_{i}\mathbb{E}|d_{i}^{(m)}|^{4}=C_{0}
\sum_{i}\mathbb{E}|d_{i}|^{4}.
\]
To get the form on the left side of $\left(  \ref{eq:weighted}\right)  ,$
suppose that the random vector $(\xi,\eta)$ has marginals $P,Q$ and its joint
coupling attains $W_{2}(P,Q)$. Such joint coupling exists. Let $\left(
S^{\prime},G^{\ast}\right)  =\left(  L^{-1/2}\xi,L^{-1/2}\eta\right)  $. Then
its marginal distributions are $\mathcal{L}(S)$ and $N(0,L)$ respectively,
and
\[
\mathbb{E}[(S^{\prime}-G^{\ast})^{\top}L(S^{\prime}-G^{\ast})]=\mathbb{E}
\left\vert \xi-\eta\right\vert ^{2}=W_{2}^{2}(P,Q)\leqslant C_{0}\sum
_{i}\mathbb{E}|d_{i}|^{4}.
\]
We now already have the desired form, since $\mathcal{L}(S)=\mathcal{L}
(S^{\prime}).$ To restore the orginally given $S.\,$Suppose $S$ is defined on
the probability space on $(\Omega,\mathcal{A},\mathbb{P})$. We disintegrate
the joint distribution of $\left(  S^{\prime},G^{\ast}\right)  $ as
$\mathcal{L}(S)(dx)\kappa_{\ast}(x,dg)$. On $\Omega\times\mathbb{R}^{d}$ we
use the measure $\mathbb{P}(d\omega)\kappa_{\ast}(S(\omega),dg)$. Then for
Borel $A,D$, the probability of $\{S\in A,g\in D\}$ is $\int_{A}\kappa_{\ast
}(x,D)\mathcal{L}(S)(dx)$. Thus $S(\omega)$ and the second marginal have the
required joint law on $\Omega\times\mathbb{R}^{d}$. So the proof of $\left(
\ref{eq:weighted}\right)  $ is completed.
\end{proof}

\section{One coupling for all coordinates}

\label{sec:separation}

We first isolate the passage from scalar inequalities to simultaneous
coordinate bounds. It uses only compactness of the set of couplings and
finite-dimensional separation \cite[Corollary~11.4.2]{rockafellar1997convex}.

\begin{lemma}
\label{lem:separate} Let $m\geqslant1$ be an integer, $P,Q$ be probabilities
on $\mathbb{R}^{m}$ with finite second moments, and $z\in\lbrack0,\infty)^{m}
$. For $\pi\in\Pi(P,Q),$ let $v(\pi)=(v_{1}(\pi),\ldots,v_{m}(\pi))$ with
\[
v_{j}(\pi)=\int(x_{j}-y_{j})^{2}d\pi(x,y).
\]
If for every $a\in(0,\infty)^{m}$,
\begin{equation}
\inf_{\pi\in\Pi(P,Q)}\sum_{j}a_{j}v_{j}(\pi)\leqslant\sum_{j}a_{j}z_{j},
\label{eq:dual}
\end{equation}
then some $\pi\in\Pi(P,Q)$ satisfies $v_{j}(\pi)\leqslant z_{j}$ for every $j$.
\end{lemma}

\begin{proof}
\textbf{1}. For $a\geqslant0$ coordinatewise, nonnegativity of the costs
gives
\[
\inf_{\pi}\sum_{j}a_{j}v_{j}(\pi)\leqslant\inf_{\pi}\sum_{j}(a_{j}+\eta
)v_{j}(\pi)\leqslant\sum_{j}(a_{j}+\eta)z_{j}\text{ for }\eta>0.
\]
Here $a+\eta(1,\ldots,1)$ is strictly positive, hence belongs to
$(0,\infty)^{m}$ even if some $a_{j}=0$. The right side is $\sum_{j}a_{j}
z_{j}+\eta\sum_{j}z_{j}$, with $\sum_{j}z_{j}<\infty$. Letting $\eta
\rightarrow0$ therefore extends $\left(  \ref{eq:dual}\right)  $ to all
nonnegative $a$.

\textbf{2}. Define
\[
C=\left\{  c\in\mathbb{R}^{m}:\exists\text{ }\pi\in\Pi(P,Q)\text{ such that
}v_{j}(\pi)\leqslant c_{j}\text{ for all }j\right\}  .
\]
Since
\[
v_{j}(\pi)\leqslant2\int x_{j}^{2}\,P(dx)+2\int y_{j}^{2}\,Q(dy),
\]
each coordinate cost is finite. The set $C$ is non-empty, since the product
coupling satisfies $v(P\otimes Q)\in C$. If $c_{0},c_{1}\in C$ with the
accompanying coupling being $\pi_{0},\pi_{1}\in\Pi(P,Q)$, then $(1-\theta
)\pi_{0}+\theta\pi_{1}$ has the same marginals and cost vector $(1-\theta
)v(\pi_{0})+\theta v(\pi_{1})$ for $0\leqslant\theta\leqslant1$. Thus $C$ is
convex. Also any increase in a coordinate bound $c_{j}$ preserves the
controlling inequality, so
\begin{equation}
C+[0,\infty)^{m}=C. \label{eq:C-cone}
\end{equation}
The set $C$ is closed. Indeed, let $c^{(k)}\in C$ converging to $c$ with
accompanying coupling $\pi_{k}\in\Pi(P,Q)$. The sequence $\pi_{k}$ tight, thus
up to a subsequence converges weakly to some $\pi\in\Pi(P,Q)$. For each $j$,
let $\psi_{j}(x,y)=(x_{j}-y_{j})^{2}$. For every $b>0$, by Fatou's lemma we
have
\[
\int(\psi_{j}\wedge b)d\pi=\lim_{k}\int(\psi_{j}\wedge b)d\pi_{k}
\leqslant\liminf_{k}v_{j}(\pi_{k})\leqslant\lim_{k\rightarrow\infty}
c_{j}^{(k)}=c_{j}.
\]
Increasing $b$ and using monotone convergence we obtain that $v_{j}
(\pi)\leqslant c_{j}$ for every $j$, so $c\in C$.

\textbf{3}. Suppose to the contrary that $z\notin C$, we can now apply
\cite[Corollary~11.4.2]{rockafellar1997convex} to the non-empty disjoint
convex sets $\{z\}$ and $C$. Both are closed and the singleton $\{z\}$ is
compact. So there exists a $a\neq0$ such that
\begin{equation}
\left\langle a,z\right\rangle <\inf_{c\in C}\left\langle a,c\right\rangle .
\label{eq:strict separation}
\end{equation}
Every $a_{j}$ is nonnegative. Indeed, for $c\in C$ and $r\geqslant0$,
$c+re_{j}\in C$ by $\left(  \ref{eq:C-cone}\right)  $, whereas $a_{j}<0$ would
force the infimum to be $-\infty$. To force a contradiction, we observe that
the set $C$ has the representation
\[
C=\left\{  v(\pi)+r:\pi\in\Pi(P,Q),\ r\in\lbrack0,\infty)^{m}\right\}  .
\]
For every pair $\pi,r$ from this representation, we have
\[
\left\langle a,v(\pi)+r\right\rangle =\left\langle a,v(\pi)\right\rangle
+\sum_{j}a_{j}r_{j}\geqslant\left\langle a,v(\pi)\right\rangle ,
\]
with equality at $r=0$. Therefore
\[
\inf_{c\in C}\langle a,c\rangle=\inf_{\pi\in\Pi(P,Q)}\left\langle
a,v(\pi)\right\rangle \leqslant\left\langle a,z\right\rangle .
\]
The equality holds due to the extension of $\left(  \ref{eq:dual}\right)  $ to
nonnegative weights proved at the beginning. This contradicts strict
separation $\left(  \ref{eq:strict separation}\right)  $. Thus we have proved
$z\in C.$
\end{proof}

We are now ready to prove the main theorem.

\begin{proof}
[Proof of Theorem~\ref{thm:main}]\textbf{1}. In the eigenbasis fixed in
Theorem~\ref{thm:main}, we expand the covariance of the normalized sum,
\[
\mathbb{E}\left(  W_{n,j}W_{n,\ell}\right)  =\frac{1}{n}\sum_{i,h=1}
^{n}\mathbb{E}\left(  X_{i,j}X_{h,\ell}\right)  .
\]
Note the second indices in the subscripts denote components under the fixed
eigenbasis. If $i\neq h$, the expectation is $(\mathbb{E}X_{i,j}
)(\mathbb{E}X_{h,\ell})=0$ by independence and centering. The remaining
diagonal terms have
\[
\frac{1}{n}\sum_{i=1}^{n}\mathbb{E}[X_{i,j}X_{i,\ell}]=\operatorname{Cov}
(W_{n})_{j,\ell}=\sigma_{j}^{2}1_{\{j=\ell\}}.
\]
If $\sigma_{j}=0$, taking $\ell=j$ yields $0=n\sigma_{j}^{2}=\sum
_{i}\mathbb{E}X_{i,j}^{2}$. So $X_{i,j}=0$ almost surely for every $i$. In
particular these coordinates can be omitted in $\left(  \ref{eq:budget}
\right)  $ and restored with value zero. Retaining all original coordinate
labels, we let
\[
\Gamma_{j,k}=\frac{1}{n^{2}}\sum_{i=1}^{n}\mathbb{E}[\Delta_{i,j}^{2}
\Delta_{i,k}^{2}],\text{ }1\leqslant j,k\leqslant d\text{. }
\]
Using $\left\vert \Delta_{i}\right\vert ^{2}=\sum_{k}\Delta_{i,k}^{2}$ in the
fixed orthonormal basis, we can then write
\[
b_{j}=\sigma_{j}^{-2}\sum_{k=1}^{d}\Gamma_{j,k}.
\]
Let
\[
J=\left\{  j:\sigma_{j}>0,C_{0}b_{j}<2\sigma_{j}^{2}\right\}  .
\]
For $\sigma_{j}>0$,
\[
\sigma_{j}^{2}b_{j}=\frac{1}{n^{2}}\sum_{i}\sum_{k}\mathbb{E}[\Delta_{i,j}
^{2}\Delta_{i,k}^{2}]=\sum_{k}\left(  \frac{1}{n^{2}}\sum_{i}\mathbb{E}
[\Delta_{i,j}^{2}\Delta_{i,k}^{2}]\right)  =\sum_{k}\Gamma_{j,k}.
\]
All entries are nonnegative, so the identity holds also for infinite $b_{j}$,
and $\Gamma_{j,k}=\Gamma_{k,j}$. For $j\in J$ and each $i$,
\[
\mathbb{E}\left\vert \Delta_{i,j}\right\vert ^{4}\leqslant\sum_{h}
\mathbb{E}\left\vert \Delta_{h,j}\right\vert ^{4}=n^{2}\Gamma_{j,j}\leqslant
n^{2}\sigma_{j}^{2}b_{j}<\infty.
\]
Since the copy $X_{i}^{\prime}$ is independent of the original $X_{i}$, so
\[
\mathbb{E}\left[  \Delta_{i,j}\mid X_{i,j}\right]  =\mathbb{E}\left[
X_{i,j}^{\prime}\mid X_{i,j}\right]  -X_{i,j}=\mathbb{E}X_{i,j}^{\prime
}-X_{i,j}=-X_{i,j}.
\]
By Jansen inequality, $\mathbb{E}|X_{i,j}|^{4}\leqslant\mathbb{E}|\Delta
_{i,j}|^{4}<\infty$. In particular, with $X_{i,J}=\left(  X_{i,j}\right)
_{j\in J},$ we have
\begin{equation}
\left\vert X_{i,J}\right\vert ^{4}=\left(  \sum_{j\in J}X_{i,j}^{2}\right)
^{2}\leqslant\left\vert J\right\vert \sum_{j\in J}X_{i,j}^{4}.
\label{eq:J-4moment}
\end{equation}
Hence, the vector $X_{i,J}$ of coordinates in $J$ has finite fourth moment for
each $i$.

\textbf{2}. If $J=\varnothing$, use the product coupling $\mathcal{L}\left(
W_{n}\right)  \otimes N\left(  0,\Sigma\right)  $. Under this coupling $G$ is
independent of $W_{n}$. For every $\sigma_{j}>0$, the definition of $J$ forces
$C_{0}b_{j}\geqslant2\sigma_{j}^{2}$, and thus
\begin{align}
\mathbb{E}(W_{n,j}-G_{j})^{2}  &  =\mathbb{E}W_{n,j}^{2}+\mathbb{E}G_{j}
^{2}-2\mathbb{E}[W_{n,j}G_{j}]=\sigma_{j}^{2}+\sigma_{j}^{2}-2(\mathbb{E}
W_{n,j})(\mathbb{E}G_{j})\label{eq:c-cost}\\
&  =2\sigma_{j}^{2}=\min\{2\sigma_{j}^{2},C_{0}b_{j}\}.\nonumber
\end{align}
It follows that
\[
W_{2}^{2}(\mathcal{L}(W_{n,j}),N(0,\sigma_{j}^{2}))\leqslant\mathbb{E}
(W_{n,j}-G_{j})^{2}=\min\{2\sigma_{j}^{2},C_{0}b_{j}\}.
\]
For those $j$ with $\sigma_{j}=0,$ $(W_{n,j},G_{j})$ vanish almost surely in
both marginals. Therefore $\left(  \ref{eq:main}\right)  $ holds when
$J=\varnothing$.

\textbf{3}. Suppose next $J\neq\varnothing$. Fix arbitrary weights $a_{j}>0$,
$j\in J$. Define $c_{j}>0$ by
\[
c_{j}^{2}=\sqrt{a_{j}}/\sigma_{j}.
\]
To apply Theorem~\ref{thm:weighted}, we let $U_{i}=(c_{j}X_{i,j}/\sqrt
{n})_{j\in J}.$ These are vectors in $\mathbb{R}^{J}$. Accordingly let
$S=\sum_{i}^{n}U_{i}$, $L$ be the corresponding $\left\vert J\right\vert $ by
$\left\vert J\right\vert $ covariance matrix. Each $U_{i}$ is a linear
function of $X_{i}$, so these summands of $S$ remain independent and centered
and have the finite fourth moments $\left(  \ref{eq:J-4moment}\right)  $ just
established . Theorem~\ref{thm:weighted} therefore yields a Gaussian vector
$G^{\ast}\in\mathbb{R}^{J},$ $G^{\ast}\sim N\left(  0,L\right)  $ and a
coupling of $S,$ $G^{\ast}$ such that $\left(  \ref{eq:weighted}\right)  $
holds. Denote the coupling by $\gamma_{J}$. For $j,k\in J$, the components of
$S$ satisfy
\[
S_{j}=\sum_{i}\frac{c_{j}X_{i,j}}{\sqrt{n}}=c_{j}W_{n,j},\text{ }
\mathbb{E}S_{j}=c_{j}\mathbb{E}W_{n,j}=0,
\]
and
\[
\operatorname{Cov}(S)_{j,k}=\mathbb{E}[S_{j}S_{k}]=c_{j}c_{k}\mathbb{E}
[W_{n,j}W_{n,k}]=c_{j}c_{k}\sigma_{j}^{2}1_{\{j=k\}}=c_{j}^{2}\sigma_{j}
^{2}1_{\{j=k\}}.
\]
Thus $L=\operatorname{diag}(\lambda_{j})_{j\in J}$ with $\lambda_{j}=c_{j}
^{2}\sigma_{j}^{2}>0$. Let $G_{j}^{\prime}=G_{j}^{\ast}/c_{j}$. Then
$(G_{j}^{\prime})_{j\in J}$ is Gaussian with covariance
\[
\mathbb{E}\left[  G_{j}^{\prime}G_{k}^{\prime}\right]  =\frac{\lambda
_{j}1_{\{j=k\}}}{c_{j}c_{k}}=\sigma_{j}^{2}1_{\{j=k\}}.
\]
Consider the scaling $x\mapsto(c_{j}x_{j})_{j\in J}$ which is invertible
because all $c_{j}>0$. Applying its inverse to both marginals of $\pi_{J}$
yields a coupling of $(W_{n,j})_{j\in J}$ and $(G_{j}^{\prime})_{j\in J}$.
Denote the resulting coupling by $\gamma_{J}^{a}$. Noting that $\lambda
_{j}c_{j}^{2}=\sigma_{j}^{2}c_{j}^{4}=\sigma_{j}^{2}(a_{j}/\sigma_{j}
^{2})=a_{j}$, we therefore have
\begin{align}
(S-G^{\ast})^{\top}L(S-G^{\ast})  &  =\sum_{j,k\in J}(S_{j}-G_{j}^{\ast
})L_{j,k}(S_{k}-G_{k}^{\ast})\label{eq:scaled-cost}\\
&  =\sum_{j\in J}\lambda_{j}\left(  c_{j}W_{n,j}-c_{j}G_{j}^{\prime}\right)
^{2}=\sum_{j\in J}a_{j}(W_{n,j}-G_{j}^{\prime})^{2}.\nonumber
\end{align}
Let $U_{i}^{\prime}=(c_{j}X_{i,j}^{\prime}/\sqrt{n})_{j\in J}$, then
$(U_{i}^{\prime}-U_{i})_{j}=c_{j}\Delta_{i,j}/\sqrt{n}$. In particular
$|U_{i}^{\prime}-U_{i}|^{2}=n^{-1}\sum_{j\in J}c_{j}^{2}\Delta_{i,j}^{2}$.
Hence
\begin{align}
\sum_{i}\mathbb{E}\left\vert U_{i}^{\prime}-U_{i}\right\vert ^{4}  &
=\frac{1}{n^{2}}\sum_{i}\mathbb{E}\left(  \sum_{j\in J}c_{j}^{2}\Delta
_{i,j}^{2}\right)  ^{2}=\frac{1}{n^{2}}\sum_{i}\sum_{j,k\in J}c_{j}^{2}
c_{k}^{2}\mathbb{E}\left(  \Delta_{i,j}^{2}\Delta_{i,k}^{2}\right)
\label{eq:scaled}\\
&  =\sum_{j,k\in J}c_{j}^{2}c_{k}^{2}\Gamma_{j,k}=\sum_{j,k\in J}\Gamma
_{j,k}\frac{\sqrt{a_{j}a_{k}}}{\sigma_{j}\sigma_{k}}.\nonumber
\end{align}
Let
\[
P=\mathcal{L}((W_{n,j})_{j\in J}),\text{ }Q=N(0,\operatorname{diag}(\sigma
_{j}^{2})_{j\in J}).
\]
The equations~$\left(  \ref{eq:weighted}\right)  $, $\left(
\ref{eq:scaled-cost}\right)  $ and $\left(  \ref{eq:scaled}\right)  $ provide
a bound on the infimum of the weighted cost over $\Pi(P,Q)$,
\[
\inf_{\pi\in\Pi(P,Q)}\sum_{j\in J}a_{j}v_{j}(\pi)\leqslant\sum_{j\in J}
a_{j}v_{j}(\gamma_{J}^{a})=\mathbb{E}[(S-G^{\ast})^{\top}L(S-G^{\ast
})]\leqslant C_{0}\sum_{i}\mathbb{E}|U_{i}^{\prime}-U_{i}|^{4}.
\]
On the other hand,
\begin{align}
\inf_{\pi\in\Pi(P,Q)}\sum_{j\in J}a_{j}v_{j}(\pi)  &  \leqslant C_{0}
\sum_{j,k\in J}\Gamma_{j,k}\frac{\sqrt{a_{j}a_{k}}}{\sigma_{j}\sigma_{k}
}\leqslant\frac{C_{0}}{2}\sum_{j,k\in J}\Gamma_{j,k}\left(  \frac{a_{j}
}{\sigma_{j}^{2}}+\frac{a_{k}}{\sigma_{k}^{2}}\right) \label{eq:allweights}\\
&  =C_{0}\sum_{j\in J}\frac{a_{j}}{\sigma_{j}^{2}}\sum_{k\in J}\Gamma
_{j,k}\leqslant\sum_{j\in J}a_{j}C_{0}b_{j}.\nonumber
\end{align}
The equality uses that $\Gamma_{j,k}=\Gamma_{k,j}.$

\textbf{4}. Applying Lemma~\ref{lem:separate} with $z_{j}=C_{0}b_{j}$ now
yields a single coupling, denoted by $\pi_{J}\in\Pi(P,Q),$ for the coordinates
in $J.$ Use $x$,$g$ to denote the first and second marginal variables. To
extend $\pi_{J}$ to the full coordinate space, we disintegrate it as
\[
\pi_{J}=P\left(  dx_{J}\right)  \kappa\left(  x_{J},dg\right)  ,
\]
where $x_{J}=\left(  x_{j}\right)  _{j\in J}.$ On the full coordinate space,
we define the probability coupling
\begin{equation}
\pi=\mathcal{L}(W_{n})(dx)\kappa(x_{J},dg_{J})\prod_{j\in J^{c}}N(0,\sigma
_{j}^{2})(dg_{j}), \label{eq:lift}
\end{equation}
with the contention that $N(0,0)$ is the unit mass at zero and an empty
product equals $1$. Under the construction $g_{J^{c}}$ is independent of
$\left(  x,g_{J}\right)  $ and the $g$ marginal has distribution $N\left(
0,\Sigma\right)  $ in the fixed eigenbasis. Thus for each $j\notin J$ with
$\sigma_{j}>0,$ as in $\left(  \ref{eq:c-cost}\right)  $ we have, by centering
and the independence of $x_{J^{c}}$ and $g_{J^{c}}$, that
\[
\mathbb{E}(W_{n,j}-G_{j})^{2}=\mathbb{E}W_{n,j}^{2}+\mathbb{E}G_{j}
^{2}-2(\mathbb{E}W_{n,j})(\mathbb{E}G_{j})=2\sigma_{j}^{2}\leqslant C_{0}
b_{j}.
\]
For $j\in J$, the $(x_{J},g_{J})$ marginal preserves the coupling $\pi_{J}$,
and $C_{0}b_{j}<2\sigma_{j}^{2}$. Hence the single coupling $\pi$ satisfies
\[
\mathbb{E}_{\pi}(W_{n,j}-G_{j})^{2}\left\{
\begin{array}
[c]{ll}
\leqslant C_{0}b_{j}=\min\{2\sigma_{j}^{2},C_{0}b_{j}\}, & j\in J,\\
=2\sigma_{j}^{2}=\min\{2\sigma_{j}^{2},C_{0}b_{j}\}, & j\notin J,\ \sigma
_{j}>0,\\
=0, & \sigma_{j}=0.
\end{array}
\right.
\]
For each $j$ with $\sigma_{j}>0$, the coordinate marginal of $\pi$ is an
admissible scalar coupling in $\Pi(\mathcal{L}(W_{n,j}),N(0,\sigma_{j}^{2}))$,
so
\[
W_{2}^{2}\left(  \mathcal{L}(W_{n,j}),N(0,\sigma_{j}^{2})\right)
\leqslant\int\left(  x_{j}-g_{j}\right)  ^{2}\pi(dx,dg)\leqslant\min\left\{
2\sigma_{j}^{2},C_{0}b_{j}\right\}  .
\]
Therefore the desired bound $\left(  \ref{eq:main}\right)  $ is proved,
including the extreme case $b_{j}=\infty$.
\end{proof}

\section{Coordinate-wise central limit theorem in Hilbert space}

\label{sec:hilbert-prediction}

This section extends Theorem~\ref{thm:main} to a real separable Hilbert space.
The resulting coupling controls every coordinate in a prescribed covariance
eigenbasis and every projection onto a subset of that basis.

For a review of notation and basic setup \cite[Chapter VI]{reed1980methods}
\cite[Chapter 6-8]{kreyszig1991introductory}, let $H$ be a real separable
Hilbert space with inner product $\langle\cdot,\cdot\rangle$ and norm
$|\cdot|$. For $x,y\in H$, $x\otimes y$ denotes the operator $u\mapsto\langle
y,u\rangle x$ by
\[
\left(  x\otimes y\right)  \left(  u\right)  :u\mapsto\langle y,u\rangle x.
\]
An $H$-valued random element $V$ is a measurable map into $H$ with its Borel
sigma-field. It is \textit{square integrable }when $\mathbb{E}\left\vert
V\right\vert ^{2}<\infty$. In a separable Hilbert space this condition implies
Bochner integrability \cite[Appendix E]{cohn2013measure}, since $\mathbb{E}
\left\vert V\right\vert \leqslant\left(  \mathbb{E}\left\vert V\right\vert
^{2}\right)  ^{1/2}$. Its mean $\mathbb{E}V\in H$ is characterized by
$\langle\mathbb{E}V,u\rangle=\mathbb{E}\langle V,u\rangle$. It is
\textit{centered} if $\mathbb{E}V=0$. The\textit{ covariance operator }of a
square integrable $V$ is
\[
A=\operatorname{Cov}(V)=\mathbb{E}[(V-\mathbb{E}V)\otimes(V-\mathbb{E}V)].
\]
Thus
\[
Au=\mathbb{E}[\langle V-\mathbb{E}V,u\rangle(V-\mathbb{E}V)].
\]
Equivalently,
\[
\left\langle Au,z\right\rangle =\mathbb{E}[\left\langle V-\mathbb{E}
V,u\right\rangle \left\langle V-\mathbb{E}V,z\right\rangle ]\text{ for
}\forall u,z\in H.
\]
For centered $V$, this becomes $A=\mathbb{E}\left(  V\otimes V\right)  $.

A bounded self-adjoint operator $A$ is \textit{positive} if $\langle
Au,u\rangle\geqslant0$ for every $u$. For positive operators, $A\preccurlyeq
\Sigma$ means $\langle Au,u\rangle\leqslant\langle\Sigma u,u\rangle$ for every
$u$. A positive operator is\textit{ trace class }if $\operatorname{tr}
A=\sum_{j}\langle Af_{j},f_{j}\rangle<\infty$ for an orthonormal basis
$(f_{j})$. The sum is independent of the basis \cite[Theorem VI.18]
{reed1980methods}. For a covariance operator, Tonelli and Parseval give
$\operatorname{tr}A=\mathbb{E}\left\vert V-\mathbb{E}V\right\vert ^{2}$. For a
positive trace class $\Sigma$, its\textit{ positive spectral subspace }is the
closed subspace
\[
H_{0}=(\ker\Sigma)^{\perp}=\overline{\operatorname{Ran}\Sigma}=\overline
{\operatorname{span}\{e\in H:\Sigma e=\lambda e\text{ for some }\lambda>0\}}.
\]
Here $\ker\Sigma=\{u:\Sigma u=0\}$, $\operatorname{Ran}\Sigma=\{\Sigma u:u\in
H\}$, the bar denotes norm closure, and $M^{\perp}$ is the set of vectors
orthogonal to every vector in $M$. Self-adjointness gives $(\operatorname{Ran}
\Sigma)^{\perp}=\ker\Sigma$, that is, orthogonality to every $\Sigma u$ is
equivalent to $\Sigma z=0$. The last equality follows from the compact
self-adjoint spectral theorem \cite{kreyszig1991introductory}. Thus $H_{0}$ is
the closed span of the eigenvectors with positive variance. Its orthogonal
complement consists of zero-variance directions. The closure is needed because
$\operatorname{Ran}\Sigma$ need not be closed in infinite dimension. If
$\Sigma=0$, then $H_{0}=\{0\}$. Write $N_{H}(0,\Sigma)$ for the centered
Gaussian distribution characterized by $\langle G,u\rangle\sim N(0,\langle
\Sigma u,u\rangle)$ for every $u\in H$.

For applications in statistical learning
\cite{hsing2015theoretical,hastie2009elements}, a \textit{feature} is an
$H$-valued random input $V$ at which a linear function is evaluated. Its
feature distribution is the Borel probability measure $\xi_{V}=\mathcal{L}
(V)$, defined by $\xi_{V}(B)=\mathbb{P}(\xi\in B)$ for Borel $B\subset H$.
Thus square integrability is $\int_{H}|v|_{V}^{2}\xi_{V}(dv)<\infty$, and the
covariance of the feature is the operator $A$ just defined.

For a coefficient $w\in H$, the \textit{linear prediction} at input $v\in H$
is $\langle v,w\rangle$. In practice, e.g.
\cite{mandt2017stochastic,maddox2019simple}, $W_{n}$ is an average of models
and $G$ is its Gaussian approximation. So the two predictions at the same
future input are $\langle V,W_{n}\rangle$ and $\langle V,G\rangle$. Their
difference is $\langle V,W_{n}-G\rangle$. We choose the coupling before
drawing $V$ and require independence of $V$ from the \textit{whole pair}
$(W_{n},G)$. Thus the joint distribution of $V,W_{n},G$ is $\pi(dw,dg)\xi
_{V}(dv)$. For centered square-integrable $V$ with covariance $A$, the mean
squared prediction error is
\[
\int\int(\langle v,w\rangle-\langle v,g\rangle)^{2}\xi_{V}(dv)\pi
(dw,dg)=\int\langle A(w-g),w-g\rangle\pi(dw,dg).
\]

For Borel probability laws $\mu,\nu$ on $H$ with finite second moments, write
$\Pi(\mu,\nu)$ for probability measures on $H\times H$ whose first and second
marginals are $\mu,\nu$. The ordinary quadratic Hilbert transport distance is
\[
W_{2}(\mu,\nu)=\left(  \inf_{\pi\in\Pi(\mu,\nu)}\int_{H\times H}|x-y|^{2}
\pi(dx,dy)\right)  ^{1/2}.
\]
In an orthonormal basis $(f_{j})$, Parseval's identity (e.g. \cite[Theorem
3.6-3]{kreyszig1991introductory}) gives $|x-y|^{2}=\sum_{j}\langle
x-y,f_{j}\rangle^{2}$

For a positive trace class covariance operator $\Sigma$, the corresponding
squared prediction transport cost is
\[
\mathcal{T}_{\Sigma}^{2}(\mu,\nu)=\inf_{\pi\in\Pi(\mu,\nu)}\int|\Sigma
^{1/2}(x-y)|^{2}\pi(dx,dy).
\]
If the covariance of $V$ is $\Sigma$, then the \textit{conditional prediction
error} under each fixed coupling $\left(  W_{n},G\right)  $ is
\begin{equation}
\mathbb{E}[\langle V,W_{n}-G\rangle^{2}\mid W_{n},G]=\langle\Sigma
(W_{n}-G),W_{n}-G\rangle=|\Sigma^{1/2}(W_{n}-G)|^{2}.
\label{eq:conditional-prediction-error}
\end{equation}
The first equality uses our assumption that $V$ is independent of the pair
$(W_{n},G)$. In the language of machine learning, conditional prediction error
is the prediction error given an already determined model. Taking expectation
on both sides of $\left(  \ref{eq:conditional-prediction-error}\right)  $
w.r.t. the coupling of $\left(  W_{n},G\right)  $ therefore yields another
form of $\left(  \ref{eq:transport-identity}\right)  $ that we will frequently
used
\begin{equation}
\mathcal{T}_{\Sigma}^{2}(\mu,\nu)=\inf_{\pi\in\Pi(\mu,\nu)}\mathbb{E}[\langle
V,W_{n}-G\rangle^{2}]. \label{eq:transport-identity}
\end{equation}
In this prediction error , we assume a linear model for simplicity, nonlinear
models can inheric the error bound through linearization under suitable
regularity assumptions.

Let $\left\Vert \Sigma\right\Vert _{\mathrm{op}}=\sup_{\left\vert u\right\vert
\leqslant1}\left\vert \Sigma u\right\vert $. By integrating $\left\vert
\Sigma^{1/2}z\right\vert ^{2}\leqslant\left\Vert \Sigma\right\Vert
_{\mathrm{op}}\left\vert z\right\vert ^{2}$, we have the relation
\[
\mathcal{T}_{\Sigma}^{2}\left(  \mu,\nu\right)  \leqslant\left\Vert
\Sigma\right\Vert _{\mathrm{op}}W_{2}^{2}\left(  \mu,\nu\right)  .
\]

\begin{theorem}
[Simultaneous Hilbert coupling and prediction]\label{thm:hilbert-prediction}
Let $H$ be a real separable Hilbert space, let $n\geqslant1$ be an integer,
and let $X_{1},\ldots,X_{n}$ be independent centered $H$-valued random
elements with $\mathbb{E}\left\vert X_{i}\right\vert ^{4}<\infty$. Keep
$C_{0}$ as in $\left(  \ref{eq:C0}\right)  $ and define
\[
W_{n}=n^{-1/2}\sum_{i}X_{i},\text{ }\Sigma=\mathbb{E}\left[  W_{n}\otimes
W_{n}\right]  ,\text{ }\Delta_{i}=X_{i}^{\prime}-X_{i},
\]
where $\left(  X_{i}^{\prime}\right)  _{i=1}^{n}$ is an independent copy of
$\left(  X_{i}\right)  _{i=1}^{n}$. These moment assumptions, as will be
verified, guarantees the existence of an orthonormal eigenbasis $\left(
e_{j}\right)  _{j\in I}$ of $H_{0}=\left(  \ker\Sigma\right)  ^{\perp}$, the
positive spectral subspace with $\Sigma e_{j}=\lambda_{j}e_{j}$, $\lambda
_{j}>0$. Fix any one such eigenbasis. Let $I$ be the index set that enumerates
the eigenbasis, it could be empty $\varnothing,$ a fintie set $\{1,\ldots,r\}$
or $\mathbb{N}$. Let
\[
b_{j}=\frac{1}{n^{2}\lambda_{j}}\sum_{i}\mathbb{E}\left[  \left\vert
\Delta_{i}\right\vert ^{2}\left\langle \Delta_{i},e_{j}\right\rangle
^{2}\right]  ,\text{ }\beta_{j}=\min\left\{  2\lambda_{j},C_{0}b_{j}\right\}
.
\]
(i) There exists a coupling $\pi$ of $W_{n}$ and $G\sim N_{H}(0,\Sigma)$ such
that for every $j\in I$,
\begin{equation}
\mathbb{E}_{\pi}\left[  \langle W_{n}-G,e_{j}\rangle^{2}\right]
\leqslant\beta_{j}. \label{eq:hilbert-coordinate}
\end{equation}
(ii) For every deterministic subset $J\subset I$, let $P_{J}$ be the
orthogonal projection onto the closed span of $\{e_{j}:j\in J\}$, with
$P_{\varnothing}=0$. Then
\begin{equation}
\mathbb{E}_{\pi}\left\vert \Sigma^{1/2}P_{J}\left(  W_{n}-G\right)
\right\vert ^{2}\leqslant\frac{C_{0}}{n^{2}}\sum_{i}\mathbb{E}\left[
\left\vert \Delta_{i}\right\vert ^{2}\left\vert P_{J}\Delta_{i}\right\vert
^{2}\right]  . \label{eq:hilbert-local}
\end{equation}
(iii) For any centered $H$-valued feature distribution $\xi_{V}$ with
$\int|v|^{2}\xi_{V}(dv)<\infty$ and covariance $A=\int\left(  v\otimes
v\right)  \xi_{V}(dv)\preccurlyeq\Sigma$, draw $V$ independently of $\left(
W_{n},G\right)  $. On the product extension $\pi\otimes\xi_{V}$,
\begin{equation}
\mathbb{E}\left[  \left\langle V,W_{n}-G\right\rangle ^{2}\right]
\leqslant\frac{C_{0}}{n^{2}}\sum_{i}\mathbb{E}\left\vert \Delta_{i}\right\vert
^{4}. \label{eq:hilbert-risk}
\end{equation}
The coupling $\pi$ is chosen independently of $J$ and of $\xi_{V}$. Note the
coefficient is independent of $\xi_{V}$. In particular, in the i.i.d.\ case
the right hand side is $C_{0}\mathbb{E}\left\vert X_{1}^{\prime}
-X_{1}\right\vert ^{4}/n$, that is, the error bound (after taking square root)
has an order $O\left(  n^{-1/2}\right)  .$
\end{theorem}

\begin{proof}
\textbf{1}. We first verify that the spectral subspace and all $b_{j}$ are
well defined, and remove directions in which the sum has zero variance. Using
independence and centering we have, for $u,v\in H$, that
\begin{equation}
\left\langle \Sigma u,v\right\rangle =\frac{1}{n}\sum_{i=1}^{n}\mathbb{E}
\left[  \left\langle X_{i},u\right\rangle \left\langle X_{i},v\right\rangle
\right]  . \label{eq:covariance-identity}
\end{equation}
Moreover
\[
\mathbb{E}|W_{n}|^{2}=\frac{1}{n}\sum_{i}\mathbb{E}|X_{i}|^{2}<\infty.
\]
Thus $\Sigma$ is bounded, self-adjoint, and positive. For any orthonormal
basis $(f_{j})$ of $H$, using Tonelli theorem and Parseval identiy we have
\[
\sum_{j}\left\langle \Sigma f_{j},f_{j}\right\rangle =\mathbb{E}\sum
_{j}\left\langle W_{n},f_{j}\right\rangle ^{2}=\mathbb{E}\left\vert
W_{n}\right\vert ^{2}.
\]
In particular $\Sigma$ is trace class. It is also compact. Indeed, if $F_{k}$
projects onto the first $k$ vectors $f_{j}$, with $F_{k}=I_{H}$ once the whole
basis is included if $H$ is finite dimensional, where $I_{H}$ is the identity
on $H$, then $F_{k}\Sigma F_{k}=\mathbb{E}\left[  F_{k}W_{n}\otimes F_{k}
W_{n}\right]  $ is finite rank and $|W_{n}-F_{k}W_{n}|^{2}\leqslant|W_{n}
|^{2}$. The dominated convergence theorem therefore yields
\[
\left\Vert \Sigma-F_{k}\Sigma F_{k}\right\Vert _{\mathrm{op}}\leqslant
\mathbb{E}\left[  \left(  \left\vert W_{n}\right\vert +\left\vert F_{k}
W_{n}\right\vert \right)  \left\vert W_{n}-F_{k}W_{n}\right\vert \right]
\leqslant2\left(  \mathbb{E}\left\vert W_{n}\right\vert ^{2}\right)
^{1/2}\left(  \mathbb{E}\left\vert W_{n}-F_{k}W_{n}\right\vert ^{2}\right)
^{1/2}\rightarrow0.
\]
Here we used that $x\otimes x-y\otimes y=(x-y)\otimes x+y\otimes(x-y)$ and
$\Vert x\otimes y\Vert_{\mathrm{op}}=\left\vert x\right\vert \left\vert
y\right\vert $. The compact self-adjoint spectral theorem
\cite{kreyszig1991introductory} now gives the eigenbasis of $H_{0}=(\ker
\Sigma)^{\perp}$ with positive eigenvalues $\lambda_{j}>0$ satisfying $\Sigma
e_{j}=\lambda_{j}e_{j}$ and $\sum_{j}\lambda_{j}=\operatorname{tr}\Sigma$.
Since $\Sigma$ is bounded, $\ker\Sigma$ is closed subspace of $H.$ Thus
\[
H=H_{0}\oplus\left(  \ker\Sigma\right)  .
\]
The closed subspace $\ker\Sigma$ is separable. We can choose a countable dense
sequence in it and apply Gram-Schmidt procedure to obtain a countable
orthonormal basis of $\ker\Sigma.$ Extending $\left\{  e_{j}\right\}  $ with
this basis of $\ker\Sigma$ yields a countable orthonormal basis of $H.$ We
will henceforth use this basis. By $\left(  \ref{eq:covariance-identity}
\right)  $, $\sum_{i}\mathbb{E}\langle X_{i},u\rangle^{2}=0$ for $u\in
\ker\Sigma$. So $\langle X_{i},u\rangle=0$ almost surely. Thus intersecting
the probability one events over the countable orthonormal basis of $\ker
\Sigma$, and over $i\in\left\{  1,...,n\right\}  $, we see that $\left(
X_{i}\right)  _{i},$ $\left(  \Delta_{i}\right)  _{i},$ $W_{n}\in H_{0}$
almost surely. If $\Sigma=0$, then $I=\varnothing$, $H_{0}=\{0\}$, we could
take the deterministic coupling $(0,0)$ and all these vectors and costs
vanish. So all desired inequalities are satisfied trivially. Additionally, for
each $j$ with $\lambda_{j}>0$, by the moment assumption we have
\[
b_{j}\leqslant\frac{1}{n^{2}\lambda_{j}}\sum_{i}\mathbb{E}|\Delta_{i}
|^{4}\leqslant\frac{16}{n^{2}\lambda_{j}}\sum_{i}\mathbb{E}|X_{i}|^{4}
<\infty.
\]
Hence by the definition of $\beta_{j}$ and that $\Sigma$ is trace class, we
get $\sum_{j}\beta_{j}\leqslant2\sum_{j}\lambda_{j}<\infty$.

\textbf{2}. We next construct the fixed Gaussian marginal and obtains
couplings of its finite projections from Theorem~\ref{thm:main}. Let $Z_{j}$
be independent standard normals and
\[
G^{(k)}=\sum_{j\in I,\,j\leqslant k}\sqrt{\lambda_{j}}Z_{j}e_{j}.
\]
It is finite almost surely since, by Tonelli theorem, $\mathbb{E}\sum
_{j}\lambda_{j}Z_{j}^{2}=\sum_{j}\lambda_{j}<\infty$. By orthogonality,
$(G^{(k)})_{k}$ is almost surely Cauchy in $H_{0}$, thus has a limit in
$H_{0}$, since $H_{0}$ is closed (see e.g. \cite[p150 Problem 8]
{kreyszig1991introductory}). Denote the limit by $G.$ Then it satisfies
\[
\mathbb{E}\left\vert G-G^{(k)}\right\vert ^{2}=\sum_{j>k}\lambda
_{j}\rightarrow0\text{ as }k\rightarrow\infty.
\]
For $u\in H$ and $t\in\mathbb{R}$, by domianted convergence theorem we have
\[
\mathbb{E}e^{\mathrm{i}t\langle G,u\rangle}=\lim_{k}\exp\left(  -\frac{t^{2}
}{2}\sum_{j\leqslant k}\lambda_{j}\left\langle e_{j},u\right\rangle
^{2}\right)  =\exp(-t^{2}\left\langle \Sigma u,u\right\rangle /2).
\]
Thus the limit $G$ gives a Gaussian $\gamma=N_{H}(0,\Sigma)$ on $H_{0}.$ For
finite rank the above construction is a finite sum, and Theorem~\ref{thm:main}
applies on $H_{0}$ directly. To address the general infinite rank, we write
$\mu=\mathcal{L}(W_{n})$. Let $P_{k}$ project onto the finite spanned by the
first $k$ eigenvectors of $H_{0},$ $\operatorname{span}\{e_{1},\ldots,e_{k}
\}$. The vectors $P_{k}X_{i}$ retain independence and centering, have fourth
moments bounded by those of $X_{i}$. Moreover their normalized sum
$W_{n}^{\left(  k\right)  }=n^{-1/2}\sum_{i}P_{k}X_{i}$ has covariance
$\operatorname{diag}(\lambda_{1},\ldots,\lambda_{k})$. For $j\leqslant k$, the
finite dimensional counterpart of $b_{j}$ is given by
\[
b_{j}^{(k)}=\frac{1}{n^{2}\lambda_{j}}\sum_{i}\mathbb{E}[|P_{k}\Delta_{i}
|^{2}\langle\Delta_{i},e_{j}\rangle^{2}]\leqslant b_{j}.
\]
Applying Theorem~\ref{thm:main}, we therefore obtain a coupling $\rho_{k}$ of
$P_{k}W_{n}$ and $P_{k}G^{(\infty)}$ such that
\begin{equation}
\int\langle x-y,e_{j}\rangle^{2}\rho_{k}(dx,dy)\leqslant\min\{2\lambda
_{j},C_{0}b_{j}^{(k)}\}\leqslant\beta_{j},\text{ }\forall j\leqslant k.
\label{eq:firstk-cost}
\end{equation}

\textbf{3}. The projected couplings $\rho_{k}$ need not be consistent as $k$
increases. We extend each $\rho_{k}$\ to $H\times H$ with the same full
marginals $\mu,$ $\gamma$, so that weak compactness can be applied. Let
$\mu^{(k)}=\mathcal{L}(P_{k}W_{n})$ and $\gamma^{(k)}=\mathcal{L}(P_{k}G)$.
Write $x,u$ for the variables of $P_{k}W_{n},W_{n},$ and $\left(  y,g\right)
$ for those of $P_{k}G,G.$ By disintegrating respectively $\mathcal{L}
(P_{k}W_{n},W_{n}),$ $\mathcal{L}(P_{k}G,G)$ \cite[Theorem~8.5]
{kallenberg2021foundations} we obtain conditional probability kernels $\mu
_{k}(x,du)$ and $\gamma_{k}(y,dg)$ for the full vectors given their $P_{k}$
projections. Define the joint distribution $\pi_{k}$ as follows, for Borel
$A,B\subset H$,
\[
\pi_{k}(A\times B)=\int\mu_{k}(x,A)\gamma_{k}(y,B)\rho_{k}(dx,dy).
\]
The kernels are unique in the sense that there are Borel sets $N_{\mu},$
$N_{\gamma}$ such that $\mu^{(k)}\left(  N_{\mu}\right)  =0,$ $\gamma
^{(k)}\left(  N_{\gamma}\right)  =0$ and outside these null sets all versions
of kernels coincide. In particular
\[
\rho_{k}\left(  \left(  N_{\mu}\times H\right)  \cup\left(  H\times N_{\gamma
}\right)  \right)  \leqslant\mu^{(k)}\left(  N_{\mu}\right)  +\gamma
^{(k)}\left(  N_{\gamma}\right)  =0.
\]
Therefore the defintion of $\pi_{k}$ is independent of the versions of kernels
that are used in its construction. Since the two marginals of $\rho_{k}$ are
$\mu^{(k)},\gamma^{(k)}$, the two marginals of $\pi_{k}$ are
\begin{align*}
\pi_{k}(A\times H)  &  =\int\mu_{k}(x,A)\mu^{(k)}(dx)=\mu(A),\\
\pi_{k}(H\times B)  &  =\int\gamma_{k}(y,B)\gamma^{(k)}(dy)=\gamma(B).
\end{align*}
Thus $\pi_{k}\in\Pi(\mu,\gamma).$ The joint distribution of $(P_{k}u,P_{k}g)$
under $\pi_{k}$ is precisely $\rho_{k}.$ Indeed, the joint law of $(P_{k}
W_{n},W_{n})$ is supported on the closed set $\{(x,u):P_{k}u=x\}$. Its
disintegration has that $\mu_{k}(x,\{u:P_{k}u=x\})=1$ for $\mu^{(k)}$-almost
every $x$. The same applies to the Gaussian marginal. Hence $\pi_{k}$ is
compatible with $\rho_{k}$. Choices of the kernels on null sets have zero
integral because the conditioning marginals are fixed. Thus the first $k$
coordinate costs are at most $\beta_{j}$ by $\left(  \ref{eq:firstk-cost}
\right)  $,
\[
\int\left\langle x-y,e_{j}\right\rangle ^{2}\pi_{k}(dx,dy)=\int\left\langle
x-y,e_{j}\right\rangle ^{2}\rho_{k}(dx,dy)\leqslant\beta_{j},\forall
j\leqslant k.
\]
Since $\pi_{k}\in\Pi(\mu,\gamma)$ is weakly compact, up to a subsequence
$\pi_{k}$ converges weakly to some $\pi\in\Pi(\mu,\gamma).$ For fixed $j$,
$M>0$, the function $\psi_{j,M}\left(  u,g\right)  =\min\{\left\langle
u-g,e_{j}\right\rangle ^{2},M\}$ is bounded and continuous. We have
\[
\int\psi_{j,M}\left(  u,g\right)  \pi(du,dg)=\lim_{k\rightarrow\infty}\int
\psi_{j,M}\left(  u,g\right)  \pi_{k}(du,dg)\leqslant\beta_{j}.
\]
Sending $M\rightarrow\infty,$ the monotone convergence theorem therefore
yields $\left(  \ref{eq:hilbert-coordinate}\right)  .$

\textbf{4}. We now turn the coordinate bounds into the desired error bounds.
For every $k,$ both marginals of $\pi_{k}$ are supported in $H_{0}$ by the
construction and the already proved fact that $W_{n},G\in H_{0}.$ Hence both
marginals of $\pi$ are supported on $H_{0}$. Thus
\begin{equation}
D\triangleq W_{n}-G\in H_{0},\text{ }\pi\text{-almost surely}.
\label{eq:DH0-almost}
\end{equation}
In particular
\begin{equation}
P_{I}D=D,\text{ }\pi\text{-almost surely}. \label{eq:projection-D}
\end{equation}
Hence by Parseval identity (e.g. \cite[Theorem 3.6-3]
{kreyszig1991introductory}), $\pi$-almost surely,
\begin{equation}
\left\vert D\right\vert ^{2}=\left\vert P_{I}D\right\vert ^{2}=\sum_{j\in
I}\langle P_{I}D,e_{j}\rangle^{2}=\sum_{j\in I}\langle D,e_{j}\rangle^{2}.
\label{eq:square-D}
\end{equation}
Since $\Sigma$ is trace class,
\[
\mathbb{E}_{\pi}\left\vert D\right\vert ^{2}\leqslant2\mathbb{E}\left\vert
W_{n}\right\vert ^{2}+2\mathbb{E}\left\vert G\right\vert ^{2}
=4\operatorname{tr}\Sigma<\infty,
\]
$D$ is square integrable under $\pi.$ Fix any index set $J\subset I.$ Using
the spectral representation of $\Sigma^{1/2}$ and the just proved coordinate
bound $\left(  \ref{eq:hilbert-coordinate}\right)  $, we get
\begin{align*}
\mathbb{E}_{\pi}\left\vert \Sigma^{1/2}P_{J}\left(  W_{n}-G\right)
\right\vert ^{2}  &  =\sum_{j\in J}\lambda_{j}\mathbb{E}_{\pi}\left\langle
W_{n}-G,e_{j}\right\rangle ^{2}\leqslant C_{0}\sum_{j\in J}\lambda_{j}b_{j}\\
&  =\frac{C_{0}}{n^{2}}\sum_{i}\mathbb{E}\left[  |\Delta_{i}|^{2}\sum_{j\in
J}\langle\Delta_{i},e_{j}\rangle^{2}\right]  =\frac{C_{0}}{n^{2}}\sum
_{i}\mathbb{E}\left[  |\Delta_{i}|^{2}|P_{J}\Delta_{i}|^{2}\right]  .
\end{align*}
The last equality uses Parseval identity. Since $|P_{J}\Delta_{i}|^{2}$ is
bounded by $|\Delta_{i}|^{2},$ the above expectation on the right side is
finite and thus $\left(  \ref{eq:hilbert-local}\right)  $ holds.

\textbf{5}. For any centered $H$-valued feature distribution $\xi_{V}$ with
$\int|v|^{2}\xi_{V}(dv)<\infty,$ extend $\pi$ by the independent product
$\xi_{V}(dv)$. Since $\xi_{V}$ is centered, its covariance $A$ has
$\int\langle v,z\rangle^{2}\xi_{V}(dv)=\langle Az,z\rangle$ for fixed $z$.
Using Tonelli theorem and that $A\preccurlyeq\Sigma$, we obtain
\[
\mathbb{E}\left\langle V,D\right\rangle ^{2}=\int\left(  \int\langle
v,w-g\rangle^{2}\xi_{V}(dv)\right)  \pi(dw,dg)=\mathbb{E}_{\pi}\left\langle
AD,D\right\rangle \leqslant\mathbb{E}_{\pi}\left\vert \Sigma^{1/2}D\right\vert
^{2}\leqslant\frac{C_{0}}{n^{2}}\sum_{i}\mathbb{E}\left\vert \Delta
_{i}\right\vert ^{4}.
\]
The last inequality is justified. Indeed, using $\left(  \ref{eq:projection-D}
\right)  $ we have
\begin{align*}
\mathbb{E}_{\pi}|\Sigma^{1/2}D|^{2}  &  =\mathbb{E}_{\pi}|\Sigma^{1/2}
P_{I}D|^{2}\leqslant\frac{C_{0}}{n^{2}}\sum_{i}\mathbb{E}[|\Delta_{i}
|^{2}|P_{I}\Delta_{i}|^{2}]\\
&  =\frac{C_{0}}{n^{2}}\sum_{i}\mathbb{E}[|\Delta_{i}|^{2}|\Delta_{i}
|^{2}]=\frac{C_{0}}{n^{2}}\sum_{i}\mathbb{E}|\Delta_{i}|^{4}.
\end{align*}

\end{proof}

\section{Sharp prediction error for Hilbert-valued i.i.d. observations}

\label{sec:hilbert-example}

Averaging in various form is an important regularization method in modern
machine learning
\cite{mandt2017stochastic,maddox2019simple,hsing2015theoretical}. To
demonstrate an application of our coordiante-wise central limit theorems in
obtaining prediction error bound, we use the Hilbert space example in
\cite[Example~4.1 and Remark~2.2(b)]{fang2026sharp}, originally from
\cite{gotze2011estimates}, with all bounds worked out and compared at the end.

\subsection{Problem setup}

We first specify the observations and the independent future input. Let
$H=\ell^{2}$, the space of real sequences $x=(x_{j})_{j\geqslant1}$ with
$\sum_{j}x_{j}^{2}<\infty$, equipped with the usual inner product $\langle
x,y\rangle=\sum_{j}x_{j}y_{j}$. Let $e_{j}$ be its $j$-th unit coordinate
vector. Let
\[
p_{j}=\frac{c}{j\left[  \log(j+1)\right]  ^{2}}\text{ with }c=\left(
\sum_{j=1}^{\infty}\frac{1}{j\left[  \log(j+1)\right]  ^{2}}\right)  ^{-1}
\]
and
\[
\mathbb{P}(X=e_{j})=\mathbb{P}(X=-e_{j})=\frac{p_{j}}{2}.
\]
The series in $c$ is convergent, thus $0<c<\infty$ and $\sum_{j}p_{j}=1.$ Here
and below $\log$ is the natural logarithm. The distribution remains fixed as
the sample size increases.

Since
\[
\sum_{j=1}^{\infty}\left(  \left\vert \frac{p_{j}}{2}e_{j}\right\vert
+\left\vert \frac{p_{j}}{2}(-e_{j})\right\vert \right)  =\sum_{j}p_{j}=1,
\]
the expectation series converges absolutely in the Hilbert norm, thus grouping
is legitimate which yields
\[
\mathbb{E}X=\sum_{j=1}^{\infty}\left(  \frac{p_{j}}{2}e_{j}+\frac{p_{j}}
{2}(-e_{j})\right)  =0.
\]
By the definition, we have $\left\vert X\right\vert =1\ $almost surely and
\[
\mathbb{E}[\langle X,e_{j}\rangle\langle X,e_{k}\rangle]=p_{j}1_{\{j=k\}}.
\]
It follows that the covariance
\[
\Sigma\triangleq\operatorname{Cov}(X)
\]
has eigen vectors $e_{j}$ and eigen values $p_{j},$
\begin{equation}
\Sigma e_{j}=p_{j}e_{j}, \label{eq:sigma-eigen}
\end{equation}
thus $\operatorname{tr}\Sigma=1.$ Moreover $\mathbb{E}\left\vert X\right\vert
^{4}=1$. In particular all assumptions of Theorem~\ref{thm:hilbert-prediction} hold.

We keep the notations and variables as defined in the proof of
Theorem~\ref{thm:hilbert-prediction}. For independent standard normals $Z_{j}
$, the Gaussian series $\sum_{j}\sqrt{p_{j}}Z_{j}e_{j}$ is $H$-valued by the
construction in that proof.

Throughout this example, $X_{1},\ldots,X_{n}$ are independent copies of $X$,
so these are i.i.d. random elements. Write $W_{n}=n^{-1/2}\sum_{i}X_{i}$. For
its coordinates under the eigenbasis $\left(  \ref{eq:sigma-eigen}\right)  ,$
write $W_{n,j}=\langle W_{n},e_{j}\rangle$, $G_{j}=\langle G,e_{j}\rangle$.
Write
\[
W_{n}\sim\mu\triangleq\mathcal{L}(W_{n}),\text{ }G\sim\gamma\triangleq
N_{H}(0,\Sigma).
\]
Define
\[
D=W_{n}-G.
\]

For any coupling $\pi$ of $\mu$ and $\gamma$, we draw one further independent
copy $X_{\ast}$ of $X$, independently of the whole pair $\left(
W_{n},G\right)  $.\ The The original and Gaussian predictions are respectively
$\langle X_{\ast},W_{n}\rangle$ and $\langle X_{\ast},G\rangle$. Since
$X_{\ast}$ selects coordinate $j$ with probability $p_{j}$, their mean squared
difference is
\begin{equation}
\mathbb{E}\langle X_{\ast},D\rangle^{2}=\mathbb{E}_{\pi}\sum_{j=1}^{\infty
}p_{j}\langle D,e_{j}\rangle^{2}=\mathbb{E}_{\pi}|\Sigma^{1/2}D|^{2}.
\label{eq:prediction-error}
\end{equation}
Note that the first expectation is over the joint distribution of the three
random elements $\left(  X_{\ast},W_{n},G\right)  $ . Taking infimum over all
couplings $\pi$ in $\left(  \ref{eq:prediction-error}\right)  $ gives
$\mathcal{T}_{\Sigma}^{2}\left(  \mu,\gamma\right)  $ by the definition of
$\left(  \ref{eq:transport-identity}\right)  $. All our prediction comparisons
below refer to this quantity.

\subsection{One coupling for the full prediction and every spectral cutoff}

We now compute the coefficient $b_{j}$ in Theorem \ref{thm:hilbert-prediction}
. Let $\left(  X_{i}^{\prime}\right)  _{i}$ be an independent copy of $\left(
X_{i}\right)  _{i}$. Note in this example the two families $X_{i}^{\prime},$
$X_{i}$ are not only independent both also have identical disrtibution. Recall
that $\Delta_{i}=X_{i}^{\prime}-X_{i}.$

Let $s_{2}=\sum_{j}p_{j}^{2}$, so $0<s_{2}\leqslant1$. For each $j$, the event
that both copies $X_{i}^{\prime},$ $X_{i}$ select index $j$ with opposite
signs has probability $p_{j}^{2}/2$. On this event $|\Delta_{i}|^{2}
\langle\Delta_{i},e_{j}\rangle^{2}=16$. The event that exactly one copy
selects $j$ has probability $2p_{j}(1-p_{j})$. On it the same product is $2$.
All other events have zero probability. Therefore
\[
q_{j}\triangleq\mathbb{E}[|\Delta_{i}|^{2}\langle\Delta_{i},e_{j}\rangle
^{2}]=16\frac{p_{j}^{2}}{2}+2\cdot2p_{j}(1-p_{j})=4p_{j}+4p_{j}^{2}.
\]
Hence
\[
b_{j}=\frac{1}{n^{2}\lambda_{j}}\sum_{i}q_{j}=\frac{q_{j}}{np_{j}}
=\frac{4(1+p_{j})}{n},\text{ }\beta_{j}=\min\left\{  2p_{j},\frac
{4C_{0}(1+p_{j})}{n}\right\}  .
\]
Using Parseval identity and Tonelli theorem, we have
\[
\mathbb{E}|\Delta_{i}|^{4}=\mathbb{E}\left[  |\Delta_{i}|^{2}\sum
_{j=1}^{\infty}\left\langle \Delta_{i},e_{j}\right\rangle ^{2}\right]
=\sum_{j=1}^{\infty}\mathbb{E}\left[  |\Delta_{i}|^{2}\left\langle \Delta
_{i},e_{j}\right\rangle ^{2}\right]  =\sum_{j=1}^{\infty}q_{j}=4+4s_{2}
\leqslant8,
\]
We now invoke Theorem~\ref{thm:hilbert-prediction}, there exists a coupling
$\pi_{0}$ of $\left(  W_{n},G\right)  $ with
\[
\mathbb{E}_{\pi_{0}}\langle D,e_{j}\rangle^{2}\leqslant\beta_{j}\text{ for
every }j.
\]
For a deterministic $J\subseteq\mathbb{N}$, let $V_{J}=P_{J}X_{\ast}$. This
$V_{J}$ is centered, square integrable, independent of the pair $\left(
W_{n},G\right)  $, and has covariance $P_{J}\Sigma P_{J}\preccurlyeq\Sigma$.
Proceeding as in $\left(  \ref{eq:prediction-error}\right)  $ and using
$\left(  \ref{eq:hilbert-coordinate}\right)  $ we get
\[
\mathbb{E}\langle V_{J},D\rangle^{2}=\sum_{j\in J}p_{j}\mathbb{E}_{\pi_{0}
}\langle D,e_{j}\rangle^{2}\leqslant\sum_{j\in J}p_{j}\beta_{j}\leqslant
\frac{4C_{0}}{n}\sum_{j\in J}p_{j}(1+p_{j}).
\]
The bound is also a direct result of $\left(  \ref{eq:hilbert-risk}\right)  $
and the bound of $\mathbb{E}|\Delta_{i}|^{4}$. The bound is zero if
$J=\varnothing$. Under $\pi_{0},$ $W_{n}-G$ belongs to $H_{0}$ almost surely
by $\left(  \ref{eq:DH0-almost}\right)  $. Taking $J=\mathbb{N}$ in the above
inequality therefore yields$,$
\begin{equation}
\mathcal{T}_{\Sigma}^{2}(\mu,\gamma)\leqslant\mathbb{E}\langle X_{\ast}
,W_{n}-G\rangle^{2}=\mathbb{E}\langle V_{I},W_{n}-G\rangle^{2}\leqslant
\frac{4C_{0}(1+s_{2})}{n}\leqslant\frac{8C_{0}}{n}. \label{eq:hilbert-fast}
\end{equation}
The first inequality uses $\left(  \ref{eq:conditional-prediction-error}
\right)  $ together with the contraction property of conditional expectation.
Taking $J=\{1,\ldots,k\}$ bounds the retained prediction, while $J_{k}
=\{j:j>k\}$ bounds the omitted part of its error. For $k\geqslant2$,
monotonicity of $x\mapsto1/(x\log^{2}x)$ gives
\[
\sum_{j>k}p_{j}\leqslant c\sum_{j>k}\frac{1}{j\log^{2}j}\leqslant c\int
_{k}^{\infty}\frac{dx}{x\log^{2}x}=\frac{c}{\log k},
\]
\[
\mathbb{E}\langle P_{J_{k}}X_{\ast},D\rangle^{2}\leqslant\frac{4C_{0}}{n}
\sum_{j>k}p_{j}(1+p_{j})\leqslant\frac{8C_{0}c}{n\log k}.
\]
The same joint law of $(W_{n},G,X_{\ast})$ gives these inequalities for every
deterministic $J$ and every cutoff $k$.

\subsection{Sharpness for prediction error bound}

We claim that the prediction order in $\left(  \ref{eq:hilbert-fast}\right)  $
has a lower bound that is independent of any coupling of $(W_{n},G),$ thus
prediction order cannot be improved. Each $\langle X_{i},e_{j}\rangle$ belongs
to $\{-1,0,1\}$, so $W_{n,j}\in h_{n}\mathbb{Z}$, where $h_{n}=n^{-1/2}$. For
$S\subset\mathbb{R}$, write $\operatorname{dist}(x,S)=\inf_{s\in S}|x-s|$. For
every integer $l$,
\[
\int_{lh_{n}}^{(l+1)h_{n}}\operatorname{dist}(x,h_{n}\mathbb{Z})^{2}
dx=2\int_{0}^{h_{n}/2}u^{2}du=\frac{h_{n}^{3}}{12}.
\]
Assume $np_{1}\geqslant4$. The density of $G_{1}\sim N(0,p_{1})$ on
$[0,\sqrt{p_{1}}]$ is at least $e^{-1/2}/\sqrt{2\pi p_{1}}$. This interval
contains $\lfloor\sqrt{p_{1}}/h_{n}\rfloor\geqslant\sqrt{p_{1}}/(2h_{n})$
whole lattice intervals. Since it is always true irrespective of any coupling
of $(W_{n},G)$,
\[
\left\vert W_{n,1}-G_{1}\right\vert \geqslant\operatorname{dist}(G_{1}
,h_{n}\mathbb{Z}),
\]
we obtain
\[
\mathbb{E}_{\pi}|\Sigma^{1/2}D|^{2}\geqslant p_{1}\mathbb{E}
\operatorname{dist}(G_{1},h_{n}\mathbb{Z})^{2}\geqslant p_{1}\frac{e^{-1/2}
}{\sqrt{2\pi p_{1}}}\frac{\sqrt{p_{1}}}{2h_{n}}\frac{h_{n}^{3}}{12}
=\frac{p_{1}e^{-1/2}}{24\sqrt{2\pi}n}.
\]
The first inequality uses $\left(  \ref{eq:prediction-error}\right)  $. Taking
the infimum and combining with $\left(  \ref{eq:hilbert-fast}\right)  $ gives
\begin{equation}
\frac{p_{1}e^{-1/2}}{24\sqrt{2\pi}n}\leqslant\mathcal{T}_{\Sigma}^{2}
(\mu,\gamma)\leqslant\frac{4C_{0}(1+s_{2})}{n},\text{ }n\geqslant\frac
{4}{p_{1}}. \label{eq:sharp-root}
\end{equation}
Thus the optimal root prediction error is bounded above and below by positive
constant multiples of $n^{-1/2}$.

\subsection{The slow ordinary Hilbert transport error}

The preceding sharp prediction bound, which is weighted by $\Sigma^{1/2},$
does not control the sum of all unweighted coordinate errors at the same rate.
Here the small variance coordinates give a lower bound on that sum under every coupling.

Let $c_{1}=\mathbb{E}[Z^{2}1_{\{|Z|\leqslant1\}}]>0$ for a standard normal
$Z$. If $p_{j}\leqslant1/(4n)$, then $\sqrt{p_{j}}\leqslant h_{n}/2$, and the
nearest lattice point to $G_{j}$ is zero on $\left\vert G_{j}\right\vert
\leqslant\sqrt{p_{j}}$. Since $G_{j}\sim\sqrt{p_{j}}Z$, any coupling $\pi$ of
$(W_{n},G)$ satisfies
\[
\mathbb{E}_{\pi}(W_{n,j}-G_{j})^{2}\geqslant\mathbb{E}\operatorname{dist}
(G_{j},h_{n}\mathbb{Z})^{2}\geqslant\mathbb{E}[G_{j}^{2}1_{\{|G_{j}
|\leqslant\sqrt{p_{j}}\}}]=p_{j}\mathbb{E}[Z^{2}1_{\{|Z|\leqslant1\}}
]=c_{1}p_{j}.
\]
Let
\[
j_{n}=\max\left\{  2,\left\lceil 4cn/(\log2)^{2}\right\rceil \right\}  .
\]
For $j\geqslant j_{n}$,
\[
p_{j}\leqslant\frac{c}{j(\log2)^{2}}\leqslant\frac{1}{4n}.
\]
For $j\geqslant2$, $j+1\leqslant j^{2}$, so $\log(j+1)\leqslant2\log j$. The
function $x\mapsto1/(x\log^{2}x)$ decreases on $[2,\infty)$. It follows that
\[
\sum_{j\geqslant j_{n}}p_{j}\geqslant\frac{c}{4}\sum_{j\geqslant j_{n}}
\frac{1}{j\log^{2}j}\geqslant\frac{c}{4}\int_{j_{n}}^{\infty}\frac{dx}
{x\log^{2}x}=\frac{c}{4\log j_{n}}.
\]
As in $\left(  \ref{eq:square-D}\right)  ,$ by Parseval identity and and
Tonelli theorem we have, for any coupling $\pi$ of $(W_{n},G),$
\[
\mathbb{E}_{\pi}\left\vert W_{n}-G\right\vert ^{2}=\sum_{j\in I}
\mathbb{E}_{\pi}(W_{n,j}-G_{j})^{2}\geqslant\sum_{j\geqslant j_{n}}
\mathbb{E}_{\pi}(W_{n,j}-G_{j})^{2}\geqslant\frac{c_{1}c}{4\log j_{n}}.
\]
Hence
\begin{equation}
W_{2}^{2}(\mathcal{L}(W_{n}),N_{H}(0,\Sigma))\geqslant\frac{c_{1}c}{4\log
j_{n}}. \label{eq:hilbert-slow}
\end{equation}
Let $K_{0}=3+4c/(\log2)^{2}$. By the definition, $j_{n}\leqslant K_{0}n$ for
$n\geqslant1$. For $n\geqslant2$,
\[
\log j_{n}\leqslant\log K_{0}+\log n\leqslant\left(  1+\frac{\log K_{0}}
{\log2}\right)  \log n.
\]
Thus $\left(  \ref{eq:hilbert-slow}\right)  $ is at least $c_{2}/\log n$, with
$c_{2}=c_{1}c/[4(1+\log K_{0}/\log2)]>0$.

On the other hand, for any fixed $J\subset I,$ the coordinate bound $\left(
\ref{eq:hilbert-coordinate}\right)  $ gives
\begin{equation}
\mathbb{E}_{\pi}|P_{J}D|^{2}\leqslant\sum_{j\in J}\beta_{j}
\label{eq:original-norm}
\end{equation}
Using $\left(  \ref{eq:square-D}\right)  $ and taking $J=\mathbb{N}$ in
$\left(  \ref{eq:original-norm}\right)  $, we get
\[
W_{2}^{2}(\mu,\gamma)=\mathbb{E}_{\pi}\left\vert D\right\vert ^{2}
=\mathbb{E}_{\pi}\left\vert P_{I}D\right\vert ^{2}\leqslant\sum_{j=1}^{\infty
}\min\left\{  2p_{j},\frac{4C_{0}(1+p_{j})}{n}\right\}  \rightarrow0.
\]
The above limit is guaranteed dominated convergence theorem, since by summand
tends to zero and is bounded by the summable sequence $(2p_{j})_{j}$.

Consequently the ordinary Hilbert distance $W_{2}(\mu,\gamma)$ tends to zero
but is at least $\sqrt{c_{2}/\log n}$ for $n\geqslant2$, whereas the root
prediction error $\mathcal{T}_{\Sigma}(\mu,\gamma)$ has the sharper order
$n^{-1/2}$ by $\left(  \ref{eq:sharp-root}\right)  $.

\subsection{Earlier bounds for the same prediction loss}

Before comparing the above bounds with \cite{fang2026sharp}, we first compute
$\mathcal{T}_{\Sigma}(\mu,\gamma)$. We apply the transport identity $\left(
\ref{eq:transport-identity}\right)  $ with $\left(
\ref{eq:conditional-prediction-error}\right)  $ and compute each of its two
transformed marginals. Write $T=\Sigma^{1/2}$, $Y=TX$, and $Y_{i}=TX_{i}$.
Because $\Sigma e_{j}=p_{j}e_{j}$ by $\left(  \ref{eq:sigma-eigen}\right)  $,
the spectral representation of $T$ gives
\[
T(\pm e_{j})=\pm\sqrt{p_{j}}e_{j}.
\]
Thus the explicit pushforward law $\mathcal{L}(Y)=T_{\#}\mathcal{L}(X)$ is
specified by
\[
\mathbb{P}(Y=\sqrt{p_{j}}e_{j})=\mathbb{P}(X=e_{j})=\frac{p_{j}}{2},\text{
}\mathbb{P}(Y=-\sqrt{p_{j}}e_{j})=\mathbb{P}(X=-e_{j})=\frac{p_{j}}{2}.
\]
Since $T$ is linear, we also have
\[
TW_{n}=\frac{1}{\sqrt{n}}\sum_{i=1}^{n}TX_{i}=\frac{1}{\sqrt{n}}\sum_{i=1}
^{n}Y_{i}\text{ and }T_{\#}\mu=\mathcal{L}\left(  n^{-1/2}\sum_{i}
Y_{i}\right)  .
\]
The $Y_{i}$ remain independent because they are deterministic functions of the
independent $X_{i}$.

Since $T$ is a bounded linear map, it commutes with the Bochner expectation
\cite[E.11]{cohn2013measure}. So $\mathbb{E}Y=T\mathbb{E}X=0$. This can also
be computed directly by noting the two sign contributions of $Y$ cancel and
the series resulting from the expection converges absolutely, since
$\mathbb{E}\left\vert Y\right\vert =\sum_{j}p_{j}^{3/2}\leqslant1$. The second
and fourth moments are
\[
\mathbb{E}\left\vert Y\right\vert ^{2}=\sum_{j}\left(  \frac{p_{j}}{2}
p_{j}+\frac{p_{j}}{2}p_{j}\right)  =\sum_{j}p_{j}^{2}=s_{2},
\]
\[
\mathbb{E}\left\vert Y\right\vert ^{4}=\sum_{j}\left(  \frac{p_{j}}{2}
p_{j}^{2}+\frac{p_{j}}{2}p_{j}^{2}\right)  =\sum_{j}p_{j}^{3}\leqslant1.
\]
For $K_{Y}=\operatorname{Cov}(Y)=\mathbb{E}\left(  Y\otimes Y\right)  $ and a
fixed $u\in H$, compute its action directly by
\[
K_{Y}u=\sum_{j=1}^{\infty}\left\{  \frac{p_{j}}{2}\langle\sqrt{p_{j}}
e_{j},u\rangle\sqrt{p_{j}}e_{j}+\frac{p_{j}}{2}\langle-\sqrt{p_{j}}
e_{j},u\rangle(-\sqrt{p_{j}}e_{j})\right\}  =\sum_{j=1}^{\infty}p_{j}
^{2}\langle e_{j},u\rangle e_{j}.
\]
The expectation is Bochner integrable because $\mathbb{E}\left\vert \langle
Y,u\rangle Y\right\vert \leqslant\left\vert u\right\vert \mathbb{E}\left\vert
Y\right\vert ^{2}<\infty$. Taking $u=e_{l}$ yields $K_{Y}e_{l}=p_{l}^{2}
e_{l}.$ Thus
\[
\mathbb{E}[\langle Y,e_{j}\rangle\langle Y,e_{l}\rangle]=p_{j}^{2}
1_{\{j=l\}},
\]
and $\operatorname{tr}K_{Y}=s_{2}.$ In particular $K_{Y}=\Sigma^{2}$. The same
identity also follows from the general covariance transformation, for $u,v\in
H$,
\[
\langle K_{Y}u,v\rangle=\mathbb{E}[\langle TX,u\rangle\langle TX,v\rangle
]=\mathbb{E}[\langle X,Tu\rangle\langle X,Tv\rangle]=\langle\Sigma
Tu,Tv\rangle=\langle T\Sigma Tu,v\rangle=\langle\Sigma^{2}u,v\rangle.
\]
By centering and independence, $\operatorname{Cov}(TW_{n})=n^{-1}\sum
_{i}\operatorname{Cov}(Y_{i})=\Sigma^{2}$.

For the other marginal, if $G\sim\gamma=N_{H}(0,\Sigma)$, then $TG$ is
Gaussian with covariance $T\Sigma T=\Sigma^{2}$ and has distribution
\[
\mathcal{L}\left(  TG\right)  =\mathcal{L}\left(  \sum_{j=1}^{\infty}
p_{j}Z_{j}e_{j}\right)  \text{ and }T_{\#}\gamma=N_{H}(0,\Sigma^{2}),
\]
where the $Z_{j}$ are independent standard normals. The series converges in
mean square and almost surely, because $\sum_{j}p_{j}^{2}<\infty$, by the same
orthogonal series argument as in Theorem~\ref{thm:hilbert-prediction}. By the
definition $\left(  \ref{eq:transport-identity}\right)  $,
\[
\mathcal{T}_{\Sigma}(\mu,\gamma)=W_{2}(T_{\#}\mu,T_{\#}\gamma)=W_{2}\left(
\mathcal{L}\left(  n^{-1/2}\sum_{i}Y_{i}\right)  ,N_{H}(0,\Sigma^{2})\right)
.
\]
For an independent feature $X_{\ast}$, the identity $\left(
\ref{eq:conditional-prediction-error}\right)  $ specializes to
\[
\mathbb{E}[\langle X_{\ast},W_{n}-G\rangle^{2}\mid W_{n},G]=\sum_{j}
p_{j}\langle W_{n}-G,e_{j}\rangle^{2}=|T(W_{n}-G)|^{2}.
\]

We first apply the Hilbert CLT \cite[Theorem~2.2]{fang2026sharp} to $Y$. \ Let
$Y^{\prime}$ be an independent copy of $Y$ and let
\[
Q_{Y}=\mathbb{E}\left[  \left\vert Y^{\prime}-Y\right\vert ^{2}\left(
Y^{\prime}-Y\right)  \otimes\left(  Y^{\prime}-Y\right)  \right]  .
\]
The expectation exists since $\left\vert Y\right\vert \leqslant1$. The finite
fourth moment assumption there is also satisfied. To compute the numerator
$\operatorname{tr}(Q_{Y}^{1/2})$ in the cited theorem$,$ it is sufficient to
determine $\left\langle Q_{Y}e_{j},e_{j}\right\rangle $ for each $j.$

The off diagonal entries of $Q_{Y}$ vanish. Indeed, a simultaneous sign change
of the $j$-th coordinate of $Y,Y^{\prime}$ preserves the law and the squared
norm but changes the sign of its $(j,k)$ entry for $k\neq j$. For diagonal
entries, if both $Y,Y^{\prime}$ select coordinate $j$ with opposite signs, the
probability is $2(p_{j}/2)^{2}=p_{j}^{2}/2$, then $Y^{\prime}-Y=\pm
2\sqrt{p_{j}}e_{j}$, thus
\[
\left\vert Y^{\prime}-Y\right\vert ^{2}\left\langle Y^{\prime}-Y,e_{j}
\right\rangle ^{2}=(4p_{j})(4p_{j})=16p_{j}^{2}.
\]
If one of $Y,Y^{\prime}$ selects $j$ and the other selects $k\neq j$, the
probability is $2p_{j}p_{k}$ with $\left\vert Y^{\prime}-Y\right\vert ^{2}$
equal to $p_{j}+p_{k}$ and $\left\langle Y^{\prime}-Y,e_{j}\right\rangle ^{2}$
equal to $p_{j}$, so the product is $p_{j}(p_{j}+p_{k})$. Hence
\begin{align*}
\langle Q_{Y}e_{j},e_{j}\rangle &  =8p_{j}^{4}+2p_{j}^{2}\sum_{k\neq j}
p_{k}(p_{j}+p_{k})=8p_{j}^{4}+2p_{j}^{3}(1-p_{j})+2p_{j}^{2}(s_{2}-p_{j}
^{2})\\
&  =2p_{j}^{3}+2s_{2}p_{j}^{2}+4p_{j}^{4}.
\end{align*}
Therefore
\[
\operatorname{tr}(Q_{Y}^{1/2})=\sum_{j}p_{j}(2p_{j}+2s_{2}+4p_{j}^{2}
)^{1/2}\leqslant\sqrt{8}\sum_{j}p_{j}=\sqrt{8}.
\]

The transformed element $TG$ has finite fourth moment. Recall that
$\mathbb{E}Z_{j}^{4}=3$ and $\mathbb{E}\left(  Z_{j}^{2}Z_{l}^{2}\right)  =1$
for $j\neq l.$ For finite $m$, we get
\[
\mathbb{E}\left(  \sum_{j\leqslant m}p_{j}^{2}Z_{j}^{2}\right)  ^{2}
=3\sum_{j\leqslant m}p_{j}^{4}+2\sum_{j<l\leqslant m}p_{j}^{2}p_{l}
^{2}=\left(  \sum_{j\leqslant m}p_{j}^{2}\right)  ^{2}+2\sum_{j\leqslant
m}p_{j}^{4}.
\]
The expressions inside the expectations increase with $m$, so using monotone
convergence theorem we have
\[
\mathbb{E}\left(  \sum_{j}p_{j}^{2}Z_{j}^{2}\right)  ^{2}=s_{2}^{2}+2\sum
_{j}p_{j}^{4}\leqslant3s_{2}^{2}<\infty,
\]
where $\sum_{j}p_{j}^{4}\leqslant(\sum_{j}p_{j}^{2})^{2}$. Denote $W_{4}$ by
Wasserstein distance w.r.t. the cost $|x-y|^{4}$. Since $W_{2}\leqslant W_{4}
$, we therefore obtain by \cite[Theorem~2.2]{fang2026sharp} that
\[
\mathcal{T}_{\Sigma}(\mu,\gamma)\leqslant W_{4}(\mathcal{L}(TW_{n}
),N_{H}(0,\Sigma^{2}))\leqslant Cn^{-1/4}\left[  \operatorname{tr}(Q_{Y}
^{1/2})\right]  ^{1/2}\leqslant C8^{1/4}n^{-1/4},
\]
where $C$ is the universal constant in that theorem. Since $\operatorname{tr}
(Q_{Y}^{1/2})$ is positive and independent of $n$, this bound gives order
$n^{-1/4}$. However, it is worth noting that if we apply the cited theorem
directly to the original $X,$ instead of the transformed $Y,$ then the upper
bound there becomes infinity. Indeed, define
\[
Q_{X}=\mathbb{E}[|X^{\prime}-X|^{2}\left(  X^{\prime}-X\right)  \otimes\left(
X^{\prime}-X\right)  ].
\]
By the same calculations as above we have $Q_{X}e_{j}=(4p_{j}+4p_{j}^{2}
)e_{j}$. Therefore
\[
\operatorname{tr}(Q_{X}^{1/2})=2\sum_{j}\sqrt{p_{j}}(1+p_{j})^{1/2}
\geqslant2\sum_{j}\sqrt{p_{j}}=\infty.
\]
This shows that direct use of that Hilbert bound \cite[Theorem~2.2]
{fang2026sharp} for the original $X$ yields an bound of infinity.

The aggregate CLT \cite[Theorem~2.1]{fang2026sharp} gives a sharper prediction
than its Hilbert CLT bound. We include the coordinate-wise minima in that
estimate. Let $P_{m}$ project onto $e_{1},\ldots,e_{m}$ and put $s_{2,m}
=\sum_{j\leqslant m}p_{j}^{2}$. For $j\leqslant m$, similar to the
calculations of $\langle Q_{Y}e_{j},e_{j}\rangle,$ we have
\[
\mathbb{E}\left[  \left\vert P_{m}\left(  Y^{\prime}-Y\right)  \right\vert
^{2}\left\langle Y^{\prime}-Y,e_{j}\right\rangle ^{2}\right]  =8p_{j}
^{4}+2p_{j}^{2}\sum_{k\leqslant m,k\neq j}p_{k}(p_{j}+p_{k})+2p_{j}^{3}
\sum_{k>m}p_{k}=2p_{j}^{3}+2s_{2,m}p_{j}^{2}+4p_{j}^{4}.
\]
The last sum accounts for the case where one of $Y,Y^{\prime}$ selects $j$ and
the other selects $k\neq j$ with $k>m$, the probability is $2p_{j}p_{k}$ with
$\left\vert Y^{\prime}-Y\right\vert ^{2}$ and $\left\langle Y^{\prime}
-Y,e_{j}\right\rangle ^{2}$ both equal to $p_{j}$. The covariance of $P_{m}Y$
has positive eigenvalues $p_{j}^{2}$. The aggregate CLT theorem therefore
yields
\[
W_{2}^{2}(\mathcal{L}(P_{m}TW_{n}),\mathcal{L}(P_{m}TG))\leqslant
2\sum_{j\leqslant m}\min\left\{  p_{j}^{2},\frac{2p_{j}+2s_{2,m}+4p_{j}^{2}
}{n}\right\}  .
\]
As $m\rightarrow\infty,$ the right hand side is dominated by the summable
sequence $(2p_{j}^{2})$, so we can use the dominated convergence theorem to
find its limit. The left hand side converges to $W_{2}^{2}(\mathcal{L}
(TW_{n}),\mathcal{L}(TG))$. Indeed, by triangle inequality we have
\begin{align*}
&  \left\vert W_{2}(\mathcal{L}(P_{m}TW_{n}),\mathcal{L}(P_{m}TG))-W_{2}
(\mathcal{L}(TW_{n}),\mathcal{L}(TG))\right\vert \\
&  \leqslant W_{2}(\mathcal{L}(P_{m}TW_{n}),\mathcal{L}(TW_{n}))+W_{2}
(\mathcal{L}(P_{m}TG),\mathcal{L}(TG)).
\end{align*}
Write $U$ for $TW_{n}$ or $TG.$ Both choices of $U$ are centered and have
covariance $\Sigma^{2}$. Projecting $U$ to $P_{m}U$ creates a natural coupling
between $P_{m}U$ and $U.$ Using this natural coupling and Parseval identity,
we have (recalling $\left(  \ref{eq:sigma-eigen}\right)  $),
\[
W_{2}^{2}(\mathcal{L}(P_{m}U),\mathcal{L}(U))\leqslant\mathbb{E}\left\vert
U-P_{m}U\right\vert ^{2}=\mathbb{E}\sum_{j>m}\left\langle U,e_{j}\right\rangle
^{2}=\mathbb{E}\sum_{j>m}\left\langle U,e_{j}\right\rangle ^{2}=\sum
_{j>m}p_{j}^{2}\rightarrow0.
\]
We therefore have that
\[
\mathcal{T}_{\Sigma}^{2}(\mu,\gamma)\leqslant B_{n}\triangleq2\sum
_{j=1}^{\infty}\min\left\{  p_{j}^{2},\frac{2p_{j}+2s_{2}+4p_{j}^{2}}
{n}\right\}  .
\]
We now compute the order of $B_{n}$. For all sufficiently large $n$, let
$k_{n}=\lfloor\sqrt{n}/(\log n)^{2}\rfloor\geqslant2$. The definition formula
for $p_{j}$ implies
\[
\frac{k_{n}(\log n)^{2}}{\sqrt{n}}\rightarrow1,\text{ }\frac{\log(k_{n}
+1)}{\log n}\rightarrow\frac{1}{2}\text{ and }np_{k_{n}}^{2}\rightarrow
16c^{2}.
\]
Indeed, the floor changes $\sqrt{n}/(\log n)^{2}$ by less than one, and its
logarithm is $\tfrac{1}{2}\log n-2\log\log n+o(1)$. Since $(p_{j})$ decreases
and $2s_{2}\leqslant2p_{j}+2s_{2}+4p_{j}^{2}\leqslant8$, we have
\[
B_{n}\geqslant2k_{n}\min\left\{  p_{k_{n}}^{2},\frac{2s_{2}}{n}\right\}  ,
\]
and
\[
B_{n}\leqslant\frac{16k_{n}}{n}+2\sum_{j>k_{n}}p_{j}^{2}\leqslant\frac
{16k_{n}}{n}+\frac{2c^{2}}{k_{n}\{\log(k_{n}+1)\}^{4}}.
\]
For the last inequality, use $p_{j}^{2}\leqslant c^{2}/[j^{2}\{\log
(k_{n}+1)\}^{4}]$ for $j>k_{n}$ and $\sum_{j>k_{n}}j^{-2}\leqslant\int_{k_{n}
}^{\infty}x^{-2}\,dx$. The $\min\{np_{k_{n}}^{2},2s_{2}\}$ in the lower bound
of $B_{n}$ tends to the positive number $\min\{16c^{2},2s_{2}\}$. The terms in
the upper bound of $B_{n}$ each have order $n^{-1/2}(\log n)^{-2}$. Thus there
are constants $0<a\leqslant b<\infty$, independent of $n$, such that
\[
\frac{a}{\sqrt{n}(\log n)^{2}}\leqslant B_{n}\leqslant\frac{b}{\sqrt{n}(\log
n)^{2}}
\]
for all sufficiently large $n$. In particular $\sqrt{B_{n}}$ has order
$n^{-1/4}/\log n$.

We therefore have the improvements,
\[
\begin{array}
[c]{rl}
\text{Theorem~\ref{thm:hilbert-prediction}:} & n^{-1/2}\text{ (matching lower
bound).}\\
\text{aggregate CLT \cite[Theorem~2.1]{fang2026sharp}:} & n^{-1/4}/\log n,\\
\text{Hilbert CLT \cite[Theorem~2.2]{fang2026sharp}:} & n^{-1/4}.
\end{array}
\]

\bibliographystyle{plain}
\bibliography{main}

@article{fang2026sharp,
  title={Sharp Fourth-Moment Wasserstein CLT Bounds},
  author={Fang, Xiao and Su, Zi-Yao and Wang, Jiaxue and Wang, Xiaolin},
  journal={arXiv preprint arXiv:2610.35354},
  year={2026}
}

@article{bonis2020stein,
  title={Stein’s method for normal approximation in Wasserstein distances with application to the multivariate central limit theorem},
  author={Bonis, Thomas},
  journal={Probability Theory and Related Fields},
  volume={178},
  number={3},
  pages={827--860},
  year={2020},
  publisher={Springer}
}

@article{shreve1991brownian,
  title={Brownian motion and stochastic calculus},
  author={Shreve, Steven E and Karatzas, Ioannis},
  journal={Graduate texts in mathematics},
  volume={113},
  year={1991}
}

@book{revuz2013continuous,
  title={Continuous martingales and Brownian motion},
  author={Revuz, Daniel and Yor, Marc},
  year={2013},
  publisher={Springer Science \& Business Media}
}

@book{kallenberg2021foundations,
  title={Foundations of modern probability Third Edition},
  author={Kallenberg, Olav},
  year={2021},
  publisher={Springer}
}

@article{el2022information,
  title={An information-theoretic view of stochastic localization},
  author={El Alaoui, Ahmed and Montanari, Andrea},
  journal={IEEE Transactions on Information Theory},
  volume={68},
  number={11},
  pages={7423--7426},
  year={2022},
  publisher={IEEE}
}

@article{eldan2020clt,
  title={The CLT in high dimensions},
  author={Eldan, Ronen and Mikulincer, Dan and Zhai, Alex},
  journal={The Annals of Probability},
  volume={48},
  number={5},
  pages={2494--2524},
  year={2020},
  publisher={JSTOR}
}

@inproceedings{rio2009upper,
  title={Upper bounds for minimal distances in the central limit theorem},
  author={Rio, Emmanuel},
  booktitle={Annales de l'IHP Probabilit{\'e}s et statistiques},
  volume={45},
  number={3},
  pages={802--817},
  year={2009}
}

@article{eldan2016skorokhod,
  title={Skorokhod embeddings via stochastic flows on the space of Gaussian measures},
  author={Eldan, Ronen},
  journal={Annales de l'Institut Henri Poincaré, Probabilités et Statistiques},
  volume={52},
  pages={1259-1280},
  year={2016}
}

@book{rockafellar1997convex,
  title={Convex analysis},
  author={Rockafellar, R Tyrrell},
  volume={28},
  year={1997},
  publisher={Princeton university press}
}

@article{comon2008symmetric,
  title={Symmetric tensors and symmetric tensor rank},
  author={Comon, Pierre and Golub, Gene and Lim, Lek-Heng and Mourrain, Bernard},
  journal={SIAM Journal on Matrix Analysis and Applications},
  volume={30},
  number={3},
  pages={1254--1279},
  year={2008},
  publisher={SIAM}
}

@book{durrett2019probability,
  title={Probability: theory and examples},
  author={Durrett, Rick},
  volume={49},
  year={2019},
  publisher={Cambridge university press}
}

@article{Wihler2009Holder ,
  author  = {Wihler, Thomas P.},
  title   = {On the {H}\"older continuity of matrix functions for normal matrices},
  journal = {Journal of Inequalities in Pure and Applied Mathematics},
  volume  = {10},
  number  = {4},
  year    = {2009}
}

@article{zhai2018high,
  title={A high-dimensional CLT in W 2 distance with near optimal convergence rate},
  author={Zhai, Alex},
  journal={Probability Theory and Related Fields},
  volume={170},
  number={3},
  pages={821--845},
  year={2018},
  publisher={Springer}
}

@article{eldan2013thin,
  title={Thin shell implies spectral gap up to polylog via a stochastic localization scheme},
  author={Eldan, Ronen},
  journal={Geometric and Functional Analysis},
  volume={23},
  number={2},
  pages={532--569},
  year={2013},
  publisher={Springer}
}

@article{fang2021high,
  title={High-dimensional central limit theorems by Stein’s method},
  author={Fang, Xiao and Koike, Yuta},
  journal={The Annals of Applied Probability},
  volume={31},
  number={4},
  pages={1660--1686},
  year={2021},
  publisher={JSTOR}
}

@article{lopes2022central,
  title={Central limit theorem and bootstrap approximation in high dimensions: Near 1/n rates via implicit smoothing},
  author={Lopes, Miles E},
  journal={The Annals of Statistics},
  volume={50},
  number={5},
  pages={2492--2513},
  year={2022},
  publisher={Institute of Mathematical Statistics}
}

@article{bentkus2003dependence,
  title={On the dependence of the Berry--Esseen bound on dimension},
  author={Bentkus, Vidmantas},
  journal={Journal of Statistical Planning and Inference},
  volume={113},
  number={2},
  pages={385--402},
  year={2003},
  publisher={Elsevier}
}

@article{bentkus2005lyapunov,
  title={A Lyapunov-type bound in Rd},
  author={Bentkus, Vidmantas},
  journal={Theory of Probability \& Its Applications},
  volume={49},
  number={2},
  pages={311--323},
  year={2005},
  publisher={SIAM}
}

@article{esseen1956moment,
  title={A moment inequality with an application to the central limit theorem},
  author={Esseen, Carl-Gustav},
  journal={Scandinavian Actuarial Journal},
  volume={1956},
  number={2},
  pages={160--170},
  year={1956},
  publisher={Taylor \& Francis}
}

@book{kreyszig1991introductory,
  title={Introductory functional analysis with applications},
  author={Kreyszig, Erwin},
  year={1991},
  publisher={John Wiley \& Sons}
}

@book{hsing2015theoretical,
  title={Theoretical foundations of functional data analysis, with an introduction to linear operators},
  author={Hsing, Tailen and Eubank, Randall},
  year={2015},
  publisher={John Wiley \& Sons}
}

@article{gotze2011estimates,
  title={Estimates for the rate of strong approximation in Hilbert space},
  author={G{\"o}tze, Friedrich and Zaitsev, A Yu},
  journal={Siberian mathematical journal},
  volume={52},
  number={4},
  pages={628--638},
  year={2011},
  publisher={Springer}
}

@article{maddox2019simple,
  title={A simple baseline for bayesian uncertainty in deep learning},
  author={Maddox, Wesley J and Izmailov, Pavel and Garipov, Timur and Vetrov, Dmitry P and Wilson, Andrew Gordon},
  journal={Advances in neural information processing systems},
  volume={32},
  year={2019}
}

@article{mandt2017stochastic,
  title={Stochastic gradient descent as approximate bayesian inference},
  author={Mandt, Stephan and Hoffman, Matthew D and Blei, David M},
  journal={Journal of Machine Learning Research},
  volume={18},
  number={134},
  pages={1--35},
  year={2017}
}

@misc{hastie2009elements,
  title={The elements of statistical learning},
  author={Hastie, Trevor and Tibshirani, Robert and Friedman, Jerome and others},
  year={2009},
  publisher={Springer series in statistics, New-York}
}

@book{reed1980methods,
  title={Methods of modern mathematical physics: Functional analysis},
  author={Reed, Michael and Simon, Barry},
  volume={1},
  year={1980},
  publisher={Gulf Professional Publishing}
}

@book{cohn2013measure,
  title={Measure theory},
  author={Cohn, Donald L},
  volume={2},
  year={2013},
  publisher={Springer}
}
\end{document}